\documentclass[11pt]{amsart}
\usepackage{setspace}
\usepackage{amsmath, amsfonts, amsthm, amstext, amssymb, xspace, float}
\usepackage{eucal}
\usepackage{amssymb,latexsym, mathtools}
\usepackage{xcolor}
\usepackage[shortlabels]{enumitem}
\usepackage{microtype}
\usepackage{longtable,booktabs}
\usepackage{tikz-cd}
\usepackage{graphicx}
\usepackage{hyperref}
\usepackage[capitalise, nameinlink]{cleveref}
\hypersetup{
    colorlinks,
    linkcolor={red!50!black},
    citecolor={blue!50!black},
    urlcolor={blue!80!black}
}

\theoremstyle{definition}
\newtheorem{thm}{Theorem}[section]
\newtheorem{lemma}[thm]{Lemma}
\newtheorem{corollary}[thm]{Corollary}
\newtheorem{proposition}[thm]{Proposition}

\newtheorem{remark}[thm]{Remark}

\usepackage{todonotes}

\def\c{\mathbb{C}}

\def\f{\mathbb{F}}

\def\m{\mathbb{M}}
\def\n{\mathbb{N}}

\def\q{\mathbb{Q}}
\def\r{\mathbb{R}}
\def\z{\mathbb{Z}}

\def\ca{\mathcal{A}}
\def\ce{\mathcal{E}}

\def\Sp{\mathcal{S}\mathrm{p}}

\def\Spec{\operatorname{Spec}}
\def\Ext{\operatorname{Ext}}

\def\ker{\operatorname{ker}}

\def\Sq{\operatorname{Sq}}

\author[P. Dhankhar]{Prerna Dhankhar}\address{University of Kentucky}\email{prerna.dhankhar@uky.edu}
\author[R. Field]{Rebecca Field}\address{James Madison University}\email{fieldre@jmu.edu}
\author[A. Nigam]{Arjun Nigam}\address{Duke University}\email{arjun.nigam@duke.edu}
\author[J. D. Quigley]{J.D. Quigley}\address{University of Virginia}\email{mbp6pj@virginia.edu}
\author[A. J. Yang]{Albert Jinghui Yang}\address{University of Pennsylvania}\email{yangjh@sas.upenn.edu}

\title[Very effective algebraic and Hermitian K-theory of $C_2$]{Very effective algebraic and Hermitian K-theory of the cyclic group of order two}

\begin{document}

\begin{abstract}
We compute the $2$-completed effective algebraic K-theory and very effective Hermitian K-theory of the geometric classifying space of the cyclic group of order two over algebraically closed fields, real closed fields, and finite fields. \end{abstract}

\maketitle

\tableofcontents

\section{Introduction}

The Chow groups of classifying spaces of finite groups arise in counterexamples to the integral Hodge conjecture \cite{AH62}. These Chow groups were studied systematically by Totaro in \cite{Tot14}, where he suggested that motivic classifying spaces should also be studied using richer motivic cohomology theories. The Chow--Witt groups of various classifying spaces have been the subject of much recent attention, with applications in enumerative geometry and the classification of $G$-torsors, cf. \cite{BW25, DLM23, Lac25, HW19, Wen24}. 

In this work, we continue the study of motivic classifying spaces using several richer invariants. We will always work with  \emph{geometric classifying spaces} as introduced by Morel and Voevodsky in \cite[Sec.~4.2]{MV99}. Recall that effective algebraic K-theory $kgl$ is the motivic analogue of connective complex topological K-theory $ku$ and very effective Hermitian K-theory $kq$ is the analogue of connective real topological K-theory $ko$. We will compute the $kgl$- and $kq$-homology of the classifying space of the cyclic group of order two after $2$-completion over different base fields:

\begin{thm}\label{MT}
After $2$-completion, we compute the effective algebraic K-theory and very effective Hermitian K-theory of $BC_2$, as modules over the coefficients of effective algebraic K-theory and very effective Hermitian K-theory, respectively, over the following base fields:
\begin{enumerate}
\item algebraically closed fields (\cref{Sec:C});
\item real closed fields (\cref{Sec:R});
\item finite fields of odd prime order (\cref{Sec:Fq}).
\end{enumerate}
\end{thm} 

Inspecting our computations in the motivic setting, we see that there are natural surjective ring homomorphisms from the motivic homology to the  classical groups $H\z_*BC_2$ 
over algebraically closed fields, over real closed fields, and over finite fields via maps 
$$H\z^{\operatorname{mot}}_{i,\lceil i/2 \rceil}(BC_2) \twoheadrightarrow H\z_i(BC_2)$$
that send the new generators in degree $0,1$ to the identity and annihilate any other new generators. 

Moreover, the corresponding maps from $H\z^{\operatorname{mot}}_{**}(BC_2)$  over real closed fields and over finite fields to the analogous groups over algebraically closed fields commute with these homomorphisms. The same statements hold if we replace $H\z^{\operatorname{mot}}$ and $H\z$ by $kgl$ and $ku$, respectively, or by $kq$ and $ko$, respectively. 

\begin{corollary}
If $F$ is any field of positive odd characteristic, or any field of characteristic zero which admits a real embedding, then we have natural surjections 
\begin{align*}
    & H\z_{i,\lceil i/2 \rceil}^{mot}(BC_2) \twoheadrightarrow H\z_i(BC_2), \\
    & kgl_{i,\lceil i/2 \rceil} (BC_2) \twoheadrightarrow ku_i(BC_2), \\
    & kq_{i,\lceil i/2 \rceil}(BC_2)\twoheadrightarrow ko_i(BC_2),
\end{align*}
for all $i \in \z$. 
\end{corollary}

Although they have many similarities, the groups $H\z_{**}(BC_2)$, $kgl_{**}(BC_2)$, and $kq_{**}(BC_2)$ vary quite significantly as we vary our base field. We invite the reader to note the differences in $kq_{**}(BC_2)$ visible in \cref{fig:kq_ASS_E2ALL}, \cref{fig:kq_ASS_E2_R}, and \cref{Fig:kqASS_F3}.

\subsection{Motivation and future directions}

Before discussing our computations in more detail, we motivate our results and suggest some directions for future work. 

\subsubsection{Equivariant motivic homotopy theory}

In classical equivariant stable homotopy theory, the classifying space $BG$ naturally arises in the study of Borel $G$-spectra, which can be used to study genuine $G$-spectra through the isotropy separation sequence. In the motivic context, the motivic classifying space $BG$ plays a similar role. Our computations of $kgl_{**}BC_2$ and $kq_{**}BC_2$ are therefore related to Borel $C_2$-equivariant algebraic and Hermitian K-theory. By studying the motivic geometric fixed points of \cite{GH23}, it should be possible to describe the $C_2$-equivariant algebraic and Hermitian K-groups of various fields after $2$-completion. It could be interesting to compare the resulting equivariant algebraic K-groups to the computations of Chan and Vogeli \cite{CV25}. 

Alternatively, one could use our $kgl$- and $kq$-homology computations to get at even richer motivic homology theories. Using our computations for $BC_2$, the fourth author and J. Morris plan to use the $kq$-based Adams spectral sequence \cite{CQ21} to compute the low-dimensional motivic stable homotopy groups of $BC_2$ over various base fields. We note that for these computations, it is really $kq$-homology, and not $kq$-cohomology, which is the most natural input. The fourth author and S. Zhang then plan to use these $\r$-motivic computations and their analogue under equivariant Betti realization to analyze the $RO(C_2 \times C_2)$-graded equivariant stable homotopy groups of spheres in a range. 

\subsubsection{Other groups}

We have focused on $G = C_2$, but in the classical setting, Bruner and Greenlees \cite{BG03, BG10} computed the connective complex and real topological K-theory of many other finite groups. It would be interesting to adapt their techniques to the motivic setting and see which $kgl$- and $kq$-homology groups can be calculated. 

The first author is currently carrying out analogous computations for $BC_{2^n}$ for all $n$, which will appear in forthcoming work. Another natural direction is to consider $C_{p^n}$ for an odd prime $p$. These computations are expected to be simpler than the $2$-primary calculations, since mod $p$ motivic homology is often more tractable for odd primes than for $p=2$. The second author plans to leverage our techniques to extend motivic analogues of earlier results on the Chow ring of  $BSO(2n,\mathbb{C})$ \cite{F12}.

Computations for elementary abelian $2$-groups should also be interesting; the group $G = (C_2)^3$ arises in Atiyah and Hirzebruch's disproof of the integral Hodge conjecture \cite{AH62} in part because its $2$-rank is large enough for the cohomology of $BG$ to support certain nontrivial Steenrod operations. As we see below, the action of the Steenrod operations on the mod two motivic cohomology of $BC_2$ and their effect via differentials in the Atiyah--Hirzebruch spectral sequence depends heavily on the base field, suggesting interesting phenomena may arise for larger $G$.

Finally, one might consider the classifying spaces which appear in the papers \cite{BW25, DLM23, Lac25, HW19, Wen24} mentioned above, plus related work. 

\subsubsection{Slice spectral sequence}

The slice filtration is a powerful tool in motivic stable homotopy theory. The effective slice filtration of $kgl$ gives rise to the spectral sequence \cite{Lev08}
$$E_1 = \bigoplus_{i \geq 0} \Sigma^{2i,i} H\z_{**} \Rightarrow kgl_{**}$$
and the very effective slice filtration for $kq$ gives rise to the spectral sequence \cite{Bac17}
$$E_1 = \bigoplus_{i \geq 0} \Sigma^{8i,4i} H\z_{**} \oplus \Sigma^{8i+1,4i+1} H\z/2_{**} \oplus \Sigma^{8i+2,4i+2} H\z/2_{**} \oplus \Sigma^{8i+4,4i+2} H\z_{**} \Rightarrow kq_{**}.$$
Smashing the slice towers for $kgl$ and $kq$ with $BG$ yields spectral sequences converging from shifts of $H\z_{**}BG$ and $H\z/2_{**}BG$ to $kgl_{**}BG$ and $kq_{**}BG$. 

As in \cite{RSO19}, it should be possible to express the $d_1$-differentials in these spectral sequence in terms of motivic Steenrod operations. If the spectral sequence collapses in a range at $E_2$, this would give a description of $kgl_{**}BG$ and $kq_{**}BG$ in a range over general base fields in terms of subquotients of the motivic homology groups of $BG$. This approach could also be used to remove some of the $2$-completions which are necessary in our approach. If the spectral sequence does not collapse at $E_2$, it would be interesting to understand the higher differentials.

\subsection{Outline of computations}

We produce our computations using a combination of the Atiyah--Hirzebruch and the motivic Adams spectral sequences. We review the classical computations of $H\z_*BC_2$, $ku_*BC_2$, and $ko_*BC_2$ via these methods in \cref{Sec:Classical}. Although these computations are already well-known, there are a variety of possible approaches and we have adopted a strategy which we believe generalizes most efficiently to the motivic context. 

Based on our classical computations, we outline a general strategy for these computations in the motivic context in \cref{Sec:Strategy}. \cref{MT} is obtained by implementing this strategy over specific base fields. The general strategy can be summarized as follows.

First, we analyze the Atiyah--Hirzebruch spectral sequence converging to $E_{**}BC_2$, whose differentials are determined by the attaching maps in $BC_2$ and multiplicative structures (products, Toda brackets) in $E_{**}$. By computing the first few pages of the Atiyah--Hirzebruch spectral sequence, we obtain bounds on the ranks of the associated graded to the groups $E_{**}BC_2$ as modules over $\m_2^F$, the mod two motivic homology groups of the base field $F$. 

Second, we analyze the motivic Adams spectral sequence converging to $E_{**}BC_2$. It turns out that this spectral sequence collapses at $E_2$ for all of our choices of $E$ and all of our base fields, so the main difficulty is computing the $E_2$-term. We develop a general strategy for decomposing $H^{**}(BC_2)$ as a module over various subalgebras of the motivic Steenrod algebra in \cref{Sec:Strategy}. The bounds obtained from our Atiyah--Hirzebruch spectral sequence computations allow us to compute various connecting homomorphisms without resorting to potentially tricky homological algebra involving the motivic Steenrod algebra.

Third, we compare the Atiyah--Hirzebruch and the Adams spectral sequences. These spectral sequences together provide us with information about additive extensions in $E_{**}BC_2$. 

The rest of the paper is spent specializing to different base fields:

In \cref{Sec:C}, we work over the complex numbers. Except for the fact that computations are bigraded and the first Hopf map $\eta$ is non-nilpotent in $kq_{**}$, the computations are essentially the same as in the classical case. A motivic Lefschetz principle \cite[Prop. 5.2.1]{BCQ25} extends our calculations to all algebraically closed fields. 

In \cref{Sec:R}, we work over the real numbers. The coefficient rings $H\z_{**}$, $kgl_{**}$, and $kq_{**}$ are significantly more complicated than their classical counterparts, and some genuinely new phenomena (i.e., involving elements which base-change to zero over $\c$) appears in the Atiyah--Hirzebruch and Adams spectral sequences. The complexity of the computations also necessitates that we divide our charts into coweights, a bookkeeping trick we initially learned from \cite{GHIR19}. A newly established ``real Lefschetz principle”, introduced at the beginning of this section, allows us to extend our calculations to all real closed fields.

In \cref{Sec:Fq}, we work over finite fields $\f_q$ with $q$ odd. Unlike in all of the previous contexts, there are differentials in the Adams spectral sequences converging to $H\z_{**}$, $kgl_{**}$, and $kq_{**}$ (first observed in \cite{Kyl15}) which makes the $2$-torsion in these groups fairly complicated. As a result, we find that there is a family of exotic Atiyah--Hirzebruch spectral sequence differentials which have no counterparts over algebraically closed and real closed fields. 

\subsection{Conventions and notation}

We adopt the following contentions and notation. 
\begin{enumerate}
\item Everything is implicitly $2$-complete.
\item $BC_2$ denotes the geometric classifying space of the cyclic group $C_2$.
\item $\ca$ denotes the Steenrod algebra and the motivic Steenrod algebra.
\item $H := H\f_2$ denotes the mod two  Eilenberg--MacLane spectrum and its motivic analogue, and $H\z$ denotes their integral analogues. 
\item $E(x_1, x_2, \ldots)$ denotes the exterior algebra generated by $x_1, x_2, \ldots$.
\item $\Ext_{\ca}^{***}(M) := \Ext_{\ca}^{***}(M,k)$ denotes $\Ext$ over the Hopf algebroid $(\mathcal{A},k)$ from an $\mathcal{A}$-module $M$ to the ground ring $k$.
\item $kgl$ denotes the effective cover of the algebraic K-theory spectrum.
\item $kq$ denotes the very effective cover of the Hermitian K-theory spectrum.
\item In this paper, we will use the cohomological grading $E^{**}$ and the homological grading $E_{**}$ interchangeably. Since their indices differ only by the standard sign convention (i.e., $E_{s,t} = E^{-s,-t}$), this distinction will not affect our computations.

\end{enumerate}

\subsection{Acknowledgments}

This work was supported by the 4-VA Collaborative Research Grant \emph{Collaborative Workshops in Topology}. The authors thank Bert Guillou and Zhenyue Guan for helpful discussions, Amy Hu{\'a} Li{\'a}ng for assistance with some charts, Arun Debray for pointing out a mistake and for additional helpful suggestions, and Dan Isaksen for helpful comments on an earlier version of this paper. The authors also thank the anonymous referee for numerous helpful comments and suggestions. AN was supported by NSF grant DMS-2231514. JDQ was supported by NSF grants DMS-2414922 and DMS-2441241.

\section{Classical calculations}\label{Sec:Classical}

All (co)homology theories in this paper are reduced. In this section, we recall the mod two cohomology of $BC_2$ and use it to compute $H\z_*BC_2$, $ku_*BC_2$, and $ko_*BC_2$. 

\subsection{Mod two cohomology} 

Recall that $BC_2 \simeq \mathbb{R}P^\infty$. We have
$$H^*(BC_2) \cong \ker\left(\f_2[u] \to \f_2, \quad u \mapsto 0\right)$$
with $|u| = 1$. The action of the first two Steenrod operations is determined by
$$\Sq^1(u) = u^2, \quad \Sq^2(u) = 0, \quad \Sq^2(u^2) = u^4.$$

\subsection{Integral homology}

To obtain $H\z_*(BC_2)$ from $H^*(BC_2)$, we may use the Adams spectral sequence with signature
$$E_2^{**} = \Ext_{\ca}^{**}(H^{*}(BC_2) \otimes H^*(H\z), \f_2) \Rightarrow H_*(BC_2).$$
Here, we are using the K{\"u}nneth isomorphism for $H\f_2$-cohomology to identify
$$H^*(BC_2 \wedge H\z) \cong H^*(BC_2) \otimes H^*(H\z).$$ 

There is an isomorphism of $\ca$-modules \cite[Proposition 19.1.3]{Margolis83}
$$H\f_2^*H\z \cong \ca//\ca(0),$$
where $\ca(0) = E(\Sq^1)$. Using a change-of-rings isomorphism, the above $E_2$-term has the form
$$E_2^{**} = \Ext_{\ca}^{**}(H^{*}(BC_2) \otimes \ca//\ca(0), \f_2) \cong \Ext_{\ca(0)}^{**}(H^*(BC_2), \f_2).$$
The action of $\Sq^1$ on $H^*(BC_2)$ described above implies that there is an isomorphism of $\ca(0)$-modules
$$H^*(BC_2) \cong  
\bigoplus_{i \geq 0} \Sigma^{2i+1} \ca(0),$$
so we may rewrite the $E_2$-term as
$$E_2^{**} = \Ext_{\ca(0)}^{**}(H^*(BC_2),\f_2) \cong 
\bigoplus_{i \geq 0} \Ext_{\ca(0)}^{**}(\Sigma^{2i+1} \ca(0), \f_2).$$
Since $\ca(0)$ is free over itself, we have
$$\bigoplus_{i \geq 0} \Ext_{\ca(0)}^{**}(\Sigma^{2i+1} \ca(0), \f_2) \cong \bigoplus_{i \geq 0} \Sigma^{2i+1} \f_2.$$
The Adams spectral sequence collapses at $E_2$ for degree reasons, and there is no room for hidden extensions. We conclude that
\[
H\z_*(BC_2) \cong 
\begin{cases}
    \f_2 \quad & \text{ if } *=2i+1, i\geq 0, \\
    0 \quad & \text{ otherwise.}
\end{cases}
\]

\subsection{Connective complex K-theory} \label{section: connective complex K-theory}

To obtain $ku_*(BC_2)$ from $H^*(BC_2)$, we proceed in two steps. First, we will run the Atiyah--Hirzebruch spectral sequence to calculate the rank of the associated graded of $ku_*(BC_2)$ in all degrees. Second, we will use our Atiyah--Hirzebruch calculations to determine the $E_2$-term of the Adams spectral sequence, which will then resolve any questions about additive extensions.

\begin{figure}[H]
\includegraphics{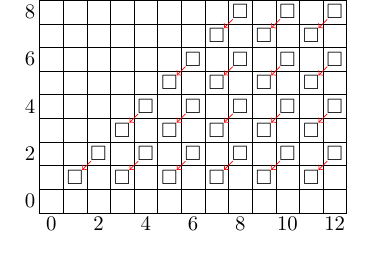}
\includegraphics{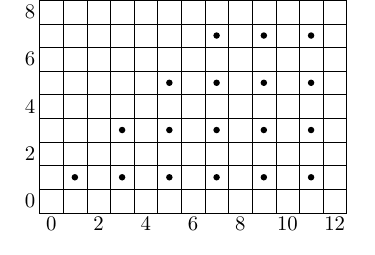}
\caption{The Atiyah--Hirzebruch spectral sequence converging to $ku_*BC_2$ (left) and its $E_\infty$-page. Each $\square$ represents $\z$ and each differential (drawn in red) represents multiplication by $2$. Each bullet represents $\f_2$.}\label{fig:kuAH}
\end{figure}

The Atiyah--Hirzebruch spectral sequence has the form
$$E_1 = \bigoplus_{\text{cells of } BC_2} ku_* \Rightarrow gr^{*}_{AH}ku_*BC_2.$$
Bott periodicity implies that  $ku_* \cong \z[\beta]$ with $|\beta| = 2$. Since each even-dimensional cell of $BC_2$ is attached to the odd-dimensional cell one dimension lower by a degree-$2$ attaching map, the nontrivial $d_1$-differentials are given by
$$d_1(x[2n]) = 2x[2n-1]$$
for all $n \geq 1$, where $x[i]$ is our notation for an element in filtration $i$ with $x \in ku_*$. See \cref{fig:kuAH}.

The $E_2$-term consists of
a copy of $ku_*/2$ in each positive odd filtration  and is zero otherwise. The spectral sequence collapses at $E_2$ for degree reasons, showing
$$gr^*_{AH} ku_n(BC_2) \cong 
\begin{cases}
    (\f_2)^{i+1} \quad & \text{ if } n = 2i+1 \geq 0, \\
    0 \quad & \text{ otherwise.}
\end{cases}$$
Here, we write $gr^*_{AH}$ for the associated graded of the Atiyah--Hirzebruch filtration. 

We resolve these hidden extensions using the Adams spectral sequence. This has signature
$$E_2^{**} = \Ext_{\ca}^{**}(H^*(BC_2) \otimes H^*(ku), \f_2) \Rightarrow ku_*(BC_2).$$
There is an isomorphism of $\ca$-modules
$$H\f_2^*(ku) \cong \ca//\ce(1),$$
where $\ce(1) = E(Q_0, Q_1)$ with $Q_0 = \Sq^1$ and $Q_1 = [\Sq^1, \Sq^2] = \Sq^1 \Sq^2 + \Sq^2\Sq^1.$
After a change-of-rings isomorphism, we may rewrite the $E_2$-term as
$$E_2^{**} = \Ext_{\ce(1)}^{**}(H^*(BC_2), \f_2).$$

We will calculate these $\Ext$-groups as follows. There is a short exact sequence of $\ce(1)$-modules
$$
\Sigma^2 \f_2 \hookrightarrow H^*(BC_2) \to V,$$
where $V$ is defined to be the quotient of the inclusion. This induces a long exact sequence 
\begin{equation}\label{eqn:kuBC2}
    \cdots \leftarrow \Ext^{s,t}_{\ce(1)}(
    \Sigma^2 \f_2) \leftarrow \Ext^{s,t}_{\ce(1)}(H^*(BC_2)) \leftarrow \Ext^{s,t}_{\ce(1)}(V) \xleftarrow{\delta} \Ext^{s-1,t}_{\ce(1)}(
    \Sigma^2 \f_2) \leftarrow \cdots
\end{equation}
which we may use to calculate $\Ext_{\ce(1)}(H^*(BC_2))$. 

Since $\ce(1)$ is an exterior algebra on two generators, we have
$$\Ext_{\ce(1)}(\Sigma^2\f_2) \cong \Sigma^2\f_2[h_0,v_1],$$
where $|h_0| = (1,1)$ and $|v_1| = (3,1)$. 

On the other hand, we may compute $\Ext_{\ce(1)}(V)$ as follows. Observe that $\ce(1)//\ce(0) = E(Q_1)$. Define a filtration on $V$ by 
$$F_i = \f_2\{x,x^4,\ldots,x^{2i+1},x^{2i+4} \}$$ 
so that 
$$F_i/F_{i-1} \cong \Sigma^{2i+1}E(Q_1).$$
Applying $\Ext_{\ce(1)}(-)$ gives a spectral sequence with signature
$$E_1 = \bigoplus_{i \geq 0} \Ext_{\ce(1)}(\Sigma^{2i+1}E(Q_1)) \Rightarrow \Ext_{\ce(1)}(V).$$
By a change-of-rings isomorphism, we may rewrite the $E_1$-term as
$$E_1 = \bigoplus_{i \geq 0} \Ext_{\ce(0)}(\Sigma^{2i+1}\f_2) \cong \bigoplus_{i \geq 0} \Sigma^{2i+1} \f_2[h_0].$$
Since these $h_0$-towers are all at least $2$ stems apart from each other, the spectral sequence collapses for degree reasons, giving 
$$\Ext_{\ce(1)}(V) \cong \bigoplus_{i \geq 0} \Sigma^{2i+1} \f_2[h_0].$$

\begin{figure}[H]
\includegraphics{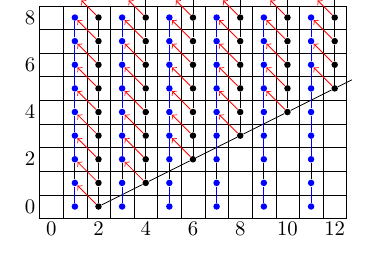}
\includegraphics{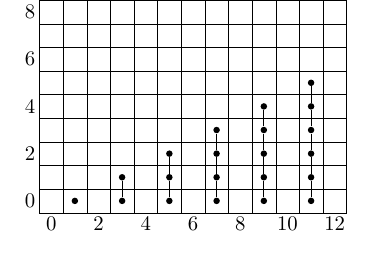}
\caption{The long exact sequence computing $\Ext_{\ce(1)}^{**}(H^*(BC_2))$ (left) and the groups $\Ext_{\ce(1)}(H^*(BC_2))$ (right). In the long exact sequence figure, contributions from $\Ext_{\ce(1)}(V)$ appear in blue and contributions from $\Ext_{\ce(1)}(\Sigma^2 \f_2)$ appear in black. In both figures, a bullet represents $\f_2$.}\label{fig:kuBC2Ext}
\end{figure}

Having calculated the groups in \eqref{eqn:kuBC2}, it only remains to determine the connecting homomorphism. We claim the connecting homomorphism is nontrivial, with
$$\delta(\Sigma^{2} v_1^i h_0^j) = \Sigma^{2i+1} h_0^{i+j+1}.$$
Indeed, if this were not the case, then \eqref{eqn:kuBC2} would split, decomposing $\Ext_{\ce(1)}(H^*(BC_2))$ as the direct sum of $\Ext_{\ce(1)}(\f_2 \oplus \Sigma^2 \f_2)$ and $\Ext_{\ce(1)}(V)$. However, such a decomposition would force the smallest possible size for $ku_{2i+1}(BC_2)$ to be $(\f_2)^{i+2}$, contradicting our Atiyah--Hirzebruch spectral sequence computations. Therefore the connecting homomorphism must be as claimed. The calculation is depicted in \cref{fig:kuBC2Ext}. 

The $E_2$-term of the Adams spectral sequence is thus given by
$$E_2 = \Ext_{\ce(1)}(H^*(BC_2)) \cong 
\bigoplus_{i \geq 0} \Sigma^{2i+1}\f_2[h_0]/(h_0)^{i+1}.$$
The spectral sequence collapses for degree reasons. Since $h_0$ detects multiplication by $2$, we conclude
$$ku_n(BC_2) \cong 
\begin{cases}
    \z/2^{i+1} \quad & \text{ if } n = 2i+1 \geq 0, \\
    0 \quad & \text{ otherwise.}
\end{cases}$$
where the generators are given by the $ku$-fundamental class of the $(2i+1)$-skeleton $\mathbb{RP}^{2i+1}$ of $BC_2=\mathbb{RP}^\infty$.
\subsection{Connective real K-theory}

We compute $ko_*(BC_2)$ using a similar strategy to our computation of $ku_*(BC_2)$. 

\begin{figure}[H]
\includegraphics[scale=.9]{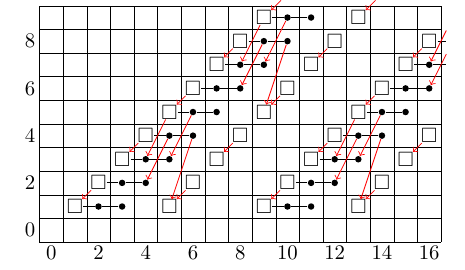}
\includegraphics[scale=.9]{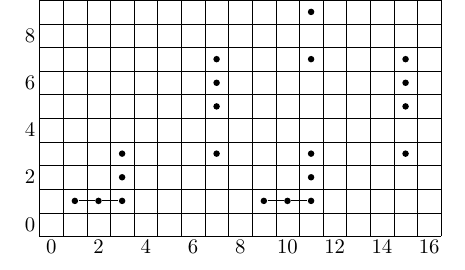}
\caption{The Atiyah--Hirzebruch spectral sequence converging to $ko_*BC_2$ (left) and the $E_\infty$-term (right). Each $\square$ represents $\z$ and each $\bullet$ represents $\f_2$. Horizontal lines represent multiplication by $\eta$. Differentials $d_1$, $d_2$, and $d_3$ (combined) are shown in red.}\label{fig:koAH}. 
\end{figure}

We begin with the Atiyah--Hirzebruch spectral sequence
$$E_1 = \bigoplus_{i \geq 0} \Sigma^i ko_* \Rightarrow ko_*(BC_2).$$
The groups $ko_*$ can be computed using the classical Adams spectral sequence. The potentially nontrivial differentials are determined by
$$d_1(x[n]) = 2x[n-1]$$
for all  $n = 2k > 0$,
$$d_2(x[n]) = \eta x[n-2]$$
for all $n \equiv 0,1 \pmod 4$ with $n \geq 4$, and
$$d_3(\eta^2[4k]) = 2\alpha[4k-3]$$
for all $k \geq 1$. 
The spectral sequence then collapses since the only potential longer differentials are from torsion groups to torsion-free groups. We find
$$gr^*_{AH}ko_n(BC_2) \cong \begin{cases}
    \f_2 \quad & \text{ if } n = 8i+1,8i+2 > 0, \\
    (\f_2)^{4i+3} \quad & \text{ if } n = 8i+3 > 0, \\
    (\f_2)^{4i+4} \quad & \text{ if } n = 8i+7 > 0, \\
    0 \quad & \text{ otherwise.}    
\end{cases}$$
The spectral sequence is depicted in \cref{fig:koAH}. 

To resolve hidden extensions, we use the Adams spectral sequence with signature
$$E_2^{**} = \Ext_{\ca}^{**}(H^*(BC_2;\f_2) \otimes H^*(ko;\f_2), \f_2) \Rightarrow ko_*(BC_2).$$
Using the isomorphism of $\ca$-modules 
$$H\f_2^*(ko) \cong \ca // \ca(1),$$
where $\ca(1) = \langle \Sq^1, \Sq^2 \rangle$, we may rewrite the $E_2$-term as
$$E_2^{**} = \Ext_{\ca(1)}^{**}(H^*(BC_2;\f_2),\f_2).$$
Recall that there is an isomorphism of $\ca(1)$-modules
$$\ca(1)//\ca(0) \cong \f_2\{1, \Sq^2, \Sq^1 \Sq^2, \Sq^2 \Sq^1 \Sq^2\}.$$

We will calculate $\Ext_{\ca(1)}(H^*(BC_2;\f_2))$ as follows. There is a short exact sequence of $\ca(1)$-modules
$$ Q \hookrightarrow H^*(BC_2;\f_2) \to R,$$
where $ Q \subseteq H^*(BC_2;\f_2)$ is the submodule generated by $x$, $x^2$, and $x^4$, i.e., $Q$ is the question mark complex. This gives a long exact sequence in $\Ext$
\begin{equation}\label{eqn:koBC2}
\cdots \leftarrow \Ext_{\ca(1)}(Q) \leftarrow \Ext_{\ca(1)}(H^*(BC_2;\f_2)) \leftarrow \Ext_{\ca(1)}(R) \xleftarrow{\delta} \Ext_{\ca(1)}( Q) \leftarrow \cdots.
\end{equation}

\begin{figure}[H]
\includegraphics{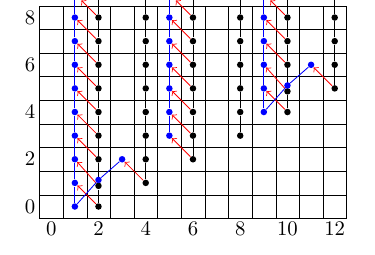}
\includegraphics{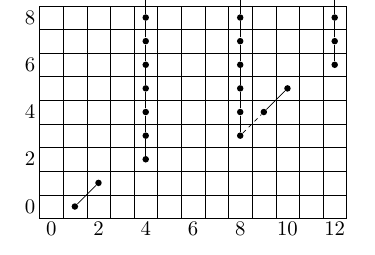}
\caption{The long exact sequence computing $\Ext_{\ca(1)}(Q)$ (left) and the groups $\Ext_{\ca(1)}(Q)$ (right). In the long exact sequence figure, contributions from $\Ext_{\ca(1)}(\f_2)$ appear in blue and contributions from $\Sigma^2 \Ext_{\ce(1)}(\f_2)$ appear in black. Each bullet represents $\f_2$.}\label{fig:Q_classical}
\end{figure}

To calculate $\Ext_{\ca(1)}(Q)$, we use the short exact sequence of $\ca(1)$-modules
$$\Sigma^2 E(\Sq^2) \to Q \to \Sigma^1 \f_2.$$
The groups $\Ext_{\ca(1)}(\f_2)$ are classical. On the other hand, there is an isomorphism of $\ca(1)$-modules $\ca(1)//\ce(1) \cong E(\Sq^2)$, so $\Ext_{\ca(1)}(E(\Sq^2)) \cong \Ext_{\ce(1)}(\f_2)$ was described in the previous subsection. The connecting homomorphism is nontrivial, either by direct computation or by arguing that if it did not happen, the size of $ko_n(BC_2)$ would be too large. The long exact sequence and final computation of $\Ext_{\ca(1)}(Q)$ are depicted in \cref{fig:Q_classical}. Note that there is a collection of hidden $h_1$-extensions, depicted by dashed line in the right-hand figure of \cref{fig:Q_classical}.

To calculate $\Ext_{\ca(1)}(R)$, we filter $R$ by
$$F_i = \{x^3, x^5, x^6, \ldots, x^{4i+3}, x^{4i+5}, x^{4i+6}, x^{4i+8}\}$$
so that
$$F_i/F_{i-1} \cong \ca(1)//\ca(0).$$
Applying $\Ext_{\ca(1)}(-)$ gives a spectral sequence with signature
$$E_1 = \bigoplus_{i \geq 0} \Ext_{\ca(1)}(F_i/F_{i-1}) = \Ext_{\ca(0)}(\Sigma^{4i+3} \f_2) \Rightarrow \Ext_{\ca(1)}(R).$$
The spectral sequence collapses for degree reasons, so we obtain
$$\Ext_{\ca(1)}(R) \cong \bigoplus_{i \geq 0} \Sigma^{4i+3} \f_2[h_0].$$

\begin{figure}[H]
\includegraphics{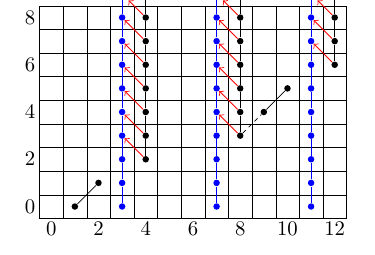}
\includegraphics{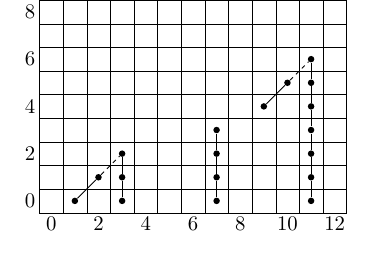}
\caption{The long exact sequence computing $\Ext_{\ca(1)}(H^*(BC_2))$ (left) and the groups $\Ext_{\ca(1)}(H^*(BC_2))$ (right). In the long exact sequence figure, contributions from $\Ext_{\ca(1)}(Q)$ appear in black and contributions from $\Ext_{\ca(1)}(R)$ appear in blue. Each bullet represents $\f_2$.}\label{fig:koBC2}

\end{figure}

We now analyze \eqref{eqn:koBC2}. 
The connecting homomorphism must be nontrivial whenever possible for it to agree with our earlier description of $gr^*_{AH}kgl_{**}BC_2$. The computation is depicted in \cref{fig:koBC2}. 

\begin{figure}[H]
\includegraphics{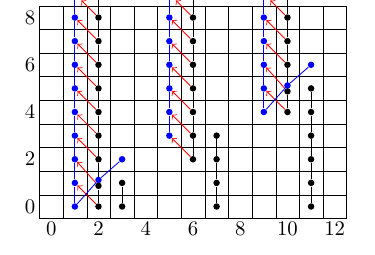}
\includegraphics{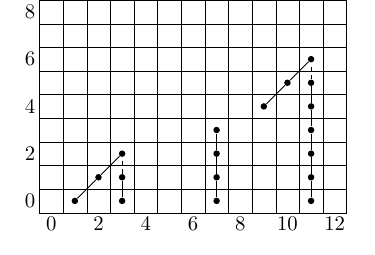}
\caption{The alternative long exact sequence computing $\Ext_{\ca(1)}(H^*(BC_2))$ (left) and the groups $\Ext_{\ca(1)}(H^*(BC_2))$ (right). In the long exact sequence figure, contributions from $\Ext_{\ca(1)}(R_1)$ appear in black and contributions from $\Ext_{\ca(1)}(\f_2)$ appear in blue. Each bullet represents $\f_2$. Note that, in this case, there is instead a collection of hidden $h_0$-extensions, depicted by dashed lines.}\label{fig:koBC2_h1}
\end{figure}

This gives the \( E_2 \)-term of the Adams spectral sequence up to one collection of hidden \( h_1 \)-extensions, depicted by dashed lines in \cref{fig:koBC2}. It turns out that these extensions do occur. Consider the short exact sequence
\[
 \Sigma^2 R_1 \hookrightarrow H^*(BC_2) \twoheadrightarrow \Sigma^1 \f_2,
\]
where $R_1$ is the kernel. It follows that $\Ext_{\ca(1)}(R_1)$ is computed in \cite[Figure 25, 26]{BC18}. Applying $\Ext_{\ca(1)}(-)$ to the aforementioned short exact sequence shows that a hidden $h_1$-extension occurs. The result is depicted in \cref{fig:koBC2_h1}.

We conclude that
$$ko_n(BC_2) \cong
\begin{cases}
\z/2 \quad & \text{ if } n = 8i+1,8i+2 > 0, \\
\z/2^{4i+3} \quad & \text{ if } n = 8i+3 > 0, \\
\z/2^{4i+4} \quad & \text{ if } n = 8i+7 > 0, \\
0 \quad & \text{ otherwise.}
\end{cases}$$

\section{General strategy}\label{Sec:Strategy}

The classical computations of the previous section provide a blueprint for analogous motivic computations. In this section, we spell out this general strategy which will be used in subsequent sections to make computations over specific base fields. 

\subsection{Mod two motivic cohomology}

Let $\m_2^{**} = H^{**}$ and let $\m_2 = H_{**}$. In $\m_2^{**}$, there are distinguished elements $\tau \in \m_2^{0,1}$ and $\rho \in \m_2^{1,1}$.  

In \cite[Thm. 6.10]{Voe03}, Voevodsky produced an isomorphism of bigraded rings
$$H^{**}BC_2 \cong \left( \ker: \m_2^{**}[x,y]/(x^2+\tau y + \rho x) \to \m_2^{**}, \quad x,y \mapsto 0 \right),$$
where $|x| = (1,1)$ and $|y| = (2,1)$. 

We will be interested in $H^{**}BC_2$ as a module over $\ca(1)$, the subalgebra of the motivic Steenrod algebra generated by $\Sq^1$ and $\Sq^2$. The actions of these Steenrod operations on $H^{**}BC_2$ are determined by 
$$\Sq^1(x) = y, \quad \Sq^2(y) = y^2,$$
together with the Cartan formula and the Adem relations:

\begin{lemma}\label{lem:Sq12C2}
Let $0 \leq i \leq 1$ and let $j \geq 0$. Then
$$\Sq^1(x^i y^j) = 
\begin{cases}
y^{j+1} \quad & \text{ if } i = 1, \\
0 \quad & \text{ if } i = 0,
\end{cases}$$
$$\Sq^2(x^i y^j) = 
\begin{cases}
x^i y^{j+1} \quad & \text{ if } i+2j \equiv 2,3 \pmod 4, \\
0 \quad & \text{ if } i+2j \equiv 0,1 \pmod 4.
\end{cases}$$
\end{lemma}

\subsection{Integral motivic homology}

The motivic Adams spectral sequence converging to $H\z_{**}BC_2$ has signature
$$E_2 = \Ext_{\ca}^{***}(H^{**}(BC_2 \wedge H\z)) \Rightarrow H\z_{**}BC_2.$$
Since $H^{**}(H\z) \cong \ca//\ca(0)$, we can use the K{\"u}nneth theorem and a change-of-rings isomorphism to rewrite
$$E_2 \cong \Ext_{\ca}^{***}(H^{**}(BC_2) \otimes \ca//\ca(0)) \cong \Ext_{\ca(0)}^{***}(H^{**}BC_2).$$
By \cref{lem:Sq12C2}, there is an isomorphism of $\ca(0)$-modules
$$H^{**}BC_2 \cong \bigoplus_{i \geq 0} \Sigma^{2i+1,i+1} \ca(0),$$
which implies that
$$\Ext_{\ca(0)}^{***}(H^{**}(BC_2)) \cong \bigoplus_{i \geq 0} \Sigma^{2i+1,i+1} \Ext_{\ca(0)}^{***}(\ca(0)) \cong \bigoplus_{i \geq 0} \Sigma^{2i+1,i+1}\m_2.$$
Since $\m_2$ is concentrated in Adams filtration zero, the Adams spectral sequence collapses with no possibility of hidden extensions, showing:

\begin{proposition}\label{Prop:HZBC2}
Over any base field, there is an isomorphism of bigraded abelian groups
$$H\z_{**}BC_2 \cong \bigoplus_{i \geq 0} \Sigma^{2i+1,i+1} \m_2.$$
\end{proposition}

\subsection{Effective algebraic K-theory}\label{SS:kgl_general}

The motivic Adams spectral sequence converging to 2-completed $kgl_{**}BC_2$ has signature
$$E_2 = \Ext_{\ca}^{***}(H^{**}(BC_2 \wedge kgl)) \Rightarrow kgl_{**}BC_2.$$
Recall that $\ce(1) \subseteq \ca$ is the subalgebra generated by $\Sq^1$ and $Q_1 \coloneqq [\Sq^1, \Sq^2]$. Using the isomorphism of $\ca$-modules $H^{**}kgl \cong \ca//\ce(1)$ (\cite[Rmk. 2.10]{ARO17}), we may rewrite the $E_2$-term as
\begin{equation*}
    E_2 \cong \Ext_{\ce(1)}^{***}(H^{**}BC_2).
\end{equation*}
To access these $\Ext$-groups, we analyze the long exact sequence obtained by applying $\Ext_{\ce(1)}^{***}(-)$ to the short exact sequence of $\ce(1)$-modules
\begin{equation} \label{eqn: les for kgl general strategy}
    \Sigma^{2,1} \m_2^{**} \hookrightarrow H^{**}BC_2 \to V,
\end{equation}
where $V$ is the cokernel.

To compute $\Ext_{\ce(1)}^{***}(V)$, we consider the spectral sequence arising from the $\ce(1)$-module filtration of $V$ defined by
$$F_i := 
\begin{cases}
    0 \quad & \text{ if } i < 0, \\
    \m_2^{**}\{x,y^2\} \quad & \text{ if } i = 0, \\
    \m_2^{**}\{x,xy,y^2,xy^2,y^3,\ldots,xy^i,y^{i+1}\} \quad & \text{ if } i > 0.
\end{cases}$$
Observing $\ce(1)//\ca(0) \cong E(Q_1)$, we see that
$$F_i / F_{i-1} \cong \Sigma^{2i+1,i+1}\ce(1)//\ca(0)$$
for each $i \geq 0$. We thus obtain a spectral sequence with signature
\begin{equation} \label{eqn: E1-page for kgl general strategy}
    E_1 = \bigoplus_{i \geq 0} \Ext_{\ce(1)}^{***}(\Sigma^{2i+1,i+1} \ce(1)//\ca(0)) \cong \bigoplus_{i \geq 0} \Sigma^{2i+1,i+1} \Ext_{\ca(0)}^{***}(\m_2^{**}) \Rightarrow \Ext_{\ce(1)}^{***}(V).
\end{equation}
We will analyze this spectral sequence over particular base fields. 

\subsection{Very effective Hermitian K-theory}\label{kq_{**}BC_2}

The motivic Adams spectral sequence converging to 2-completed $kq_{**}BC_2$ has signature
$$E_2 = \Ext_{\ca}^{***}(H^{**}(BC_2 \wedge kq)) \Rightarrow kq_{**}BC_2.$$
Recall that $\ca(1) \subset \ca$ is the subalgebra generated by $\Sq^1$ and $\Sq^2$. Using the isomorphism of $\ca$-modules $H^{**}kq \cong \ca//\ca(1)$ (\cite[Rmk. 2.14]{ARO17}), we may rewrite the $E_2$-term as
$$E_2 \cong \Ext^{***}_{\ca(1)}(H^{**}BC_2).$$
To access these $\Ext$-groups, we analyze the long exact sequence obtained by applying $\Ext^{***}_{\ca(1)}(-)$ to the short exact sequence of $\ce(1)$-modules
\begin{equation} \label{eqn: les for kq with Q and R general strategy}
    Q \hookrightarrow H^{**}BC_2 \to R,
\end{equation}
where $Q = \m_2^{**}\{x,y,y^2\}$ and  $R = H^{**}BC_2/Q \cong \m_2^{**}\{xy,xy^2,y^3,\ldots\}$. 
To compute $\Ext_{\ca(1)}^{***}(Q)$, we use the short exact sequence of $\ca(1)$-modules
\begin{equation} \label{eqn: les for kq general strategy}
    \Sigma^{2,1}C \to Q \to \Sigma^{1,1}\m_2^{**},
\end{equation}
where $\Sigma^{2,1}C = \m_2^{**}\{y,y^2\}$. There is an isomorphism of $\ca(1)$-modules $C \cong \ca(1)//\ce(1)$, so analyzing the resulting long exact sequence in $\Ext_{\ca(1)}^{***}(-)$ reduces to understanding $\Ext_{\ce(1)}^{***}(\m_2^{**})$,  $\Ext_{\ca(1)}^{***}(\m_2^{**})$, and a connecting homomorphism. We analyze these groups and the connecting homomorphism over various base fields in the following sections. 

To compute $\Ext_{\ca(1)}^{***}(R)$, we define an $\ca(1)$-module filtration of $R$ by 
\begin{equation} \label{eqn: standard filtration of R in kq general strategy}
    F_i := 
    \begin{cases}
        0 \quad & \text{ if } i<0, \\
        \m_2\{xy, xy^2,y^3,y^4\} \quad & \text{ if } i = 0, \\
        \m_2\{xy, xy^2, y^3, xy^2, y^4, \ldots, xy^{2i+1}, xy^{2i+2}, y^{2i+3}, y^{2i+4}\} \quad & \text{ if } i > 0.
    \end{cases}
\end{equation}

Note that $\ca(1)//\ca(0) \cong \m_2^{**}\{1,\Sq^2,\Sq^2\Sq^1, \Sq^2 \Sq^1 \Sq^2\}$, so we have
$$F_i / F_{i-1} \cong \Sigma^{4i+3,2i+2} \ca(1)//\ca(0).$$
We thus obtain a spectral sequence with signature
\begin{equation} \label{eqn: E1-page for kq general strategy}
    E_1 = \bigoplus_{i \geq 0} \Sigma^{4i+3,2i+2} \Ext_{\ca(0)}^{***}(\m_2^{**}) \Rightarrow \Ext_{\ca(1)}^{***}(R).
\end{equation}
We analyze this spectral sequence over different base fields in the following sections. 
\section{Computations over algebraically closed fields}\label{Sec:C}

We now compute $H\z_{**}BC_2$, $kgl_{**}BC_2$, and $kq_{**}BC_2$ from $H\f_2^{**}BC_2$ over the complex numbers, using techniques analogous to the classical setting outlined in \cref{Sec:Strategy}. We note that by the motivic Lefschetz principle (\cite[Prop. 5.2.1]{BCQ25}), our computations here imply identical results over any algebraically closed field. 

\subsection{Mod two motivic cohomology}

We have $H\f_2^{**} = \m_2^{**} = \f_2[\tau]$. A foundational computation of Voevodsky (\cite{Voe03}, Theorem 6.10) shows that 
$$H^{**}(BC_2) \cong \ker\left( \f_2[\tau][x,y]/(x^2+\tau y) \to \f_2[\tau], \quad x,y \mapsto 0\right)$$
with $|x| = (1,1)$ and $|y| = (2,1)$. 

\subsection{Integral motivic homology}

Using \cref{Prop:HZBC2}, we have
$$H\z_{**}(BC_2) \cong \bigoplus_{i \geq 0} \Sigma^{2i+1,i+1} \f_2[\tau].$$
As a bigraded abelian group, we have
$$H\z_{n,w}(BC_2) \cong 
\begin{cases}
    \f_2 \quad & \text{ if } n=2i+1, \ w\leq i+1, \\
    0 \quad & \text{ otherwise.}
\end{cases}$$

\subsection{Effective algebraic K-theory}

As in the classical case, we analyze the Atiyah--Hirzebruch and motivic Adams spectral sequences. All our charts in this section are the same as the charts for connective complex $K$-theory in \cref{section: connective complex K-theory}, except $\bullet$ now represents $\mathbb{M}_2$ instead of $\mathbb{F}_2$.

To start, we compute $kgl_{**}$ using the motivic Adams spectral sequence with signature
$$E_2 = \Ext_{\ca}^{***}(H^{**}kgl) \cong \Ext_{\ce(1)}^{***}(\m_2^{**}) \Rightarrow kgl_{**}.$$
Since $\ce(1)$ is an exterior algebra over $\m_2$ on $\Sq^1$ and $Q_1$, we have
$$E_2 = \m_2[h_0,v_1],$$
where $|h_0| = (1,1,0)$ and $|v_1| = (1,3,1)$. There is no room for differentials, so $E_2 = E_\infty$. There is also no room for hidden extensions, so we find
$$(kgl_{n,w}) \cong 
\begin{cases}
    \z_2 \quad & \text{ if } n=2i \geq 0, \ w \leq i, \\
    0 \quad & \text{ otherwise.}
\end{cases}$$

We now run the Atiyah--Hirzebruch spectral sequence, which has signature
$$E_1 = \bigoplus_{i \geq 0} \Sigma^{2i+1,i+1} kgl_{**} \oplus \bigoplus_{j \geq 1} \Sigma^{2j,j} kgl_{**} \Rightarrow kgl_{**}BC_2.$$
As in the classical case, we have
$$d_1(x[2n]) = 2x[2n-1]$$
for all $n \geq 1$. The spectral sequence collapses at $E_2$ for degree reasons, so we find
$$gr^*_{AH} kgl_{n,w} \cong 
\begin{cases}
    (\f_2)^{i+1} \quad & \text{ if } n = 2i+1 > 0, \ w \leq i+1, \\
    0 \quad & \text{ otherwise.}
\end{cases}
$$

We follow \cref{SS:kgl_general} to resolve hidden extensions. In the long exact sequence of $\Ext$ derived from \eqref{eqn: les for kgl general strategy}, we have already determined
$$\Ext_{\ce(1)}^{***}(\Sigma^{2,1} \m_2^{**}) \cong \Sigma^{2,1} \m_2[h_0,v_1].$$
On the other hand, the computation of $\Ext_{\ca(0)}^{***}(\m_2^{**})$ is well known (e.g. \cite[p. 541]{GI15b}). So in \eqref{eqn: E1-page for kgl general strategy}, we obtain
$$E_1 = \bigoplus_{i \geq 0} \Ext_{\ca(0)}^{***}(\Sigma^{2i+1,i+1} \m_2^{**}) \cong \bigoplus_{i \geq 0} \Sigma^{2i+1,i+1}\m_2[h_0] \Rightarrow \Ext_{\ce(1)}^{**}(V),$$
which collapses for degree reasons. 

Returning to the long exact sequence in $\Ext$ obtained from the short exact sequence of $\ce(1)$-modules, we see that the connecting homomorphism must be given by
$$\delta(\Sigma^{2,1} v_1^i) = \Sigma^{2i+1,i+1}h_0^{i+1}$$
for it to agree with our earlier description of the associated graded. This gives the $E_2$-term of the Adams spectral sequence. The spectral sequence collapses for degree reasons and we get
\[
E_\infty \cong \bigoplus_{i \ge 0} \Sigma^{2i+1, i+1} \frac{\mathbb{F}_2[\tau, h_0]}{(h_0^{i+1})}
\]
As a bigraded abelian group, we have
$$
kgl_{n,w}BC_2 \cong 
\begin{cases}
    \z/2^{i+1}\quad & \text{ if } n = 2i+1 > 0, \ w \leq i+1, \\
    0 \quad & \text{ otherwise.}
\end{cases}$$
\subsection{Very effective Hermitian K-theory}\label{kqC}

Following the same steps as outlined in \cref{Sec:Strategy}, we will compute $kq_{**}BC_2$. 

We may compute $kq_{**}$ using the motivic Adams spectral sequence with signature
$$E_2 = \Ext_{\ca(1)}^{***}(\m_2^{**}) \Rightarrow kq_{**}.$$
The $E_2$-term can be computed using a minimal resolution (\cite{IS11}), but we will use the motivic May spectral sequence \cite{DI10}, whose signature is
$$E_1 = \m_2[h_0,h_1,h_{20}] \Rightarrow \Ext_{\ca(1)}^{***}(\m_2^{**}),$$
with differentials induced from the classical filtration of the cobar complex. The $d_1$-differential is given by the coproduct structure of the dual motivic Steenrod algebra:
$$d_1(h_{20}) = h_0h_1, \quad d_1(h_0)=d_1(h_1)=0$$
The $d_2$-differential is determined by how it acts on $h_0, h_1,$ and $h_{20}^2$. Note that $h_0$ and $h_1$ must be killed by $d_2$ for degree reasons, and we have 
$$d_2(h_{20}^2) = \tau h_1^3.$$
Writing $b_{20} \coloneqq [h_{20}^2]$, $a \coloneqq [h_0 b_{20}]$, and $b \coloneqq [b_{20}^2]$, we obtain
$$\Ext_{\ca(1)}^{***}(\m_2^{**}) \cong \dfrac{\m_2[h_0,h_1,a,b]}{(h_0h_1,\tau h_1^3, h_1 a, a^2 - h_0^2 b)}.$$
The Adams spectral sequence collapses for degree reasons, so we get an isomorphism of $kq_{**}$ with the ring above. Let $2 = [h_0]$, $\eta = [h_1]$, $\alpha = [a]$, and $\beta = [b]$. A chart can be found at \cite[Figure A.1]{IS11}. 

This provides the input for the Atiyah--Hirzebruch spectral sequence, whose signature is
$$E_1 = \bigoplus_{i \geq 0} \Sigma^{2i+1,i+1} kq_{**} \oplus \bigoplus_{j \geq 0} \Sigma^{2j,j} kq_{**} \Rightarrow kq_{**}BC_2$$
with differentials induced by the standard cellular structure of $BC_2$. We refer to \cite[Sec. 2.3]{Qui19a} for a discussion of differentials in the Atiyah--Hirzebruch spectral sequence for $L^\infty_{-\infty}$; this spectral sequence is essentially identical to ours in positive filtration. We see that $d_1$ is given by multiplication by $2$ on even cells and $d_2$ by multiplication by $\eta$ on cells of dimensions $4k, 4k+1\geq 4$:
$$d_1(x[2n]) = 2x[2n-1], \quad n \geq 1;\quad d_2(x[n]) = \eta x[n-2], \quad n = 4k, 4k+1 \geq 4. $$
By the above description and the relations on $kq_{**}$, we see that $\tau\eta^2$ in degree $4k$ survives to the $E_3$ page. This forces 
$$d_3(\tau \eta^2[4n]) = \alpha[4n-3].$$
The spectral sequence collapses at the $E_4$ page for degree reasons. We conclude
$$gr^*_{AH}kq_{n,w}BC_2 \cong 
\begin{cases}
    (\f_2)^{\oplus i} \quad & \text{ if } n = 8i > 0, \ w \leq 8i, \\
    (\f_2)^{\oplus (i+1)}\quad & \text{ if } n = 8i+j, \ w \leq 4i+j,  j \in \{ 1, 2, 4, 5, 6 \},\\
    (\f_2)^{\oplus (5i+3)}\quad & \text{ if } n = 8i+3, \ w \leq 4i+3, \\
    (\f_2)^{\oplus (5i+5)}\quad & \text{ if } n = 8i+7, \ w \leq 4i+7.
\end{cases}$$

\begin{remark}
This Atiyah--Hirzebruch spectral sequence can be obtained by truncating \cite[Fig. 6]{Qui19a} above filtration $1$. That figure computes the geometric $C_2$-Tate construction of $kq$ with trivial $C_2$-action, which maps to the suspension of the geometric $C_2$-homotopy orbits of $kq$; the latter identifies with the suspension of $kq_{**}BC_2$. 
\end{remark}

\begin{figure}[H]
\centering
\includegraphics[scale = 0.15]{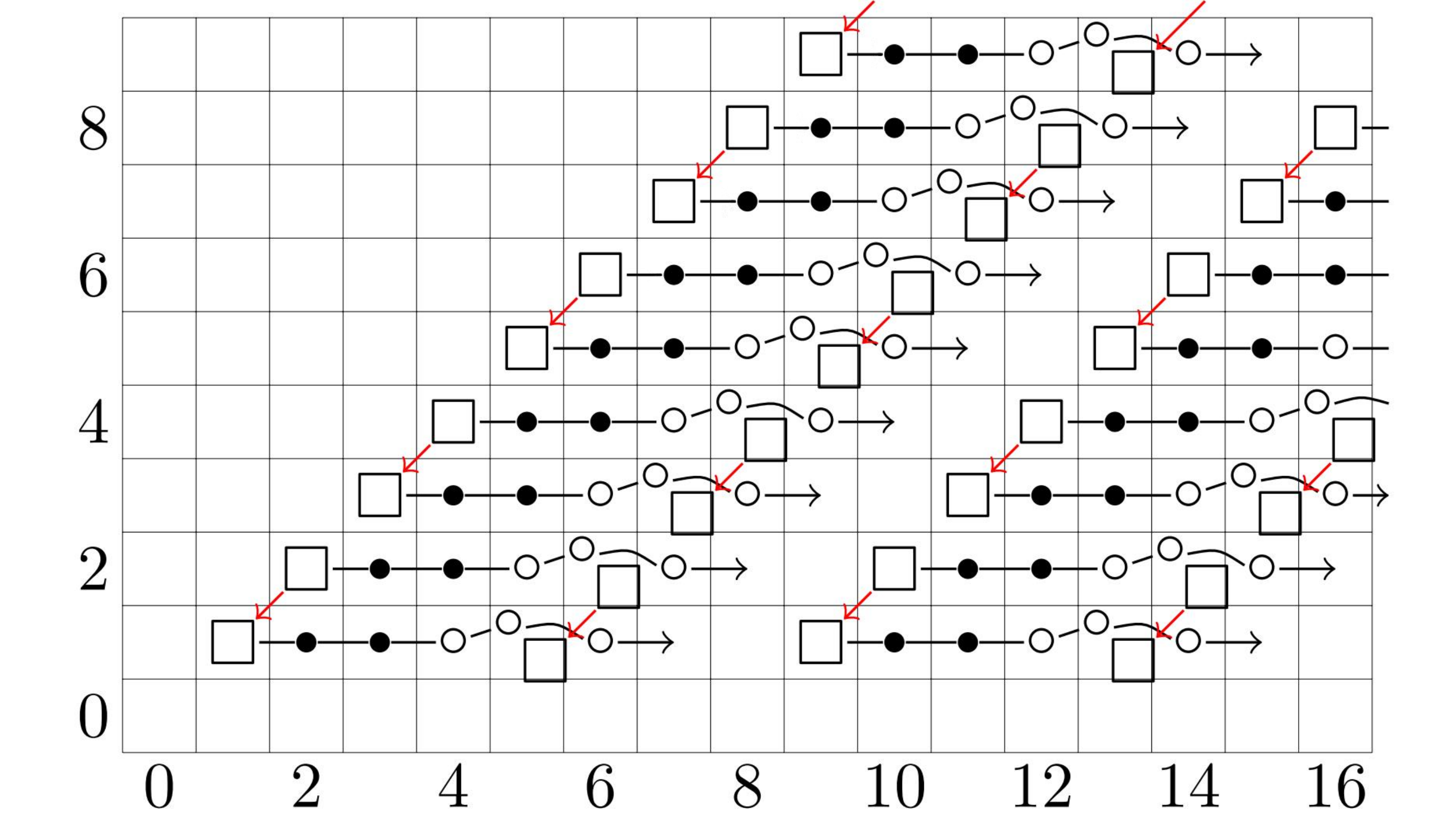}
\includegraphics[scale = 0.15]{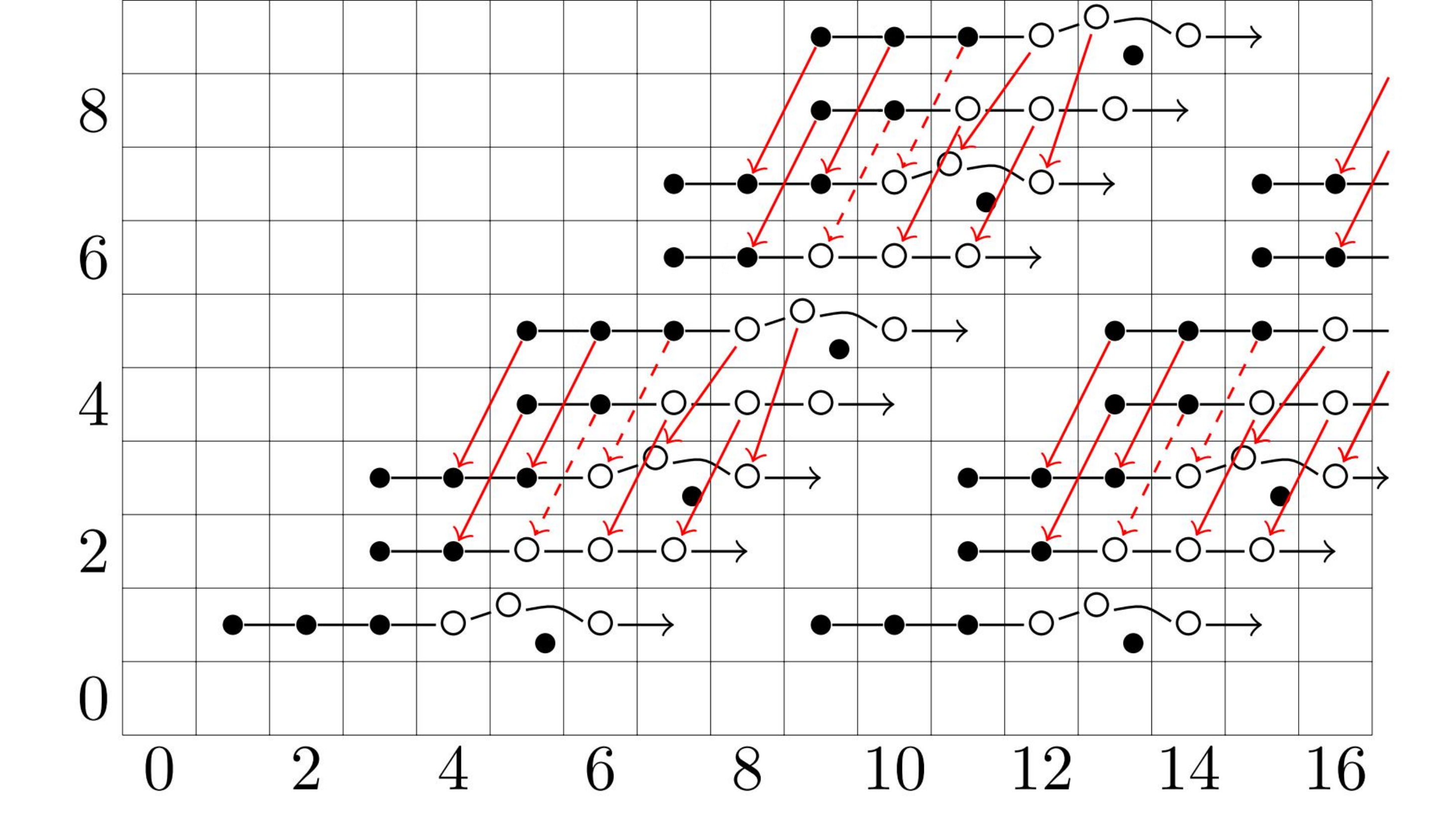}
\includegraphics[scale = 0.9]{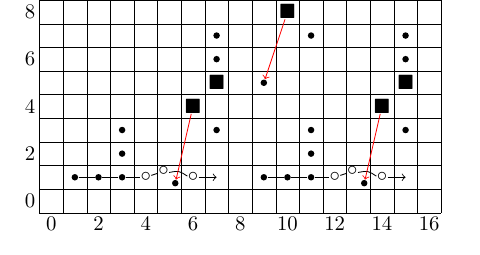}
\includegraphics[scale = 0.9]{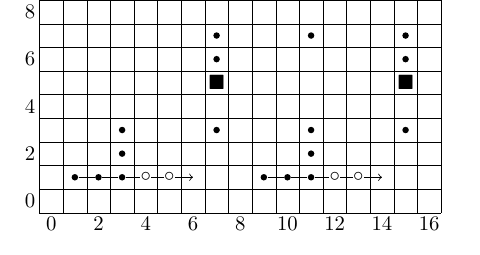}
\caption{Charts for the Atiyah--Hirzebruch spectral sequence converging to $kq_{**}BC_2$ ($E_1$ through $E_4 = E_{\infty}$ in sequential order). Here $\square = \z_2^{\wedge}[\tau]$, $\bullet=\f_2[\tau]$, $\circ=\f_2$, and $\blacksquare = \f_2[\tau]\{ \tau h_1^2 \}$. Horizontal lines represent multiplication by $\eta$. Differentials are in red. Note that the dashed arrows in the $E_2$-page means that the generators in $\f_2[\tau]$ multiplied with $\tau$ is mapped to 0, same as the generators themselves.}\label{fig:kqBC2_AH}.
\end{figure}

To solve extension problems, we use the motivic Adams spectral sequence, following the strategy described in \cref{kq_{**}BC_2}.

The connecting homomorphism, obtained by applying $\Ext_{\ca(1)}^{***}(-)$ to the short exact sequence of $\ca(1)$-modules \eqref{eqn: les for kq general strategy}, is determined by
$$\delta(\Sigma^{2,1} v_1^{2i}) = 
\begin{cases}
    b^{j/2} h_0 \quad & \text{ if } i=2j\geq 0, \\
    b^{j/2}a \quad & \text{ if } i=2j+1.
\end{cases}$$

On the other hand, we compute $\Ext_{\ca(1)}^{***}(R)$ as in \cref{kq_{**}BC_2}; namely,
$$\Ext_{\ca(1)}^{***}(R) \cong \bigoplus_{i \geq 0} \Sigma^{4i+3,2i+2} \Ext_{\ca(0)}^{***}(\m_2^{**}) \cong \bigoplus_{i \geq 0} \Sigma^{4i+3,2i+2} \m_2[h_0].$$

Returning to the long exact sequence of $\Ext$ derived from \eqref{eqn: les for kq with Q and R general strategy}, we see that in order for our final answer to agree with our Atiyah--Hirzebruch computations, the connecting homomorphism must be given by
$$\delta(\Sigma^{2,1}v_1^{2i+1}) = \Sigma^{4i+3,2i+2}h_0^{2i+2}$$
(and then extend by $h_0$-linearity). 
As in the classical case, there is a hidden $h_0$-extension, which can be resolved here using Betti realization.

\begin{figure}[H]
\includegraphics[scale = 0.15]
{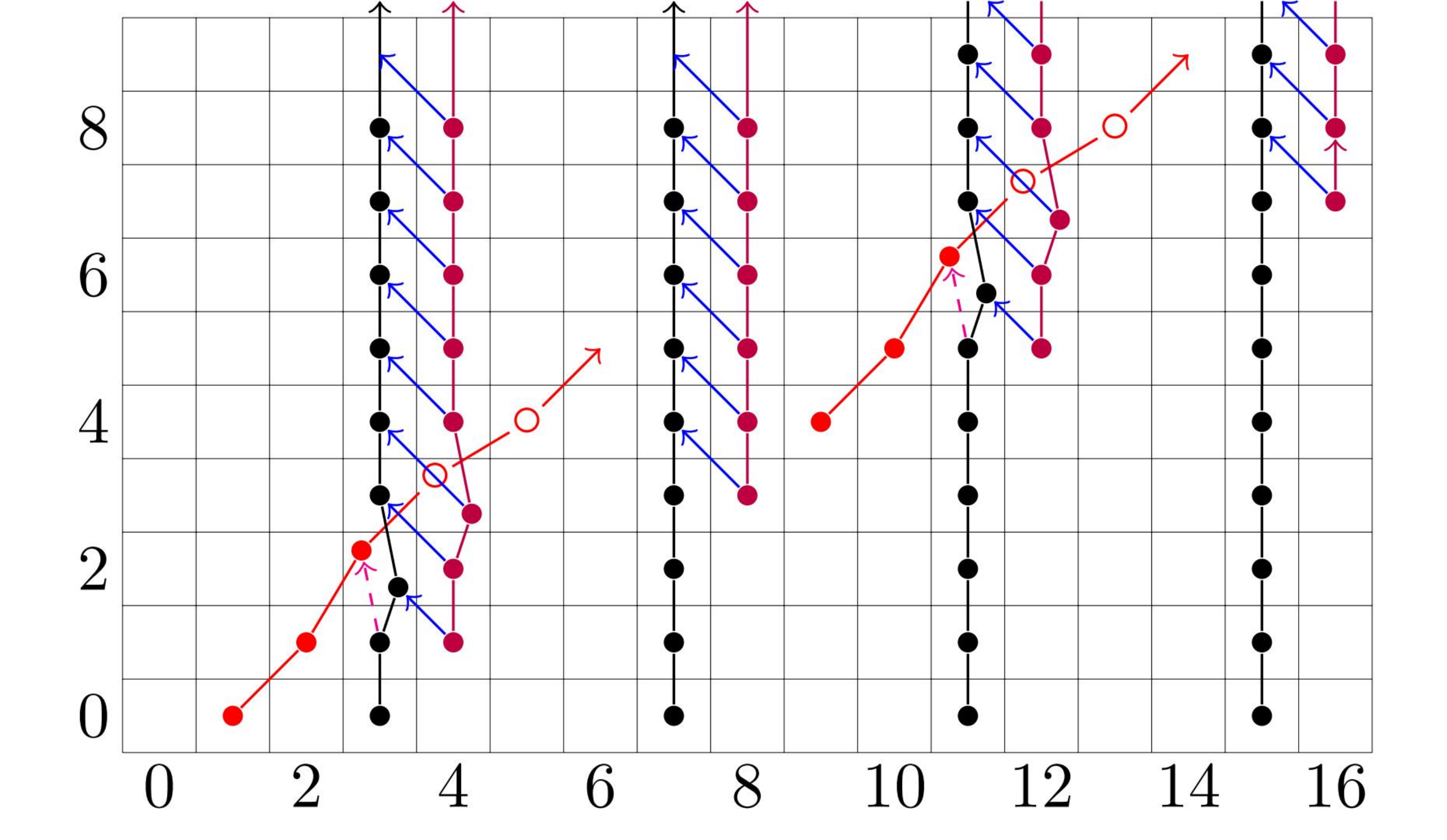}
\includegraphics{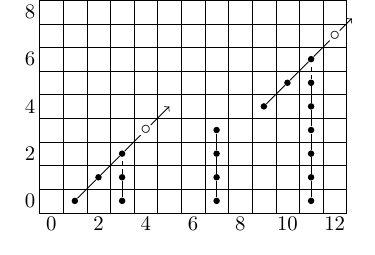}
\caption{The long exact sequence computing the $E_2$-page of the motivic Adams spectral sequence over $\c$ (left) and the resulting $E_\infty$-page (right). Here each  $\bullet$ represents $\m_2$, and each $\circ$ represents $\f_2$. Vertical lines represent multiplication by $h_0$, and the red diagonal lines represent multiplication by $h_1$. The red lines and purple lines represent the contribution from $Q$, and the black lines represent the contribution from $R$. The connecting homomorphisms are in blue. The dashed lines represent the hidden $h_0$-extensions.}\label{fig:kq_ASS_E2ALL}
\end{figure}

This completes the calculation of the $E_2$-term of the Adams spectral sequence. The spectral sequence collapses for degree reasons. The calculation and final $E_2$-term are depicted in \cref{fig:kq_ASS_E2ALL}.

We can now read off the final answer. As in the $kgl$ computations above, each instance of $\f_2[h_0]/(h_0^k)$ below detects $\z/2^k$. 

\begin{thm}
    The $E_\infty$-page of the motivic Adams spectral sequence converging to the very effective Hermitian K-theory of $BC_2$ over algebraically closed fields is given as a module over the $E_\infty$-page of the motivic Adams spectral sequence for $kq_{**}$ by
    \[
        E_\infty \cong \left( \bigoplus_{i \geq 0} \Sigma^{4i+3, 2i+2} \frac{\mathbb{F}_2[\tau, h_0]}{(h_0^{2i+2})} \right) \oplus \left ( \bigoplus_{j \geq 0} \Sigma^{8j+1, 4j+1} \frac{\f_2[\tau, h_1]}{(\tau h_1^3)} \right ).
    \]
\end{thm}

\begin{remark}
Here, and in all subsequent computations, we state the final answer in terms of the $E_\infty$-page of the Adams spectral sequence, instead of describing the exact bigraded abelian group as we did in the computation of $kgl_{**}(BC_2)$ above. This provides more information about the multiplicative structure (as a graded ring, and as a module over the coefficients of $kgl$ or $kq$), but it sometimes leaves open questions about hidden extensions. In most cases, there will be no room for hidden extensions and the bigraded abelian group can be extracted relatively easily from the $E_\infty$-term. For instance, each $\f_2[\tau,h_0]/(h_0^{2i+2})$ in the $E_\infty$-page described in the previous theorem gives rise to $\z/2^{2i+2}[\tau]$ in $kgl_{**}(BC_2)$. 
\end{remark}


\section{Computations over real closed fields}\label{Sec:R}

We now compute $H\z_{**}BC_2$, $kgl_{**}BC_2$, and $kq_{**}BC_2$ from $H^{**}BC_2$ over the real numbers. Our approach follows the same strategy as in the case over the complex numbers (or algebraically closed fields), as outlined in \cref{Sec:Strategy}, except that we now use the $\rho$-Bockstein spectral sequence \cite{Hil11} to compute $\r$-motivic $\Ext$ groups from $\c$-motivic $\Ext$ groups. The computations for $kgl$ and $kq$ are correspondingly more involved.

As in the complex case, we note that these computations yield identical results over any real closed field. To justify this, we begin by proving a real analogue of the “Lefschetz principle,” which makes the extension to real closed fields natural.

\subsection{Motivic real Lefschetz principle}

The Lefschetz principle, originally due to Solomon Lefschetz, asserts that algebraic geometry over algebraically closed fields is, in a suitable sense, equivalent to algebraic geometry over the complex numbers. Motivically, by \cite[Prop.~5.2.1]{BCQ25}, this means that working with motivic spectra over algebraically closed fields is equivalent to working over the complex numbers, as stated at the beginning of \cref{Sec:C}. We prove an analogous result comparing real closed fields and the real numbers here. The proof is parallel to that of \cite[Prop.~5.2.1]{BCQ25}, and part of the argument is inspired by \cite{WO17} and \cite{HO16}.

\begin{proposition} \label{prop: Motivic real Lefschetz principle}
    Let $F$ be a real closed field, and let $\Sp^{F}$ denote the category of modules over the $H \f_2^F$-nilpotent completion of the $F$-motivic sphere spectrum $S^{0,0}$. Let $\Sp^{F,\text{cell}} \subset \Sp^F$ denote the full subcategory of cellular spectra, i.e., the subcategory generated under colimits by the spheres $S^{a,b}$. Then there is an equivalence of categories
    \[
        \Sp^{F, \text{cell}} \simeq \Sp^{\r, \text{cell}}.
    \]
    Moreover, this equivalence is compatible with Adams spectral sequences, in the sense that it induces the isomorphism $\Ext_F \cong \Ext_{\r}$.
\end{proposition}
\begin{proof}
    Let $\r_{\text{alg}}$ denote the real closure of the rational numbers $\q$. This is the smallest real closed field. Consider the zigzag of field embeddings 
    \begin{equation}\label{eq:emb}
        F \xhookleftarrow{} \r_{\text{alg}} \xhookrightarrow{} \r.
    \end{equation}
    These field extensions give rise to base change functors
    \[
        \Sp^F \leftarrow \Sp^{\r_{\text{alg}}} \rightarrow \Sp^{\r}.
    \]
    Restricting to the cellular subcategories, we obtain
    \[
        \Sp^{F, \text{cell}} \xleftarrow{L_F} \Sp^{\r_{\text{alg}}, \text{cell}} \xrightarrow{L_\r} \Sp^{\r, \text{cell}}.
    \]
    It suffices to show that $L_F$ and $L_\r$ are equivalences. For simplicity, let $k \in {F, \r, \r_{\text{alg}}}$. Note that $\Sp^{k, \text{cell}}$ is stable and presentable, and is generated by the compact object $S_k^{0,0}$. It follows from the $\infty$-categorical Schwede--Shipley theorem \cite[Prop. 7.1.2.1, Rmk. 7.1.2.3]{Lur16} that
    \[
        \Sp^{k, \text{cell}} \simeq \text{RMod}_{\text{End}(S_k^{0,0})}.
    \]
    Since $\text{End}(S_k^{0,0}) \cong \pi_{*,*}(S_k^{0,0})$, it suffices to show that $L_F$ and $L_{\r}$ induce isomorphisms on the motivic stable stems $\pi_{*,*}(S_k^{0,0})$. 
    
    By the assumption that everything is $H\f_2^F$-nilpotent complete, the motivic Adams spectral sequence over $k$, with
    \[
        E_2 = \Ext_{\ca^k}^{***}(H^{**}(k; \f_2))
    \]
    converges strongly and unconditionally to $\pi_{*,*}(S_k^{0,0})_2^{\wedge}$. By \cite{Voe11},
    \[
        H^{**}(k; \f_2) = \text{K}_*^{M}(k)/2 \,\, [\tau],
    \]
    where $\text{K}_*^{M}(k)$ denotes the Milnor K-theory of $k$ and $\tau \in H^{0,1}(k; \f_2)$. By Artin--Schreier, any real closed field $k$ has absolute Galois group $\text{Gal}(k^{\text{sep}} \mid k) = C_2$. Thus, the mod two Galois cohomology of $k$ is 
    \[
        H_{\text{Gal}}^{*}(k; \f_2) \cong H^{*}(C_2; \f_2) \cong \f_2[x],
    \]
    with $|x|=1$, and the zigzag \eqref{eq:emb} gives rise to isomorphisms of Galois cohomology rings
    $$H^*_{\text{Gal}}(F; \f_2) \xleftarrow{\cong} H^*_{\text{Gal}}(\r_{\text{alg}};\f_2) \xrightarrow{\cong} H^*_{\text{Gal}}(\r;\f_2).$$
    Voevodsky's proof of the Milnor conjecture identifies Milnor K-theory with Galois cohomology,
    \[
        \text{K}_*^{M}(k)/2 \cong H_{\text{Gal}}^{*}(k; \f_2),
    \]
    sending $\{ -1 \}$ to $x$. Under this identification, $\{ -1 \}$ corresponds to $\rho \in H^{1,1}$ in  motivic cohomology. Hence, 
    \[
        (\m_2^k)^{**} =  H^{**}(k; \f_2) = \f_2[\rho, \tau]
    \]
    for all $k \in \{ F, \r, \r_{\text{alg}} \}$, with isomorphisms between these rings as $k$ varies induced by \eqref{eq:emb}.  By \cite[Rmk. 2.2.8]{BCQ25}, the map
    \[
        \ca^{\r} \, \otimes_{(\m_2^{\r})^{**}} (\m_2^{k})^{**} \to \ca^{k}
    \]
    is an isomorphism, which yields
    \[
        \Ext_k \cong \Ext_{\ca^\r}((\m_2^{\r})^{**}, (\m_2^{k})^{**}).
    \]
    Therefore, $L_F$ and $L_{\r}$ induce isomorphisms on the $E_2$-pages of the corresponding motivic Adams spectral sequences. Moreover, the base change functors $L_F$ and $L_{\r}$ preserve the Adams filtrations, as they are symmetric monoidal, exact, and left adjoints. Since these spectral sequences converge strongly, it follows that $L_F$ and $L_{\r}$ induce isomorphisms on the motivic stable stems. The result follows.
\end{proof}

\subsection{Mod two motivic cohomology}

We begin the computational process outlined in \cref{Sec:Strategy} over the real numbers $\r$. By the motivic real Lefschetz principle (\cref{prop: Motivic real Lefschetz principle}), these computations yield identical results over any real closed field.

We have $H\f_2^{**} = (\m_2^\r)^{**} = \f_2[\tau,\rho]$, so 
$$H^{**}(BC_2) \cong \ker\left( \f_2[\tau,\rho][x,y]/(x^2 + \tau y + \rho x) \to \f_2[\tau,\rho], \quad x,y \mapsto 0 \right)$$
where $|x| = (1,1)$ and $|y| = (2,1)$. 

\subsection{Integral motivic homology}

By \cref{Prop:HZBC2}, we have
$$H\z_{**}(BC_2) \cong \bigoplus_{i \geq 0} \Sigma^{2i+1,i+1}\m_2 \cong \bigoplus_{i \geq 0} \Sigma^{2i+1,i+1}\f_2[\tau,\rho].$$

\subsection{Effective algebraic K-theory}

We now compute $kgl_{**}BC_2$. 

We begin with $kgl_{**}$, which can be computed by first analyzing the $\rho$-Bockstein spectral sequence and then analyzing the motivic Adams spectral sequence:
$$E_1 = \Ext_{\ce(1)^\c}^{***}((\m_2^\c)^{**})[\rho] \overset{\rho-BSS}{\Longrightarrow} \Ext_{\ce(1)}^{***}((\m_2^\r)^{**}) \overset{mASS}{\Longrightarrow} kgl_{**}.$$

In the $\rho$-Bockstein spectral sequence, the nontrivial differentials are shown by \cite[Theorem 3.2]{Hil11}
$$d_1(\tau) = \rho h_0, \quad d_3(\tau^2) = \rho^3 v_1.$$
There are no hidden extensions (\cite[Theorem 3.4]{Hil11}), leaving 
\begin{equation} \label{eqn: result for Ext_E1(M2) in R}
    \Ext_{\ce(1)}^{***}(\m_2^{**}) \cong \dfrac{\f_2[\tau^4,\rho,h_0,v_1]\{1,[\tau^2 h_0]\}}{(\rho h_0, \rho^3 v_1)}.
\end{equation}

The motivic Adams spectral sequence also collapses since there is no room for nontrivial differentials from the multiplicative generators, giving \cite[Fig. 1]{Hil11}.

\begin{figure}[H]
    \includegraphics[scale=0.5]{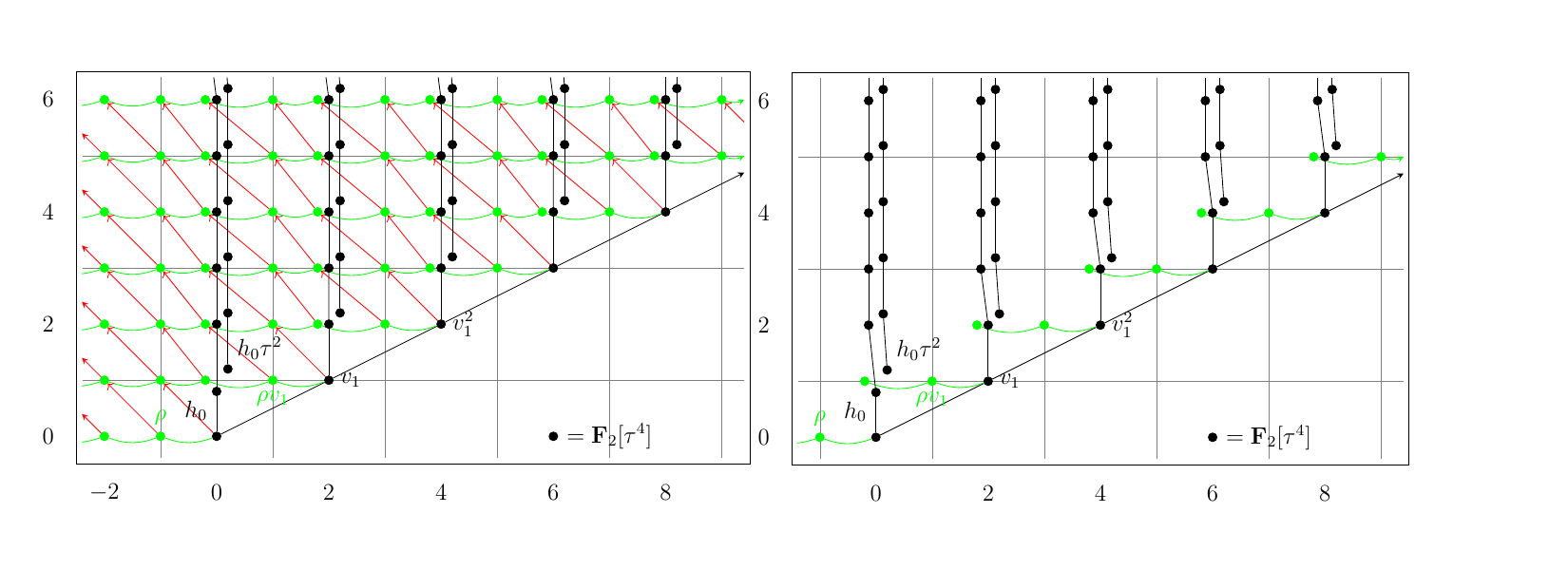}
    \vspace{-0.5cm}
    \caption{The left chart is the $E_3$-page of the $\rho$-Bockstein spectral sequence. Each $\bullet$ represents $\f_2[\tau^4]$. Differentials are in red.
    The right chart is $E_2=E_{\infty}$-page of the motivic Adams spectral sequence converging to $kgl_{**}$.} 
    \label{Fig. 1} 
\end{figure}

We now analyze the Atiyah--Hirzebruch spectral sequence. The $d_1$-differentials are determined by
$$d_1(x[2n]) = 2x[2n-1]$$
for $n \geq 0$. Unlike in the classical and $\c$-motivic cases, the spectral sequence does not collapse at $E_2$. Instead, we have $E_2 = E_3$, and the $d_3$-differentials are given by
$$d_3(x[2n]) = v_1 x[2n-3]$$
for all $n \geq 2$. These differentials follow from the observation that there is a nontrivial action of $Q_1 = [Sq^2,Sq^1]$ in degrees $2n-3$ of $H^{**}(BC_2)$, together with the fact that $Q_1$ detects multiplication by $v_1$ in $kgl_{**}$. After running these differentials, the spectral sequence collapses at $E_4$. 
\begin{figure}[H]
    \includegraphics[trim={0 0 1cm 0},clip,scale = 1.05]{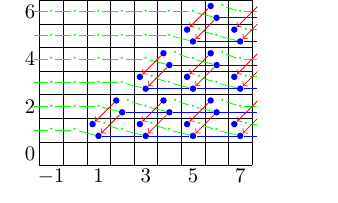} 
    \includegraphics[trim={0 0 1cm 0},clip,scale = 1.05]{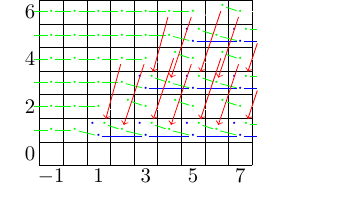}
    \includegraphics[trim={0 0 1cm 0},clip,scale = 1.05]{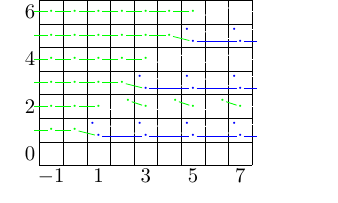}
    \vspace{-0.5cm}
    \caption{Charts for the Atiyah--Hirzebruch spectral sequence converging to $kgl_{**}BC_2$ ($E_1, E_3,$ and $E_4= E_{\infty}$ in sequential order). Here $\bullet = \z_2[\tau^4]$ and $\cdot=\f_2[\tau]$.}\label{fig:kglBC2R}.
\end{figure}
\vspace{-0.5cm}
We follow \cref{SS:kgl_general} to resolve hidden extensions. In the long exact sequence of $\Ext$ derived from \eqref{eqn: les for kgl general strategy}, the first term has already been determined in \eqref{eqn: result for Ext_E1(M2) in R}. To compute the result in \eqref{eqn: E1-page for kgl general strategy}, we still need to calculate $\Ext_{\ca(0)^{\mathbb{R}}}^{***}((\mathbb{M}_2^{\mathbb{R}})^{**})$, which can be done via the $\rho$-Bockstein spectral sequence:
\[
    E_1 = \Ext_{\ca(0)^{\c}}^{***}((\mathbb{M}_2^\c)^{**})[\rho] \Rightarrow \Ext_{\ca(0)^{\mathbb{R}}}^{***}((\mathbb{M}_2^\r)^{**}).
\]
By \cite{Hil11}, the only nontrivial differential is generated by 
\[
    d_1(\tau) = \rho h_0,
\]
and there are no hidden extensions. It follows that
\[
    \Ext_{\ca(0)^{\mathbb{R}}}^{***}((\mathbb{M}_2^\r)^{**}) \cong \frac{\f_2[\tau^2, \rho, h_0]}{(\rho h_0)}.
\]

Substituting this result into \eqref{eqn: E1-page for kgl general strategy}, we obtain
\begin{equation*}
    \Ext_{\ce(1)}^{***}(V) \cong \bigoplus_{i \geq 0} \Sigma^{2i+1,i+1} \Ext_{\ca(0)^{\r}}^{***}((\m_2^\r)^{**}) \cong \bigoplus_{i \geq 0} \Sigma^{2i+1,i+1} \dfrac{\f_2[\tau^2,\rho,h_0]}{(\rho h_0)}.
\end{equation*}
Comparing to the ranks of $gr^*_{AH}kgl_{**}BC_2$, we see that the connecting homomorphism must be nontrivial.

\begin{figure}
\includegraphics[scale=0.7]{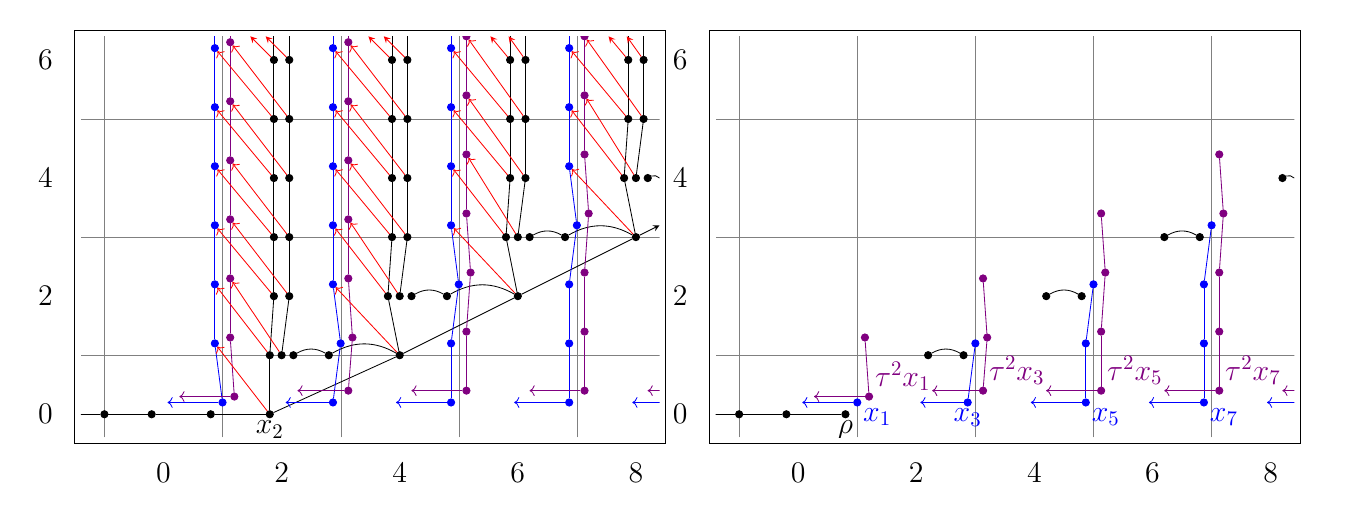}
    \caption{Charts for $kgl_{**}^{\mathbb{R}}(BC_2)$,
    where each $\bullet$ represents $\f_2[\tau^4]$. The black portion of the chart depicts $\Ext^{***}_{\ce(1)}(\Sigma^{2,1}\m_2)$, while the blue and purple portions depict $\Ext^{***}_{\ce(1)}(V)$.}
    \label{fig:ext_ce(1)(H_BC2) over R} 
\end{figure}

$\Ext^{***}_{\ce(1)}(H^{**}BC_2)$ forms the $E_2$-page of the motivic Adams spectral sequence converging to $kgl^{\r}_{**}(BC_2)$.
The spectral sequence collapses at the $E_2$-page for degree reasons, leaving no room for hidden extensions illustrated in \cref{fig:ext_ce(1)(H_BC2) over R}.

\begin{thm}
The $E_\infty$-page of the motivic Adams spectral sequence converging to the effective algebraic \(K\)-theory of \(BC_2\), as a module over the \(E_\infty\)-page of the motivic Adams spectral sequence converging to \(kgl_{**}\), is given by
\[
E_\infty \cong 
\frac{\f_2[\tau^4,\rho,v_1]}{(\rho^2 v_1)} \oplus \left ( \bigoplus_{i \geq 0} \frac{\f_2[\tau^4,\rho,h_0]\{x_{2i+1}\}}{(\rho h_0, h_0^{i+1})} \right )
 \oplus \left ( \bigoplus_{i \geq 0} \frac{\f_2[\tau^4,\rho, h_0]\{\tau^2 x_{2i+1}\}}{(\rho h_0, h_0^{i+2})} \right ).
\]
\end{thm}

\subsection{Very effective Hermitian K-theory}

We now turn to the case $kq_{**}BC_2$. 

\begin{figure}[H]
\includegraphics[scale=0.45]{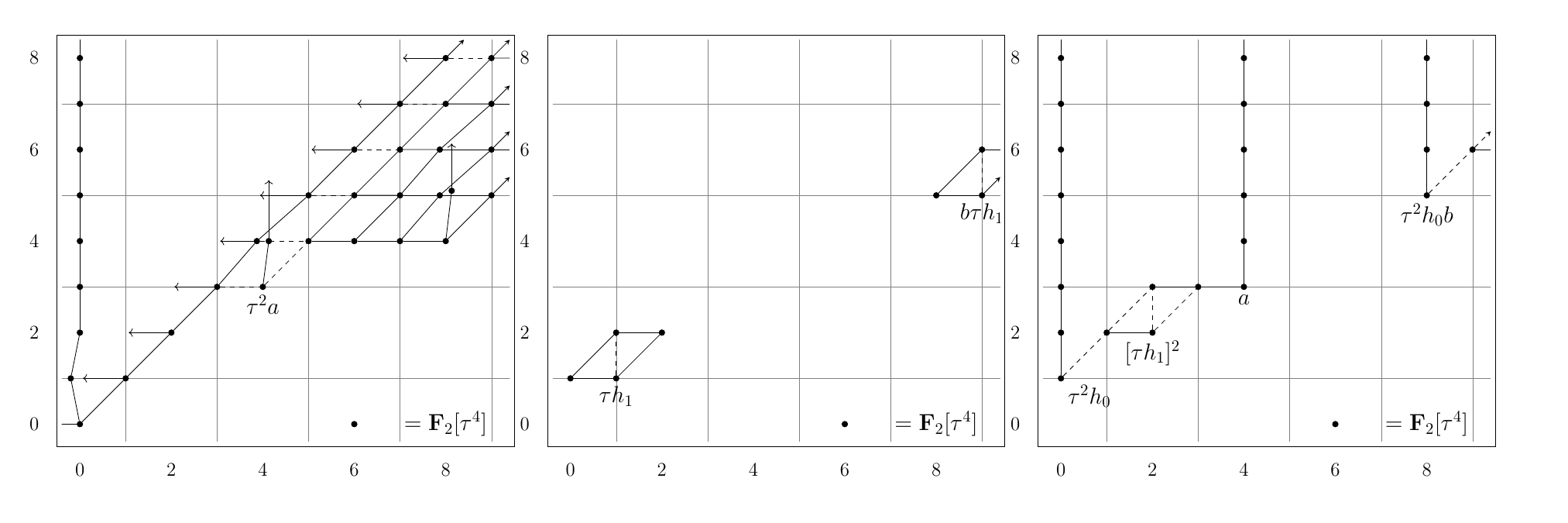}
\vspace{-1cm}
\caption{$\Ext_{\ca(1)}^{***}(\m_2^\r)$ organized by coweights modulo $4$, converging to $kq_{**}$. Each $\bullet$ represents $\f_2[\tau^4]$. The charts represent the case $0, 1, 2\,  ( \text{mod  } 4)$ from left to right, respectively. Note that in the case $3 \, (\text{mod } 4)$ the $\Ext$-term is trivial because no generators live in that coweight (see \cite[Theorem 6.2]{GHIR19}). We reproduce the charts from \cite[Section 12.1]{GHIR19}, which are extracted from \cite{Hil11}.}\label{fig:kqR}
\end{figure}

In \cite[Theorem 3.2]{Hil11}, Hill analyzed the $\rho$-Bockstein spectral sequence
$$E_1 = \Ext_{\ca(1)^{\c}}^{***}((\m_2^\c)^{**})[\rho] \Rightarrow \Ext_{\ca(1)^{\r}}^{***}((\m_2^\r)^{**}),$$
which computes the $E_2$-term of the motivic Adams spectral sequence 
$$E_2 = \Ext_{\ca(1)^{\r}}^{***}((\m_2^\r)^{**}) \Rightarrow kq_{**}.$$
As Hill explains, the motivic Adams spectral sequence collapses, yielding a complete explicit computation of $kq_{**}$ (see \cite[Theorem 4.8]{Hil11}). Since the ring structure of $kq_{**}$ is complicated and especially unwieldy for computational purposes, we follow the description in \cite[Section 12.1]{GHIR19} and use a chart as well as a theorem to present the final answer for $kq_{**}$. The $\Ext$ groups (and thus $kq_{**}$) are depicted in \cref{fig:kqR}. Here, ``coweight'' is defined as stem minus weight, i.e., $\text{coweight} = \text{stem} - \text{weight}$.

We now analyze the Atiyah--Hirzebruch spectral sequence
$$E_1 = \bigoplus_{i \geq 0} \Sigma^{2i+1,i+1} kq_{**} \oplus \Sigma^{2i+2,2i+1} kq_{**} \Rightarrow kq_{**}BC_2.$$
Rather than display the entire computation in one figure, we will take advantage of the periodicity in the cell structure of $BC_2$. 

\begin{figure}[H]
\includegraphics[scale=0.4]{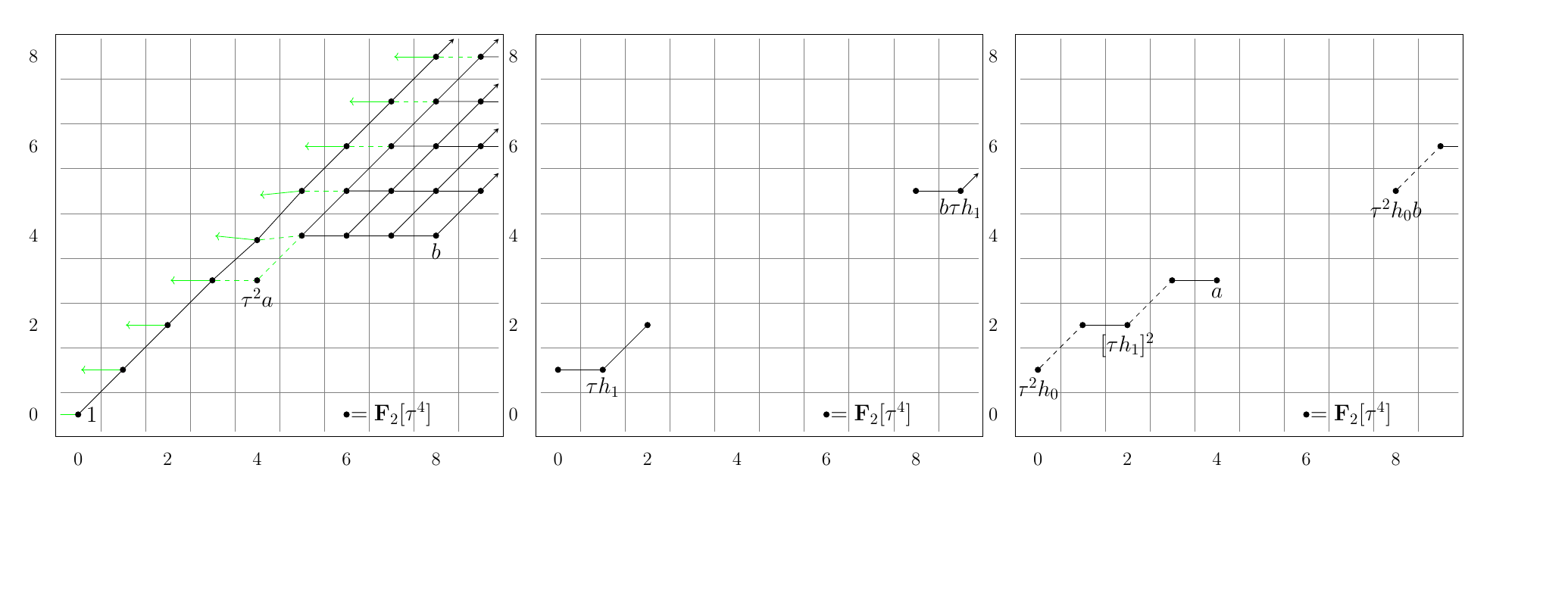}
\vspace{-1cm}
\caption{Charts representing the odd filtrations of the $E_2$-term of the Atiyah--Hirzebruch spectral sequence converging to $kq_{**}BC_2$. 
Each chart corresponds to a coweight modulo $4$: the left chart has coweight $0$, the middle $1$, and the right $2$.}\label{fig:kqAH_E2_f0}
\end{figure}

\begin{figure}
\includegraphics[scale=0.4]{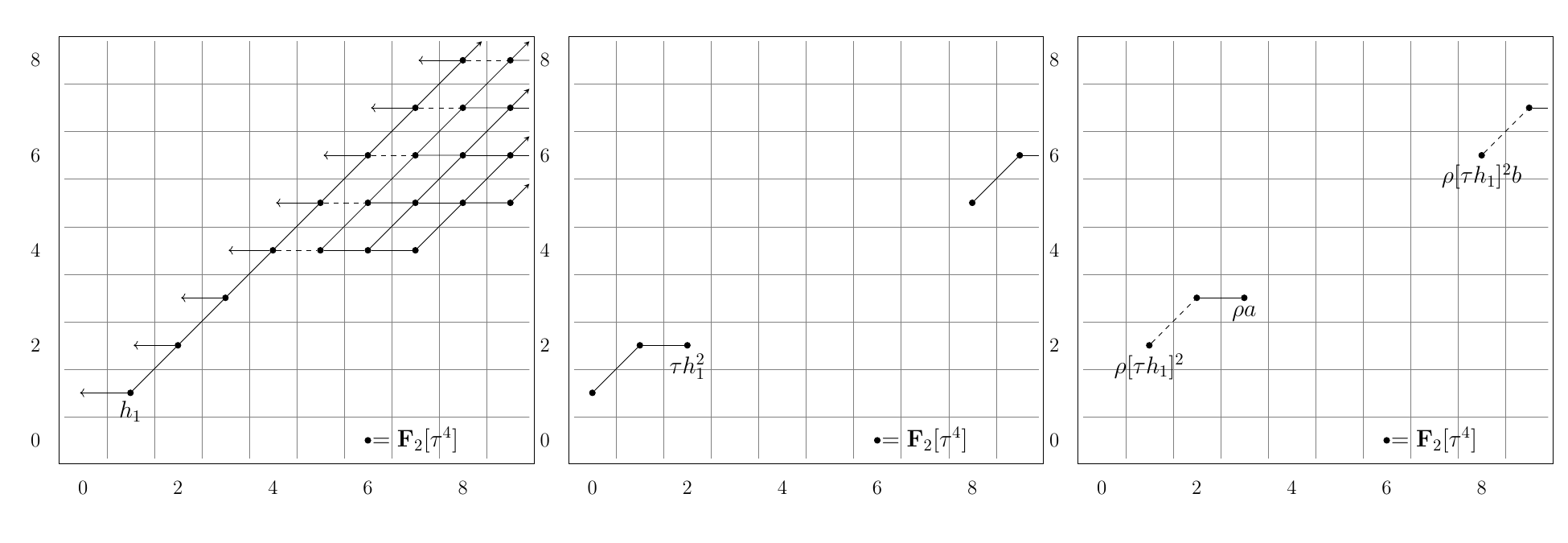}
\caption{Charts representing the even filtrations of the $E_2$-term of the Atiyah--Hirzebruch spectral sequence converging to $kq_{**}BC_2$. 
Each chart corresponds to a coweight modulo $4$: the left chart has coweight $0$, the middle $1$, and the right $2$.}\label{fig:kqAH_E2_f1}
\end{figure}

On the $E_1$-page, each positive filtration is isomorphic (up to a shift), so $kq_{**}$ appears in each positive filtration row of the spectral sequence. The $d_1$-differentials are given by multiplication by $h = [h_0]$ from even to odd filtrations, so on $E_2$, we have
$$E_2^{i,**} \cong
\begin{cases}
    kq_{**}/h \quad & \text{ if } i = 2k+1 > 0, \\
    \ker(h: kq_{**} \to kq_{**}) \quad & \text{ if } i = 2k > 0.
\end{cases}$$
The charts representing each row of the $E_2$-page can be obtained from the charts for $kq_{**}$ by erasing classes which are divisible by $h$ (if $i=2k+1$) or by erasing classes which support nontrivial multiplication by $h$ (if $i = 2k+2$). The charts for the $E_2$-term appear in \cref{fig:kqAH_E2_f0} and \cref{fig:kqAH_E2_f1}. 

To obtain the $E_3$-page, we use that the $d_2$-differentials are given by multiplication by $\eta = [h_1]$ from filtrations $4k$ and $4k+1$ to filtrations $4k-2$ and $4k-1$, respectively. As the differentials do not change the parity of the filtration and the rows in the $E_2$-term only depend on the parity of the filtration, the $E_3$-page may be obtained by computing the kernel and cokernel of multiplication by $\eta$ on each row of the $E_2$-page. More precisely, we have
$$E_3^{i,**} \cong 
\begin{cases}
    \ker(\eta: kq_{**}/h \to kq_{**}/h) \quad & \text{ if } i = 4k+1>4, \\
    \ker(h : kq_{**} \to kq_{**})/\eta \quad & \text{ if } i = 4k+2>0, \\
    (kq_{**}/h)/\eta \quad & \text{ if } i = 4k+3, \\
    \ker(\eta: \ker(h) \to \ker(h)) \quad & \text{ if } i = 4k+4 > 0.
\end{cases}$$
In particular, charts representing the rows of the $E_3$-term can be obtained from the charts representing the rows of the $E_2$-term by erasing classes which are $\eta$-divisible or which support nontrivial multiplication by $\eta$, depending on the filtration of the row. 

\begin{figure}
\includegraphics[scale=0.35]{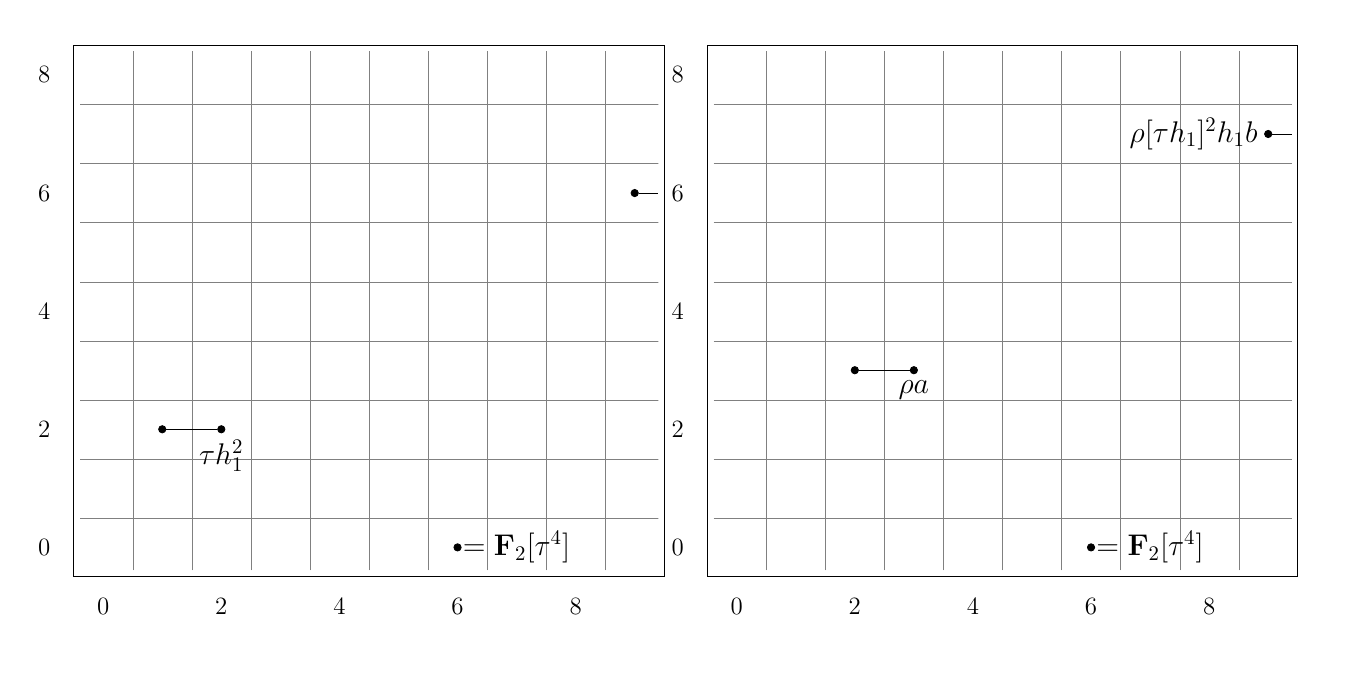}
\vspace{-0.5cm}
\caption{Charts representing filtrations congruent to $0$ mod $4$ of the $E_3$-page of the Atiyah--Hirzebruch spectral sequence converging to $kq_{**}BC_2$. Only coweights congruent to $1$ (on the left) and $2$ (on the right) are displayed since the other coweights vanish. }\label{fig:kqAH_E3_f0}

\includegraphics[scale=0.35]{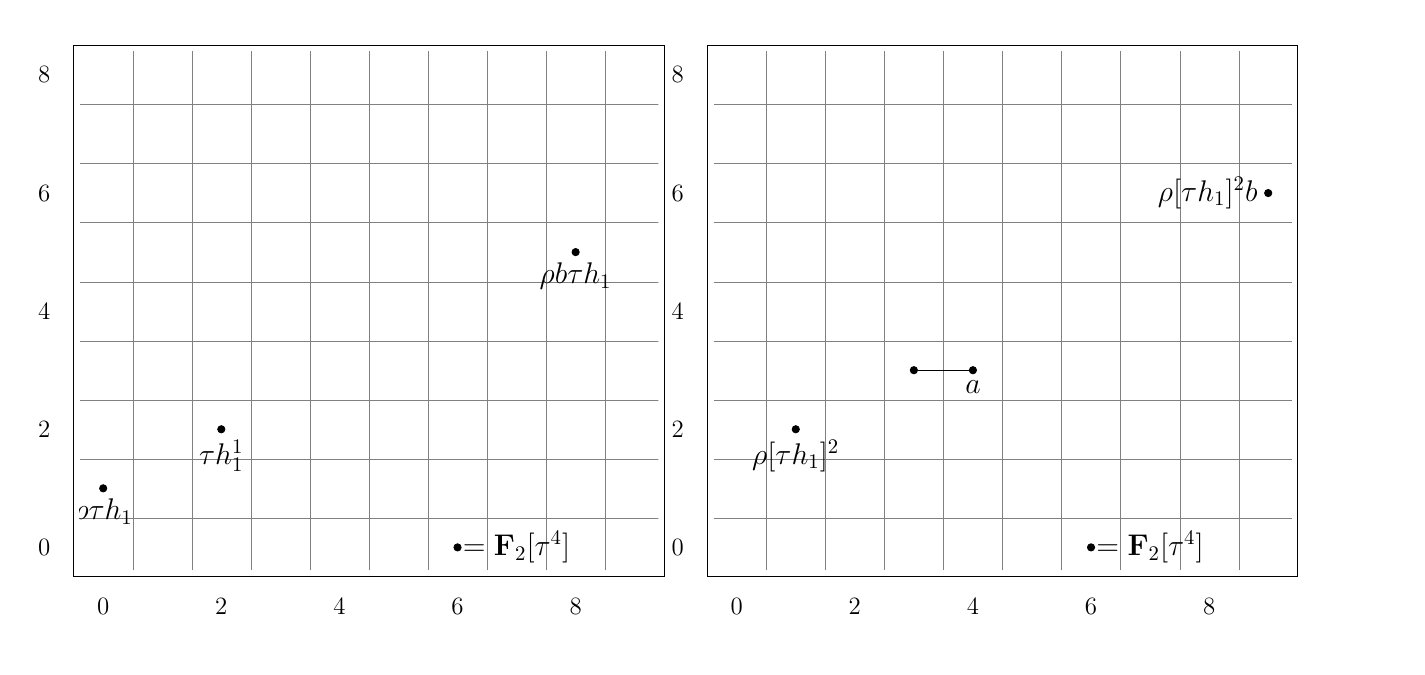}
\vspace{-0.5cm}
\caption{Charts representing filtrations congruent to $1$ mod $4$ of the $E_3$-page of the Atiyah--Hirzebruch spectral sequence converging to $kq_{**}BC_2$. Only coweights congruent to $1$ (on the left) and $2$ (on the right) are displayed since the other coweights vanish. }\label{fig:kqAH_E3_f1}

\includegraphics[scale=0.35]{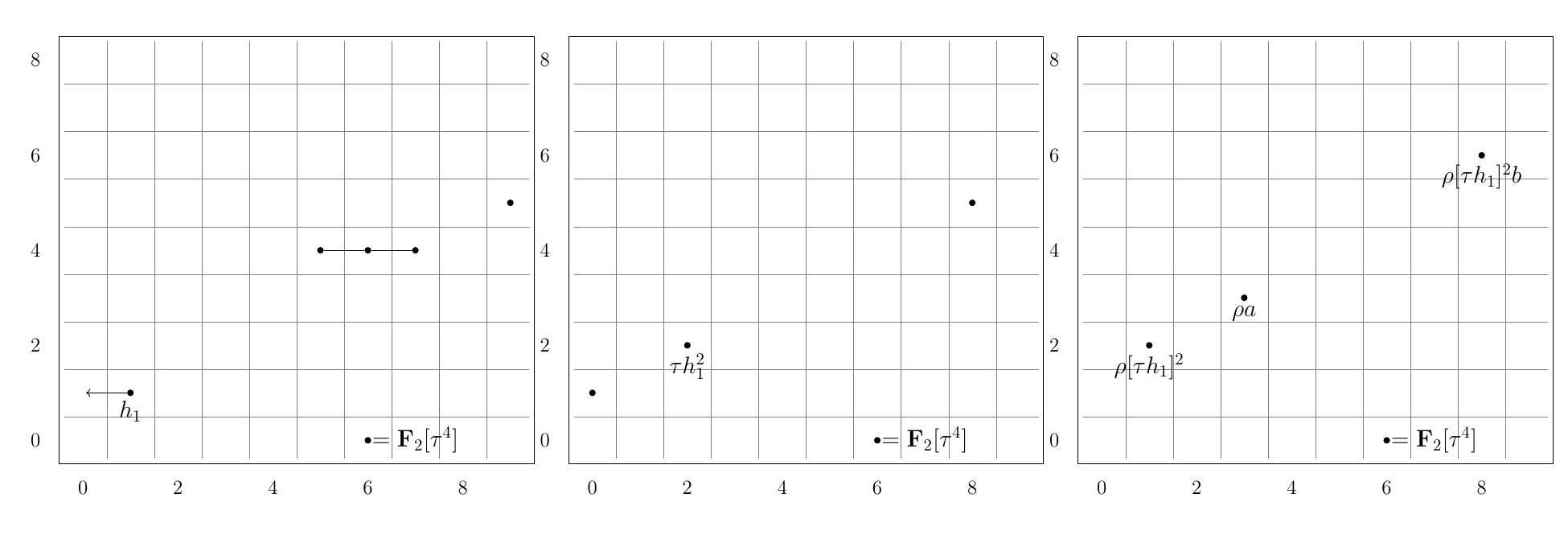}\vspace{-0.5cm}
\caption{Charts representing filtrations congruent to $2$ mod $4$ of the $E_3$-page of the Atiyah--Hirzebruch spectral sequence converging to $kq_{**}BC_2$.}\label{fig:kqAH_E3_f2}

\includegraphics[scale=0.35]{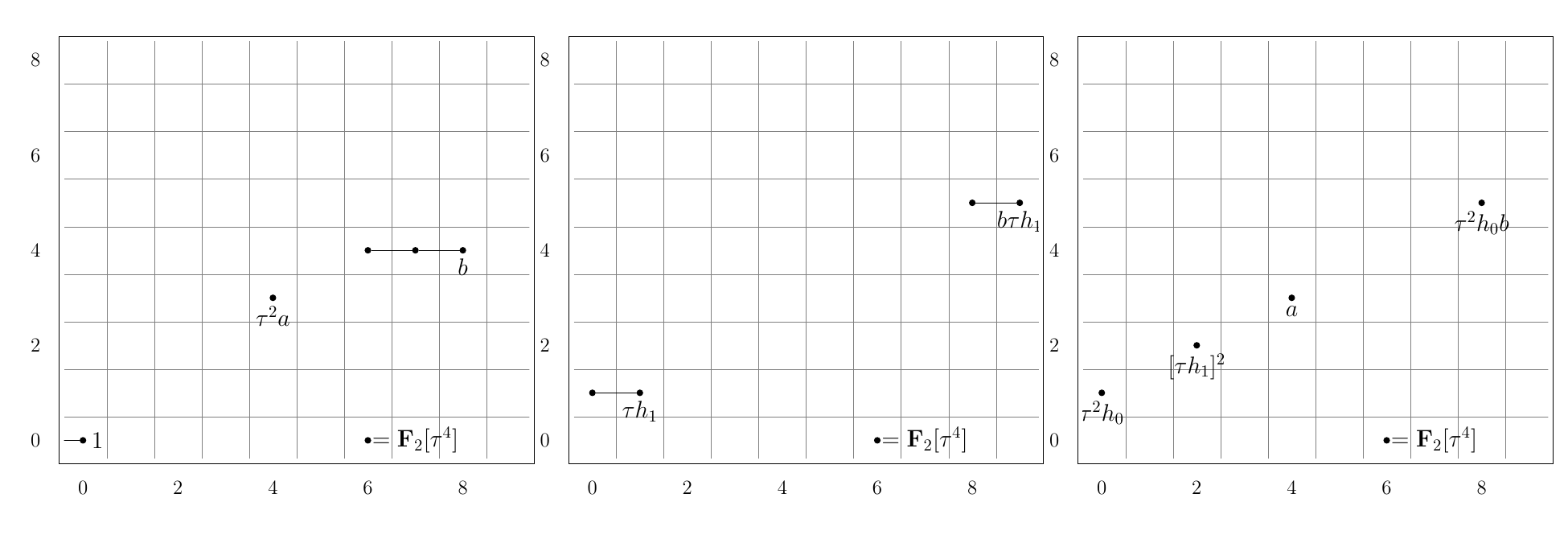}\vspace{-0.5cm}
\caption{Charts representing filtrations congruent to $3$ mod $4$ of the $E_3$-page of the Atiyah--Hirzebruch spectral sequence converging to $kq_{**}BC_2$.}\label{fig:kqAH_E3_f3}
\end{figure}

The chart in filtration $1$ is the same as for the $E_2$-page since there are no differentials entering or exiting that row. Charts representing each congruence class of filtration modulo $4$ appear in \cref{fig:kqAH_E3_f0}, \cref{fig:kqAH_E3_f1}, \cref{fig:kqAH_E3_f2}, and \cref{fig:kqAH_E3_f3}.

To obtain the $E_4$-page, we use that the $d_3$-differentials are given by Toda brackets of the form $\langle x , \eta, h \rangle$ if $x$ has filtration $i = 4k>0$, or $\langle x, h, \eta \rangle$ if $x$ has filtration $i = 4k+2 \geq 6$. Most of the relevant Toda brackets can be read off of \cite[Table 6]{GHIR19}. 

\vspace{0.5cm}
To determine the behavior of the $d_3$-differentials, we note that, by definition, $d_3$ should

\begin{enumerate}
    \item decrease the filtration from $i$ to $i-3$ (modulo 4);
    \item increase the stem $t-s$ to $t-s+2$;
    \item increase the coweight $w$ to $w+1$ (modulo 4).
\end{enumerate}

Explicitly, we have the following:

\begin{itemize}
    \item There are $d_3$-differentials from \cref{fig:kqAH_E3_f0} to \cref{fig:kqAH_E3_f1}. In particular, $d_3$ maps the left chart (coweight $1 \bmod 4$) of \cref{fig:kqAH_E3_f0} to the right chart (coweight $2 \bmod 4$) of \cref{fig:kqAH_E3_f1}, sending the dots $(1,2)$ (i.e., $\rho h_1 \cdot \tau h_1$) and $(2,2)$ (i.e., $h_1 \cdot \tau h_1$) in the source to the dots $(3,3)$ (i.e., $\rho a$) and $(4,3)$ (i.e., $a$) in the target, respectively. Similar patterns occur for the other dots in the source and target. These differentials originate from $a \in \langle \tau h_1 \cdot h_1, h_1, h_0 \rangle$ and satisfy $d_3(\rho \tau h_1 \cdot h_1) = \rho a$.

    \item There are $d_3$-differentials from \cref{fig:kqAH_E3_f2} to \cref{fig:kqAH_E3_f3}. In particular, $d_3$ maps the first chart (coweight $0 \bmod 4$) of \cref{fig:kqAH_E3_f2} to the second chart (coweight $1 \bmod 4$) of \cref{fig:kqAH_E3_f3}, sending the dot $(-1,0)$ (i.e., $\rho$) in the source to the dot $(1,1)$ (i.e., $\tau h_1$) in the target. It also maps the second chart (coweight $1 \bmod 4$) of \cref{fig:kqAH_E3_f2} to the third chart (coweight $2 \bmod 4$) of \cref{fig:kqAH_E3_f3}, sending the dot $(0,1)$ (i.e., $\rho \tau h_1$) in the source to the dot $(2,2)$ (i.e., $(\tau h_1)^2$) in the target. Similar patterns occur for the other dots in the source and target. These differentials originate from $(\tau h_1)^2 \in \langle \rho \tau h_1, h_0, h_1 \rangle$, obtained by juggling from $\tau h_1 \in \langle \rho, h_0, h_1 \rangle$ and using the linearity of $d_3$ with respect to $b$.

    \item There are no $d_3$-differentials out of filtrations $1$ or $3 \bmod 4$, because there are no primary attaching maps from filtration $1 \bmod 4$ to $2 \bmod 4$ or from filtration $3 \bmod 4$ to $0 \bmod 4$.
\end{itemize}

\begin{figure}[H]
\includegraphics{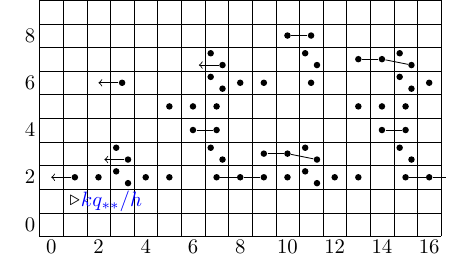}
    \caption{Chart representing the whole $E_4$-page of the Atiyah--Hirzebruch spectral sequence converging to $kq_{**}BC_2$. Note that the vertical axis represents the filtration modulo 4, and the horizontal line adorned with a right blue triangle indicates the portion corresponding to $kq_{**}/h$. }
\label{fig:kgl_AHSS_BC2_E4all}
\end{figure}

By assembling these charts after taking the $d_3$-differentials, we obtain the complete $E_4$-page appearing in \cref{fig:kgl_AHSS_BC2_E4all}. 

Note that the $d_4$-differentials decrease by four units in the vertical direction and by one unit in the horizontal direction. The only nontrivial possibility is from $(8,6)$ to $(7,2)$ in \cref{fig:kgl_AHSS_BC2_E4all}. However, the $d_4$-differential increases the coweight by two, thereby ruling out the above possibility for degree reasons. Thus, the Atiyah--Hirzebruch spectral sequence collapses at $E_4$, giving us a full calculation of $gr^{*}_{AH}kq_{**}BC_2$. The result can be read off from \cref{fig:kgl_AHSS_BC2_E4all} as follows:
\[
    gr^n_{AH}kq_{s,*}BC_2 =
    \begin{cases}
        kq_{s,*}/h & \text{if } n = 1; \\
        \f_2^{\oplus (1 + [s \equiv 3 \text{ mod } 8] - [s \equiv 6 \text{ mod } 8])}  & \text{if } n = 2, s \in \n; \\
        \f_2^{\oplus [s - n \equiv 2, 3 \text{ mod } 8]} & \text{if } n \equiv 0 \text{ mod } 4, n \geq 4, s \geq n+2; \\
        \f_2^{\oplus [s - n \equiv 0, 1, 2 \text{ mod } 8]} & \text{if } n \equiv 1 \text{ mod } 4, \, n \geq 5, \, s \geq n; \\
        \f_2^{\oplus ( [s - n/2 \equiv 0, 5, 6 \text{ mod } 8] + 2[s - n/2 \equiv 4 \text{ mod } 8] )} & \text{if } n \equiv 2 \text{ mod } 4, \, n \geq 6, \, s \geq n/2; \\
        \f_2^{\oplus (2[s - n \equiv 0 \text{ mod } 4] + [s - n \equiv 6 \text{ mod } 8] + [s - n \equiv 7 \text{ mod } 8])} & \text{if } n \equiv 3 \text{ mod } 4, \, n \ge 3, \, s \ge n; \\
        \f_2 & \text{if } n \equiv 2 \text{ mod } 4, s < n/2; \\
        \f_2 & \text{if } n \equiv 3 \text{ mod } 4, s < n; \\
        0 & \text{otherwise},

    \end{cases}
\]
where $n$ denotes the filtration, $s$ denotes the stem, and $[P]$ denotes the indicator of a proposition $P$. Namely,
\[
    [P] =
    \begin{cases}
        1 & \text{ if } P \text{ is true}; \\
        0 & \text{ if } P \text{ is false}.
    \end{cases}
\]
Here we adopt the convention that $\f_2^{\oplus 0}=0$.

To solve extension problems, we use the motivic Adams spectral sequence, following the strategy described in \cref{kq_{**}BC_2}.

\begin{figure}[H]
\includegraphics[width=14cm]{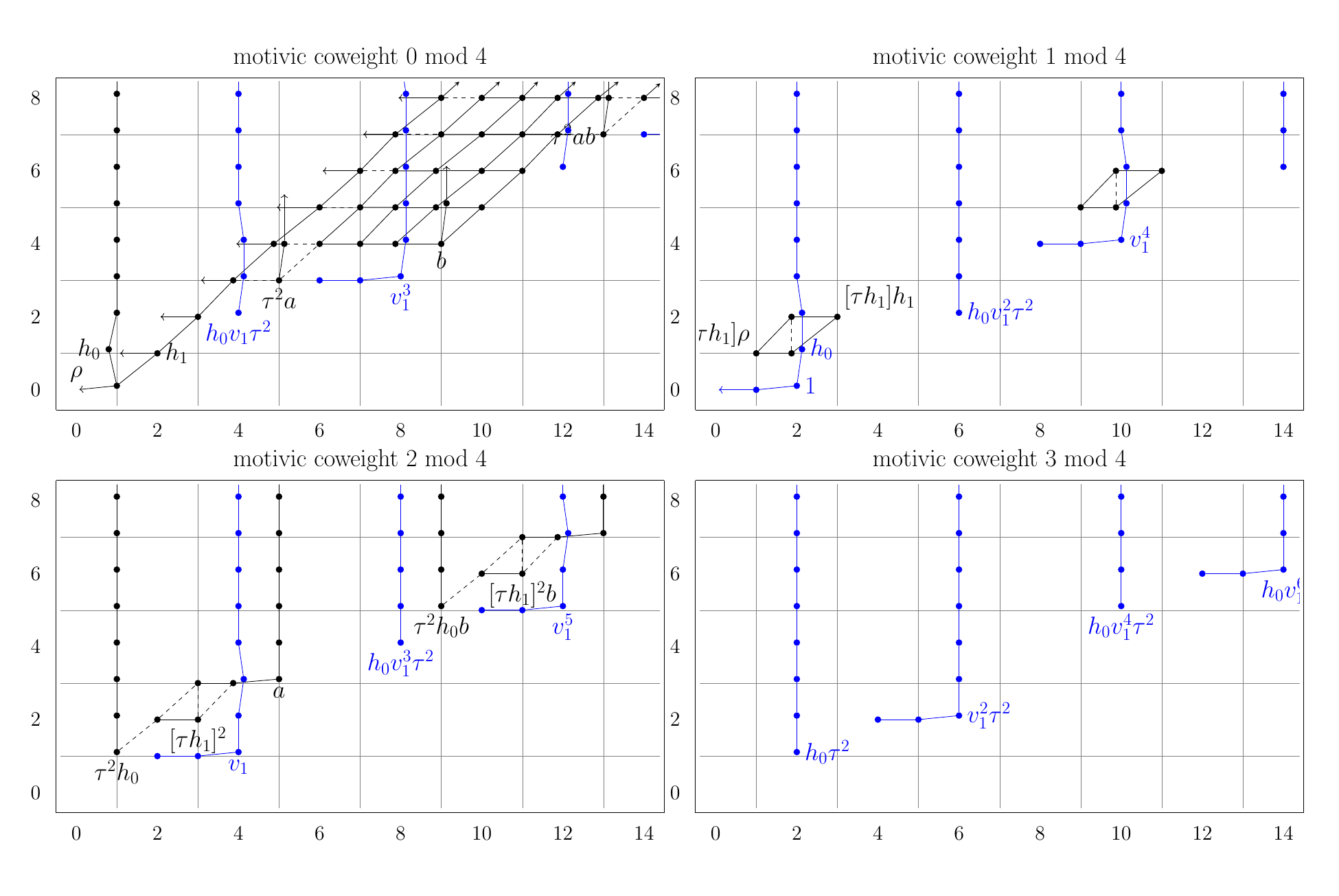}\vspace{-0.5cm} 
    \caption{The charts represents the long exact sequence for $\Ext^{***}_{\ca(1)}(Q)$, organized by coweights congruent to modulo 4. The black portion shows $\Ext^{***}_{\ca(1)}(\Sigma^{1,1}\m_2)$, and the blue portion shows $\Ext^{***}_{\ca(1)}(C)$. The connecting homomorphism is given by multiplication by $h_0$, and computing kernels and cokernels gives $\Ext_{\ca(1)}^{***}(Q)$. 
    Each $\bullet$ represents $\f_2[\tau^4]$.}
    \label{fig:LES for Ext(Q)}
    \includegraphics[scale=0.40]{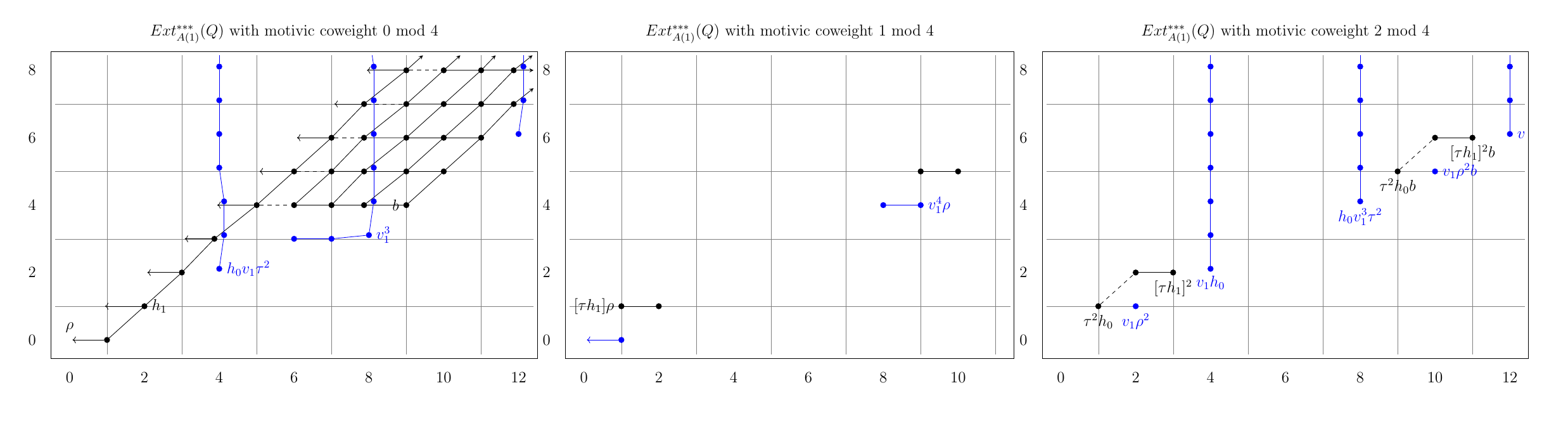}\vspace{-0.5cm}
    \caption{The charts represents $\Ext^{***}_{\ca(1)}(Q)$, organized by coweights congruent to modulo 4. 
    Each $\bullet$ represents $\f_2[\tau^4]$.} 
    \label{fig:Ext(Q) } 
\end{figure}

We begin by computing $\Ext_{\ca(1)}^{***}(Q)$ derived from the short exact sequence \eqref{eqn: les for kq general strategy}. Note that $\Ext_{\ca(1)}^{***}(\Sigma^{2,1}C) \cong  \Ext_{\ce(1)}^{***}(\Sigma^{2,1} \m_2^{**})$ was already computed in \cref{Fig. 1}, and $\Ext_{\ca(1)}^{***}(\Sigma^{1,1}\m_2^{**})$ was determined in \cref{fig:kqR}. The result for $\Ext_{\ca(1)}^{***}(Q)$ then follows directly and is shown in \cref{fig:LES for Ext(Q)}.
\vspace{-0.5cm}

Similarly, for $\Ext_{\ca(1)}^{***}(R)$, using the standard filtration of $R$ \eqref{eqn: standard filtration of R in kq general strategy} gives
$$\Ext_{\ca(1)}^{***}(R) \cong \bigoplus_{i \geq 0} \Sigma^{4i+3,2i+2} \Ext_{\ca(0)}^{***}(\m_2^{**}) \cong \bigoplus_{i \geq 0} \Sigma^{4i+3,2i+2} \frac{\f_2[\tau^2,h_0,\rho]}{(\rho h_0)}.$$

\begin{figure}[H]
    \includegraphics[width =15cm]{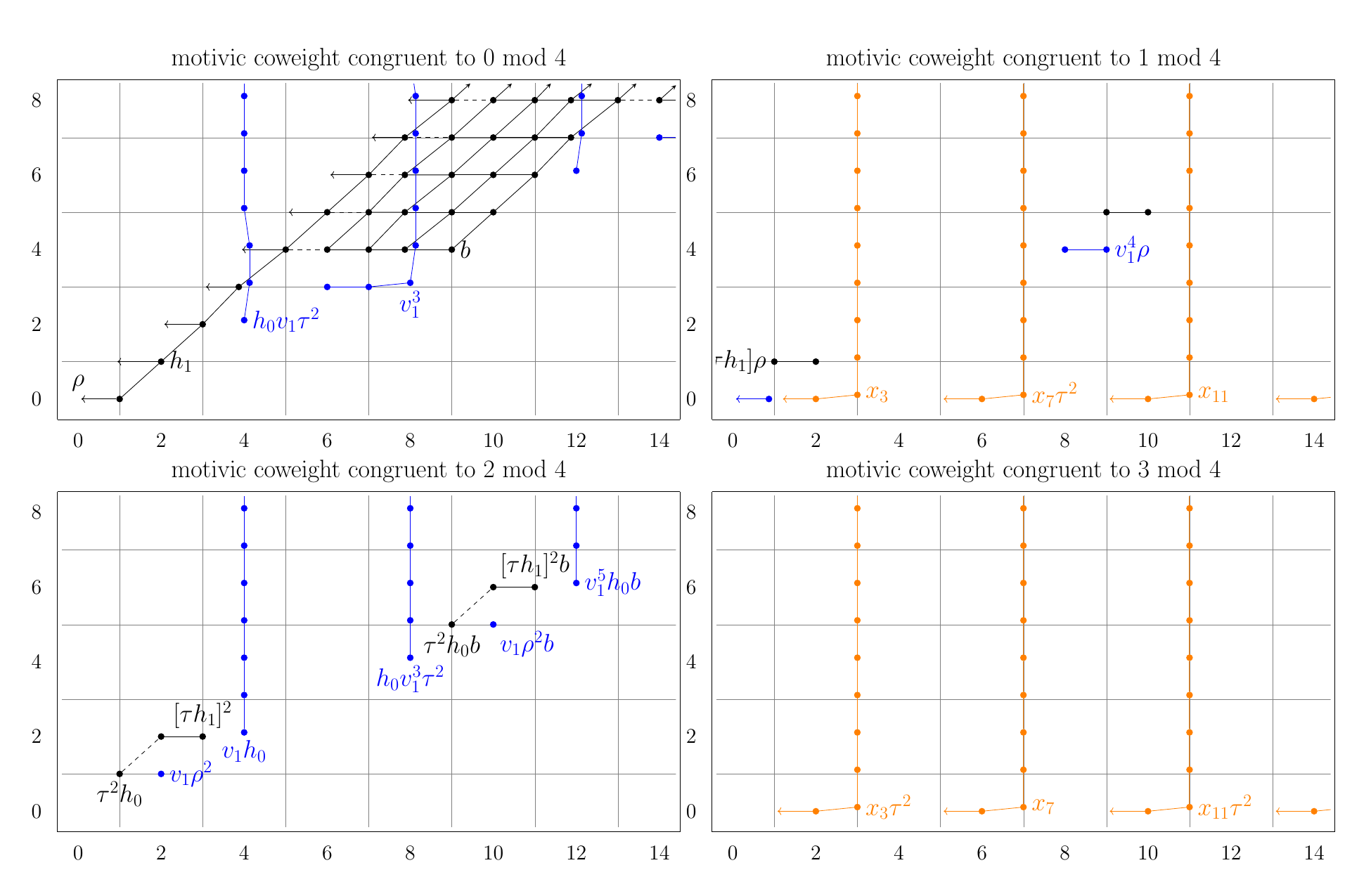}\vspace{-0.5cm}
    \caption{The charts represents the long exact sequence for $\Ext^{***}_{\ca(1)}(H^{**}BC_2)$ organized by coweights congruent to modulo 4 where the orange portion shows $\Ext^{***}_{\ca(1)}(R)$. 
    Each $\bullet$ represents $\f_2[\tau^4]$.}  
    \includegraphics[width=15cm]{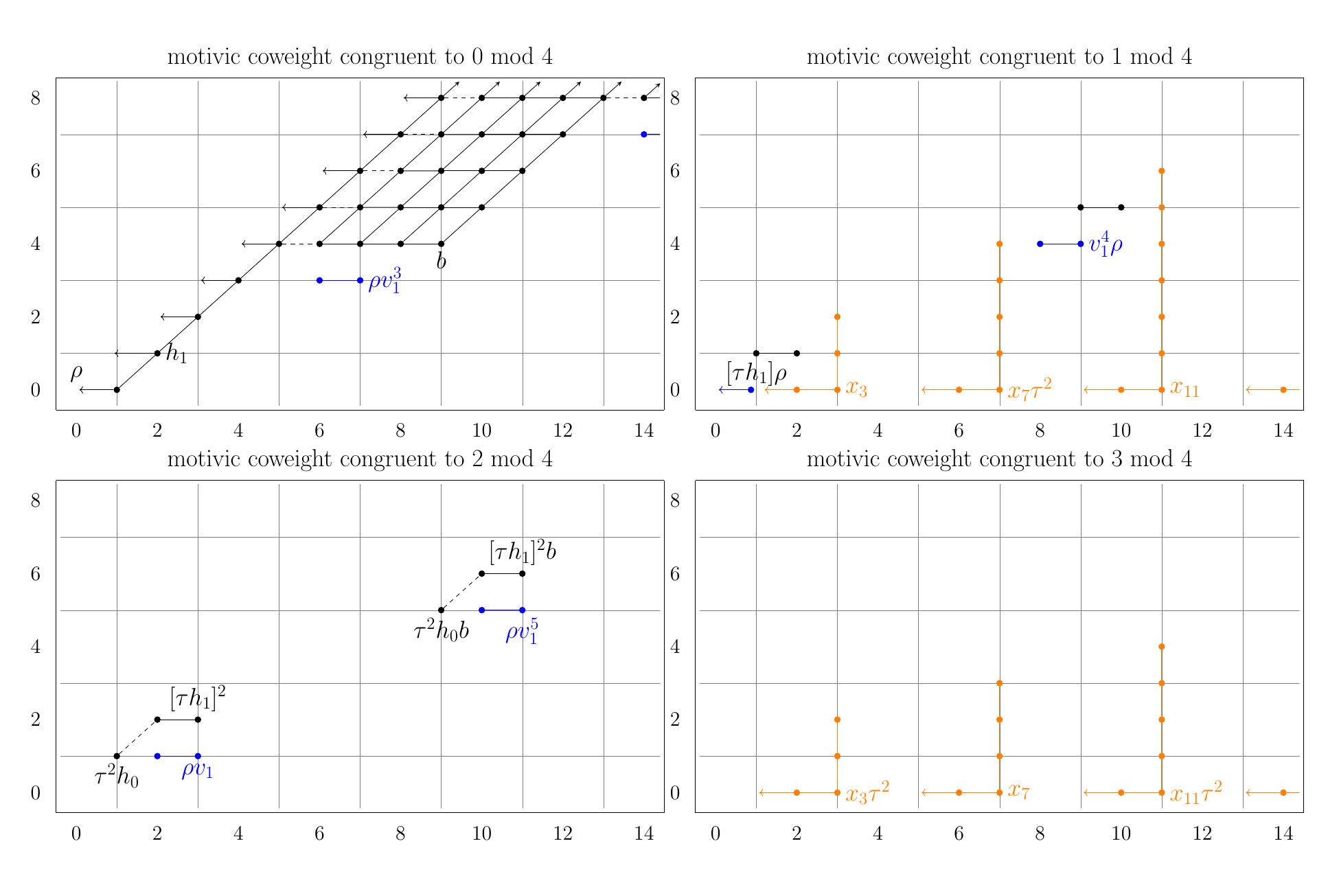}\vspace{-0.5cm}
    \caption{The charts represents $\Ext^{***}_{\ca(1)}(H^{**}BC_2)$ organized by coweights congruent to modulo 4. 
    Each $\bullet$ represents $\f_2[\tau^4]$.} 
    \label{fig:kq_ASS_E2_R}
\end{figure}

We now substitute this result into the long exact sequence in $\Ext$ derived from \eqref{eqn: les for kq with Q and R general strategy} to compute the $E_2$-term, where the Adams spectral sequence collapses as shown in \cref{fig:kq_ASS_E2_R}. 
The connecting homomorphism must be nontrivial for the sizes of $kq_{**}BC_2$ to match those from the Atiyah--Hirzebruch spectral sequence.

Combining the charts in \cref{fig:kq_ASS_E2_R}, we obtain the final result as follows.

\begin{thm}
    The $E_\infty$-page of the motivic Adams spectral sequence converging to the very effective Hermitian K-theory of $BC_2$ over real closed fields is given, as a module over the $E_\infty$-page of the motivic Adams spectral sequence for $kq_{**}$, by 
    \begin{flalign*}
        & E_\infty \cong \f_2[\tau^4, \rho, h_1, b] \oplus \frac{\rho v_1^3 \cdot \mathbb{F}_2[\tau^4, \rho, v_1^3]}{\rho^3} \\
        & \quad \oplus \frac{\tau h_1 \cdot \f_2[\tau^4, \rho, b]}{\rho^2} \oplus \frac{\rho \cdot \f_2[\tau^4, \rho, v_1^4]}{\rho^3 v_1^4} \oplus \left ( \bigoplus_{\substack{j \geq 0 \\ 2|j}} \Sigma^{4j-1, 4j-1} \frac{\tau^2 \cdot \f_2[\tau^4, \rho, h_0]}{(\rho h_0, h_0^{j+3})} \right ) \oplus \left ( \bigoplus_{\substack{j \geq 0 \\ 2 \nmid j}} \Sigma^{4j-1, 4j-1} \frac{ \f_2[\tau^4, \rho, h_0]}{(\rho h_0, h_0^{j+3})} \right ) \\
        & \quad \oplus \frac{\tau^2 (h_1^2, h_0) \cdot \f_2[\tau^4, \rho, b, h_1^2, h_0]}{(\rho^2, h_0^2, h_1^4, \rho h_0, h_1^2 h_0)} \oplus \frac{\rho v_1 \cdot \f_2[\tau^4, \rho, v_1^4]}{\rho^3 v_1^4} \\
        & \quad \oplus \left ( \bigoplus_{\substack{j \geq 0 \\ 2\nmid j}} \Sigma^{4j-1, 4j-1} \frac{\tau^2 \cdot \f_2[\tau^4, \rho, h_0]}{(\rho h_0, h_0^{j+3})} \right ) \oplus \left ( \bigoplus_{\substack{j \geq 0 \\ 2 | j}} \Sigma^{4j-1, 4j-1} \frac{ \f_2[\tau^4, \rho, h_0]}{(\rho h_0, h_0^{j+3})} \right ).
    \end{flalign*}
    Note that the result is organized into four rows, which from top to bottom correspond to the cases in \cref{fig:kq_ASS_E2_R} where the coweight is congruent to 0,1,2, and 3 (mod 4), respectively.
\end{thm}
\section{Computations over finite fields}\label{Sec:Fq}

Throughout this section, we work over a finite field $\f_q$ with $q$ an odd prime. 

\subsection{Mod two motivic cohomology}
We have
\[
H^{**} = \m_2^{**} \cong
\begin{cases}
    \f_2[\tau,\mu]/(\mu^2) \quad & \text{if } q \equiv 1 \pmod 4, \\
    \f_2[\tau,\rho]/(\rho^2) \quad & \text{if } q \equiv 3 \pmod 4,
\end{cases}
\]
where $|\tau| = (0,1)$ and $|\mu| = |\rho| = (1,1)$; see, for example, \cite{Kyl15}.

The elements $\rho$ and $\mu$ lie in
\[
H^{1,1}(\Spec \f_q; \z/2) \cong k^M_1(\f_q)/2 \cong \f_q^\times/(\f_q^\times)^2.
\]
When $q \equiv 1 \mod 4$, the class $\mu$ corresponds to a chosen non-square element of $\f_q^\times$ under the identification. For any $q$ (in fact, in the mod two motivic cohomology of any field), $\rho$ corresponds to the class of $-1$ under the identification. In the case $q \equiv 3 \pmod 4$, $\rho$ is nontrivial, whereas it vanishes when $q \equiv 1 \pmod 4$. 

Moreover, one has
\[
\Sq^1(\tau) =
\begin{cases}
0 \quad & \text{if } q \equiv 1 \pmod 4, \\
\rho \quad & \text{if } q \equiv 3 \pmod 4.
\end{cases}
\]


\subsection{Integral motivic homology}

By \cref{Prop:HZBC2}, we have
$$H\z_{**}(BC_2) \cong 
\begin{cases}
\bigoplus_{i \geq 0} \Sigma^{2i+1,i+1} \f_2[\tau,\mu]/\mu^2 \quad & \text{ if } q \equiv 1 \pmod 4, \\
\bigoplus_{i \geq 0} \Sigma^{2i+1,i+1} \f_2[\tau,\rho]/\rho^2 \quad & \text{ if } q \equiv 3 \pmod 4.
\end{cases}
$$

\subsection{Effective algebraic K-theory}

We start by computing $kgl_{**}$ using the motivic Adams spectral sequence 
$$E_2 = \Ext_{\ca}^{***}(H^{**}kgl) \cong \Ext_{\ce(1)}^{***}(\m_2^{**}) \Rightarrow kgl_{**}.$$

If  $q \equiv 1 \pmod 4$, then, since $\mu$ is central, we have
$$E_2 = \f_2[\tau, \mu]/(\mu^2)[h_0,v_1]$$
with $|h_0| = (1,1,0)$ and $|v_1| = (1,3,1)$. 

If $q \equiv 3 \pmod 4$, then there are nontrivial $\rho$-Bockstein differentials generated by $d_1(\tau) = \rho h_0$, yielding
$$E_2 = \f_2[\tau^2,\rho][h_0,v_1]/(\rho^2, \rho h_0)\{1, \rho \tau\}.$$

Unlike in all of the previous contexts, the motivic Adams spectral sequence does not collapse. Instead, there are differentials generated by
\begin{equation}\label{eqn:d_nFq}
d_{\nu(q-1)+s}(\tau^{2^s}) = \mu \tau^{2^s-1} h_0^{\nu(q-1)+s}, \quad s \geq 0, \quad q \equiv 1 \pmod 4,
\end{equation}
\begin{equation}
    d_{\nu(q^2-1)+s-1}(\tau^{2^s}) = \rho \tau^{2^s-1} h_0^{\nu(q^2-1)+s-1}, \quad s \geq 1, \quad q \equiv 3 \pmod 4,
\end{equation}
where $\nu(-)$ denotes $2$-adic valuation. See \cite[Section~4.2]{Kyl15} for details.

\begin{figure}[H]
\includegraphics[width=\linewidth]{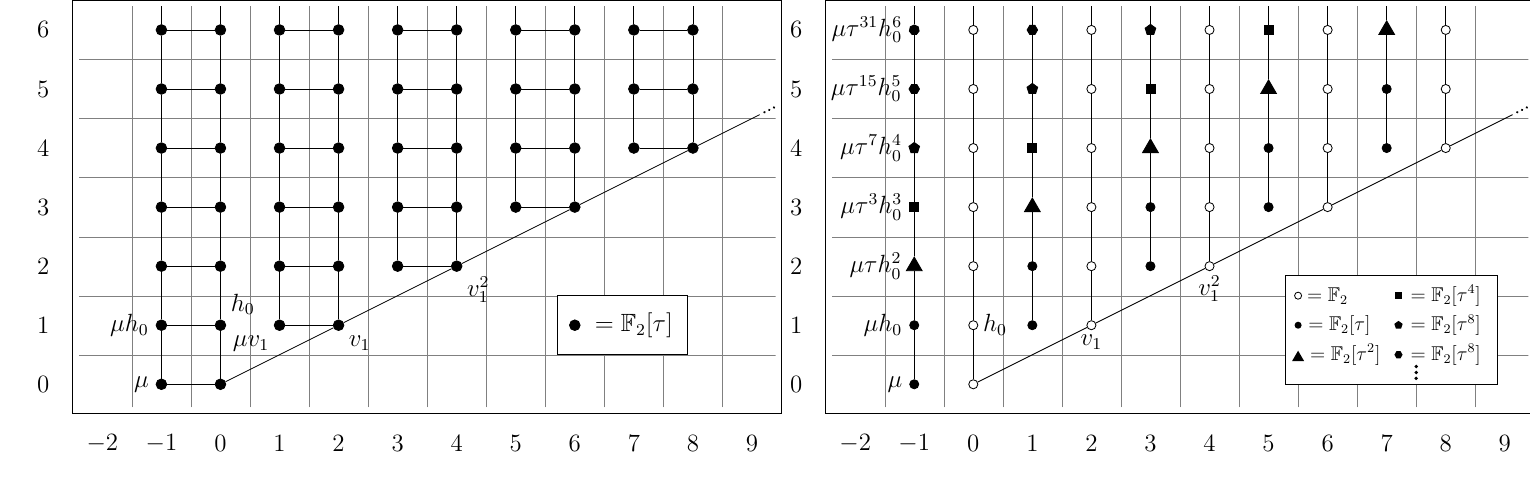}
\includegraphics[width=\linewidth]{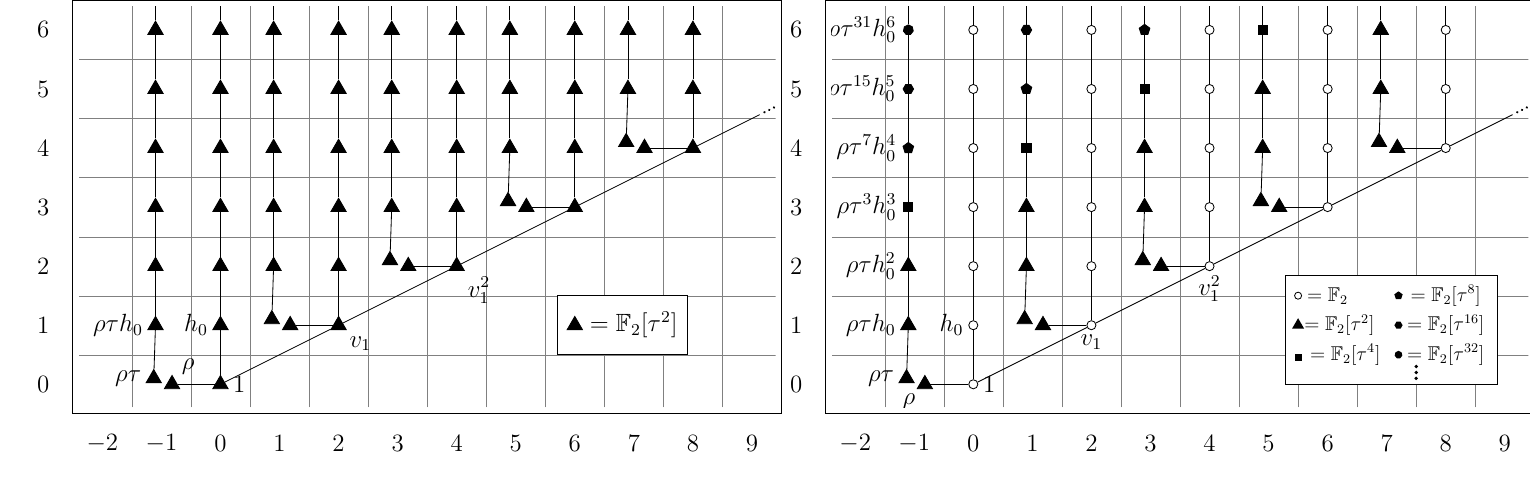}
\caption{The $E_2$ and $E_\infty$ terms of the motivic Adams spectral sequence for $kgl_{**}$ over $\f_5$ (top two charts) and $\f_3$ (bottom two charts).}
\label{kglF5}
\end{figure}


\begin{remark}
    Before proceeding, we note that all charts are drawn over $\f_5$ and $\f_3$. The same computations extend to $\f_q$ for larger $q$, with the cases $q \equiv 1 \pmod{4}$ and $q \equiv 3 \pmod{4}$ exhibiting analogous behavior. The only change is in the infinite tower in stem $-1$ of $kgl_{**}$, which takes the form
\(
\bigoplus_{i \geq 0} \z/2^{f_q(i)}\{\tau^i\},
\)
where $f_q(i)$ depends on $q$ and on the $2$-adic valuation of $i$, and satisfies $f_q(i) \geq 1$ for all $i \geq 0$. Aside from this, the overall strategy for computing $kgl_{**}BC_2$ remains unchanged.

For simplicity and clarity, we restrict our attention to computing over $\mathbb{F}_5$ and $\mathbb{F}_3$, where the functions $f_q(i)$ are given by
\begin{equation}
f_5(i)=
\begin{cases}
    2 & i \neq 2^k-1, \\
    k+2 & i = 2^k-1,
\end{cases}
\qquad
f_3(i)=
\begin{cases}
    1 & i \neq 2^k-1, \\
    k+2 & i = 2^k-1.
\end{cases}\label{f35}
\end{equation}
\end{remark}

\begin{figure}[H]
\includegraphics[width=\linewidth]{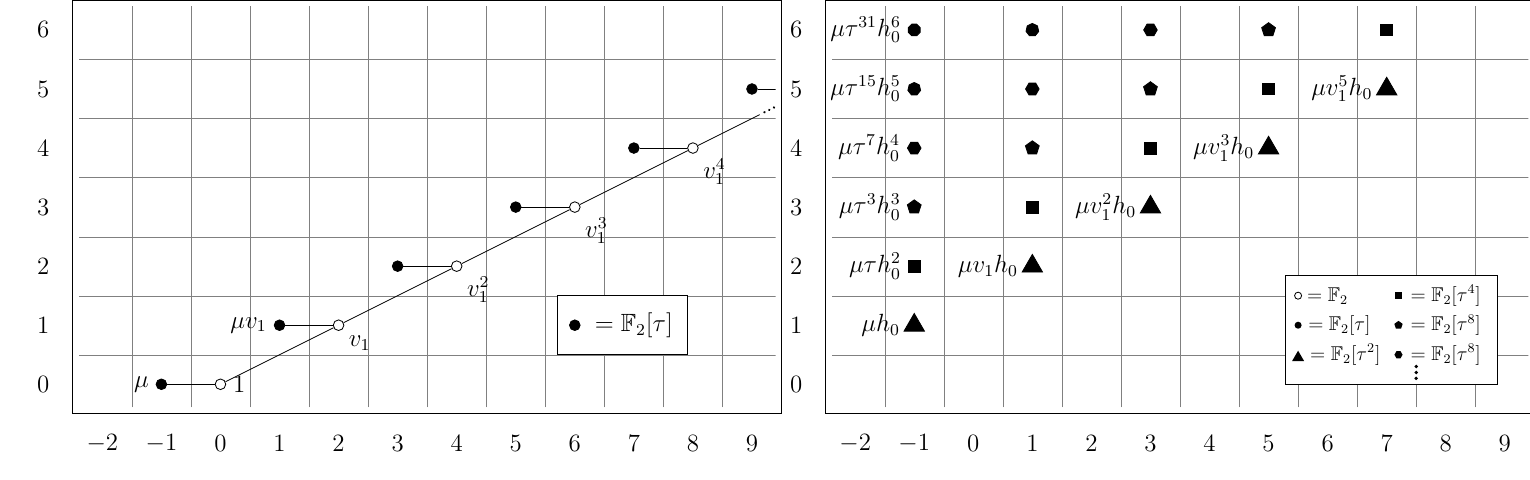}
\caption{$E_2^{2k+1,**}$ (left) and $E_2^{2k,**}$ (right) for the Atiyah--Hirzebruch spectral sequences converging to $gr_{AH}^*kgl_{**}BC_2$ over $\f_5$. A copy of this pattern appears in every positive filtration row.}\label{fig:kglBC2F5E2}
\end{figure}

We first consider the case when $q \equiv 1 \pmod{4}$. The $d_1$-differentials are given by multiplication by $2$ from the even to the odd filtrations. We depict a portion of the $E_2$-term over $\f_5$ in \cref{fig:kglBC2F5E2}. 
 
\begin{figure}[H]
\includegraphics[width=\linewidth]{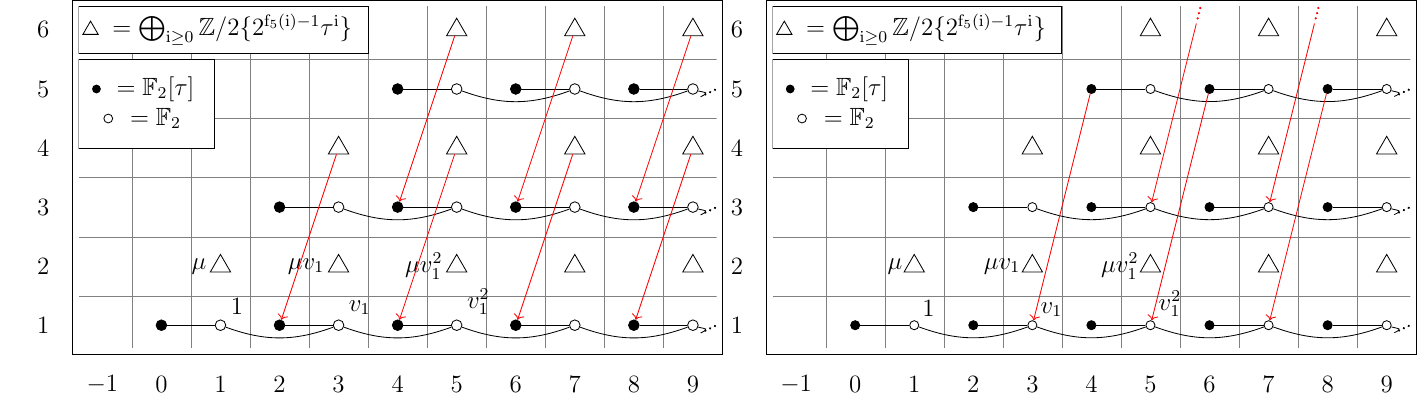}
\caption{The $E_2=E_3$ term for the Atiyah--Hirzebruch spectral sequence for $gr_{AH}^*kgl_{**}BC_2$ over $\f_5$ with potential $d_3$-differentials (left) and potential $d_4$-differentials (right).}
\label{kglBC2F5AH}
\end{figure}

Since $d_2$ is induced by multiplication by $\eta$, which is zero in $kgl$, we conclude that $E_2 = E_3$. See \cref{kglBC2F5AH}.

Recall that $d_3$ corresponds to the multiplication by $v_1$ and is $h_0$-linear. All possible sources of differentials are $h_0$-multiples, which makes the $d_3$-differentials trivial, since $h_0v_1=0$ on the $E_3$-page.
Moreover, since all differentials are linear in $\mu$, we have $d_{2k+1}=0$ and $d_{2k}=0$ for all $k$. Consequently, $E_2 = E_3 = E_4$ is already the $E_\infty$-page. 

Notice that the combinations of hollow triangles and hollow circles represent hidden $\tau$-extensions, as may be verified by comparison with the motivic Adams spectral sequence below.

Computations over $\f_q$ with $q \equiv 3 \pmod{4}$ are similar, though slightly more complicated. 

\begin{figure}[H]
\includegraphics[width=\linewidth]{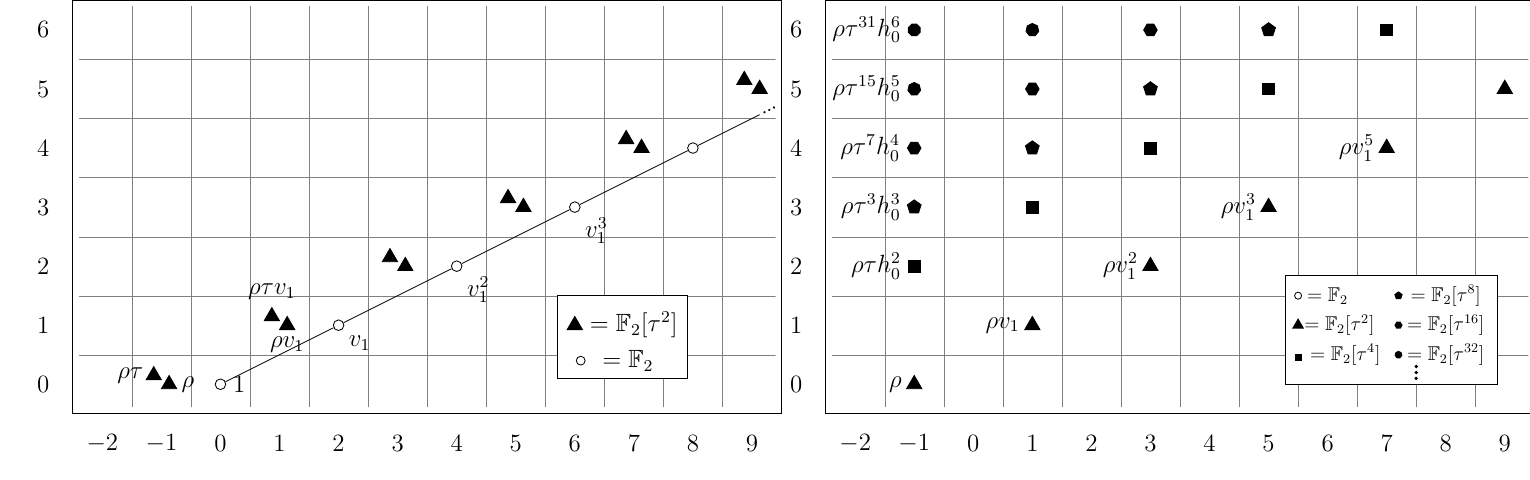}
\caption{$E_2^{2k+1,**}$ (left) and $E_2^{2k,**}$ 
(right) for the Atiyah--Hirzebruch spectral sequence converging to $kgl_{**}BC_2$ over $\f_3$. A copy of this pattern appears in every positive filtration row.}
\label{fig:kglBC2F3E2}
\end{figure}

We begin with the $E_1$-page with bottom-right chart of \cref{kglF5} in each filtration. As before, the $d_1$-differentials are given by multiplication by $2$ from even to odd filtrations. We depict the $E_2$-term over $\f_3$ in \cref{fig:kglBC2F3E2} by deleting the classes that support $h_0$ multiplication in even filtrations and the classes obtained by $h_0$-multiplication in odd filtrations. 

\begin{figure}[H]
\includegraphics[width=\linewidth]{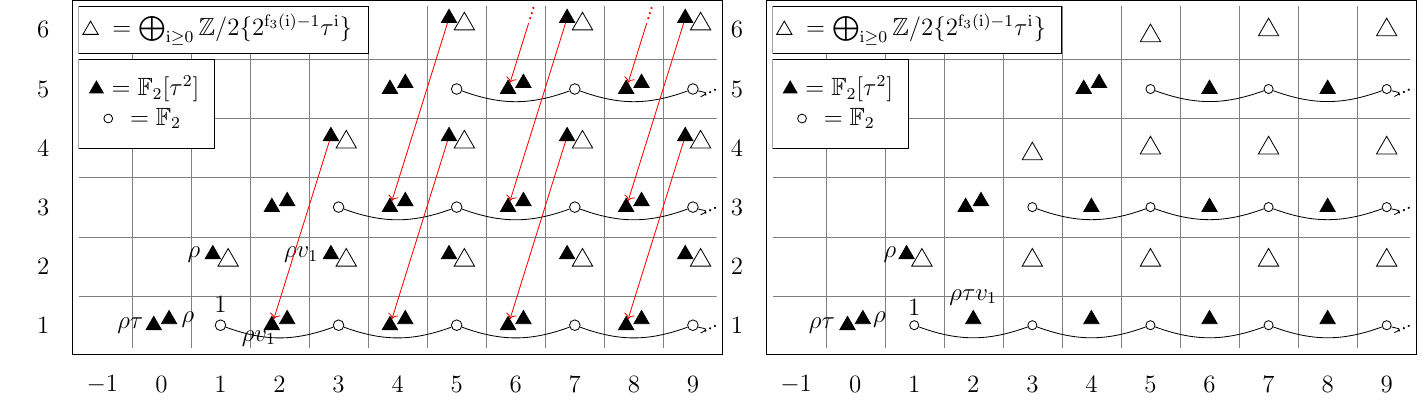}
\caption{On the left is the $E_3$ page with $d_3$-differentials in red for the Atiyah--Hirzebruch spectral sequence converging to $gr_{AH}^* kglBC_2$ at $\f_3$ and, the $E_4$=$E_{\infty}$ page is as given on the right.}
\label{kglbc2f5ahe4}
\end{figure}

There are no $d_2$-differentials as $\eta$ vanishes in $kgl$, and the $d_3$-differential is given by $d_3(z)=v_1z$, as shown in \cref{kglbc2f5ahe4}. Higher differential vanish for the same reason as in case of $\f_5$.
The spectral sequence collapses at the $E_4$-page, yielding $gr_{AH}^* kglBC_2$ over $\f_3$.

We now turn to the motivic Adams spectral sequence to solve hidden extensions, which takes the form
$$E_2 = \Ext_{\ce(1)}^{***}(H^{**}BC_2) \Rightarrow kgl_{**}BC_2,$$
following the same general strategy outlined in \cref{SS:kgl_general}. We have already computed (see the left charts in \cref{kglF5})
$$\Ext_{\ce(1)}^{***}(\Sigma^{2,1} \m_2^{**}) \cong 
\begin{cases}
\Sigma^{2,1} \f_2[\tau,h_0,v_1,\mu]/\mu^2 & \text{ if } q \equiv 1 \pmod 4, \\
\Sigma^{2,1} \f_2[\tau^2,\rho,h_0,v_1]/(\rho^2,\rho h_0)\{1,\tau \rho\}  & \text{ if } q \equiv 3 \pmod 4. 
\end{cases}$$

\begin{figure}[H]
\includegraphics[width=\linewidth]{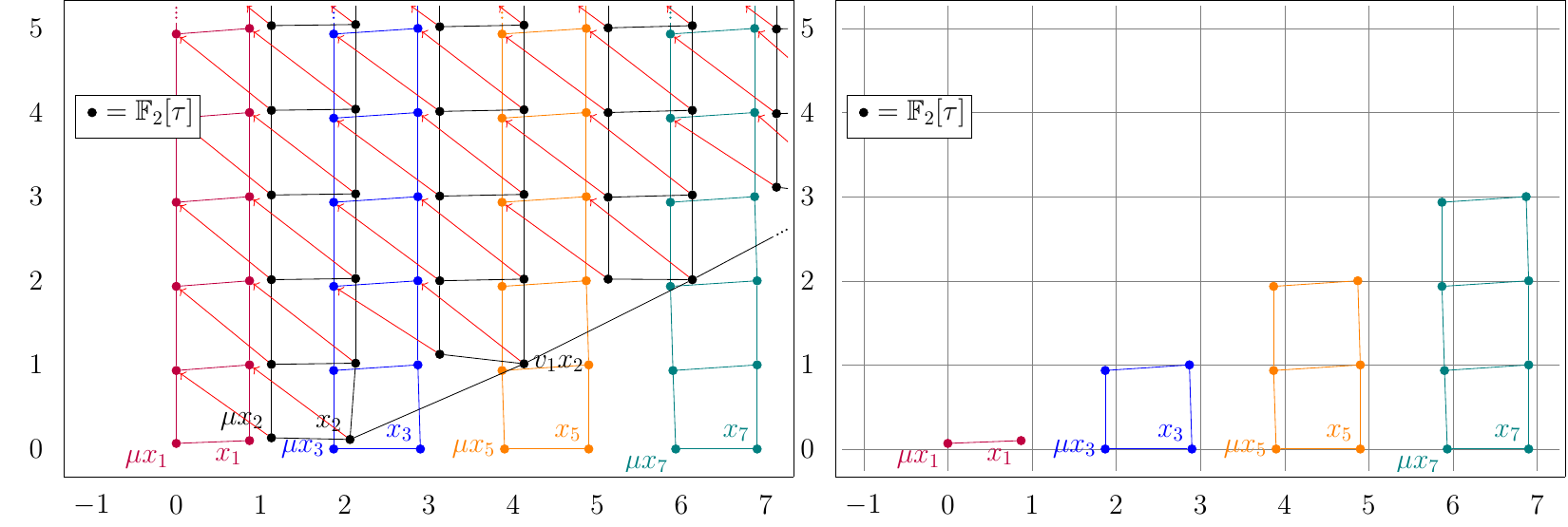}

\includegraphics[width=\linewidth]{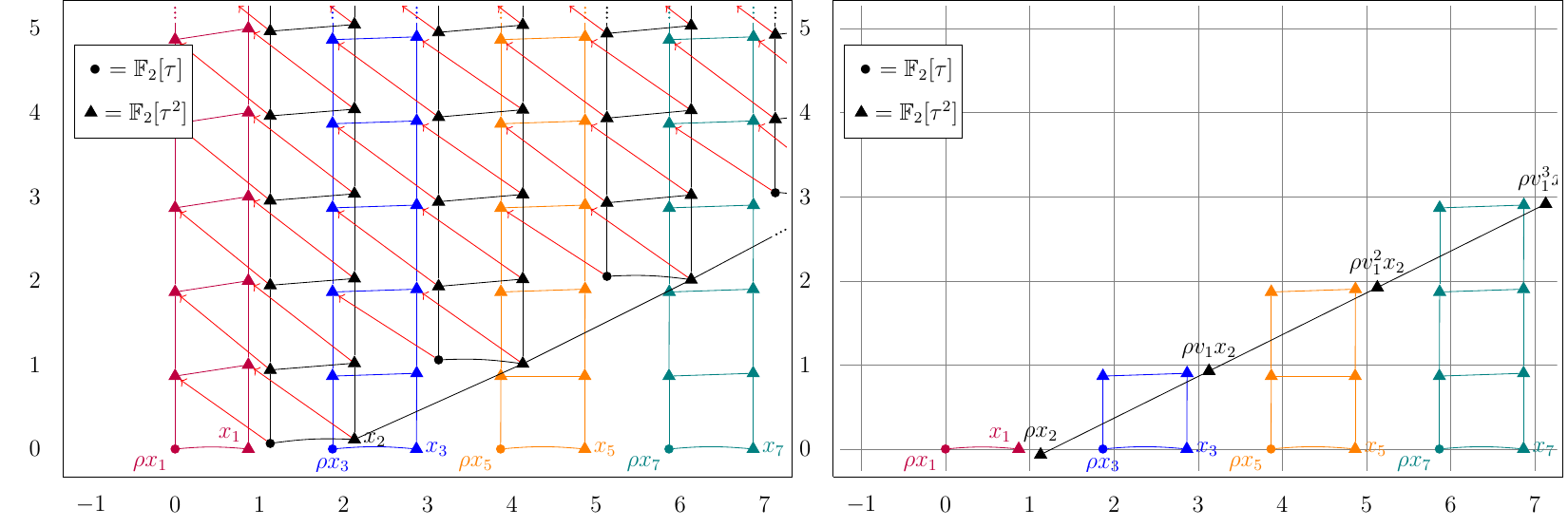}
\caption{The long exact sequence computing $\Ext_{\ce(1)}(H^{**}(BC_2))$ and the resulting $\Ext$ groups over $\f_5$ (top charts) and $\f_3$ (bottom charts). The element $x_i$ denotes the generator corresponding to the shifted copy $\Sigma^{i, \lfloor \tfrac{i + 1}{2} \rfloor}$.}\label{fig:kglASSFq}
\end{figure}

On the other hand, similar computations with the $\rho$-Bockstein spectral sequence (see \cref{SS:6.4}) show that the $E_1$-term of the spectral sequence
$$E_1 = \bigoplus_{i \geq 0} \Sigma^{2i+1,i+1} \Ext_{\ca(0)}^{***}(\m_2^{**}) \Rightarrow \Ext^{***}_{\ce(1)}(V)$$
can be identified with
$$E_1 \cong
\begin{cases}
\bigoplus_{i \geq 0} \Sigma^{2i+1,i+1} \f_2[\mu,\tau,h_0]/(\mu^2) \quad & \text{ if } q \equiv 1 \pmod 4, \\
\bigoplus_{i \geq 0} \Sigma^{2i+1,i+1} \f_2[\rho, \tau^2, h_0]/(\rho^2, \rho h_0)\{1,\tau \rho\} \quad & \text{ if } q \equiv 3 \pmod 4.
\end{cases}$$
Since the Bockstein differentials are $\mu$- and $\rho$-linear, the spectral sequence collapses, yielding the result of $\Ext^{***}_{\ce(1)}(V)$ shown in \cref{fig:kglASSFq}; the surviving classes are displayed in multiple colors.

From the long exact sequence in $\Ext$ derived from \eqref{eqn: les for kgl general strategy}, we see that, in order for the rank of the associated graded to $kgl_{**}BC_2$ as an $\m_2^{**}$-module to be correct, the connecting homomorphism must satisfy
$$\delta(\Sigma^{2,1} xv_1^i) = \Sigma^{2i+1,i+1}xh_0^{i+1}.$$ 
The resulting long exact sequences and $E_2$-terms over $\f_5$ and $\f_3$ are depicted in \cref{fig:kglASSFq}. 

Base change to the algebraic closure shows that there are no nontrivial Adams differentials, thereby yielding the final answer, depicted in the right-hand charts of \cref{fig:kglASSFq}. 

\begin{thm}
The $E_\infty$-page of the motivic Adams spectral sequence converging to the effective algebraic \(K\)-theory of \(BC_2\), as a module over the \(E_\infty\)-page for \(kgl_{**}\), is given by
\[
E_\infty \cong
\begin{cases}
\displaystyle
\bigoplus_{i\ge 0}
\Sigma^{2i+1,i+1}
\frac{\f_2[\tau,\mu,h_0]}
     {(h_0^{i+1},\mu^2)},
& q \equiv 1 \pmod 4,
\\[1.5em]
\displaystyle
\Sigma^{2,1}\f_2[\tau^2,v_1]\{\rho\} \oplus \\
\bigoplus_{i\ge 0}
\left(
\Sigma^{2i+1,i+1}
\frac{\f_2[\tau^2,\rho,h_0]}
     {(h_0^{i+1},\rho^2)}
\;\oplus\;
\Sigma^{2i,i}\f_2[\tau^2]\{\tau\rho\}\right),
& q \equiv 3 \pmod 4.
\end{cases}
\]
\end{thm}

\begin{remark}
Observe that there are hidden $\tau$-extensions in the Atiyah--Hirzebruch spectral sequence: the combinations of triangles and hollow circles in \cref{kglBC2F5AH} and \cref{kglbc2f5ahe4} form copies of $\f_2[\tau]$ over $\f_5$ and $\f_2[\tau^2]$ over $\f_3$, respectively.
\end{remark}

\subsection{Very effective Hermitian K-theory}\label{SS:6.4}

We turn to the computation of $kq_{**}(BC_2)$ over finite fields. To do so, we use the same strategy by first computing $kq_{**}$. The motivic Adams spectral sequence converging to $kq_{**}$  has signature
$$E_2=\Ext_{\ca}^{***}(H^{**}kq)\cong\Ext^{***}_{\ca(1)}(\m_2^{**})\Rightarrow kq_*.$$
This spectral sequence was analyzed by Kylling in \cite{Kyl15}; we summarize the relevant details here. 

We will compute this $E_2$-term using the $\mu$-Bockstein spectral sequence (for $q \equiv 1 \text{ mod } 4$) or the $\rho$-Bockstein spectral sequence (for $q \equiv 3 \text{ mod } 4$).  The $\mu$-Bockstein spectral sequence has signature
$$E_1=\Ext_{\ca(1)}((\m_2^\c)^{**})[\mu]/\mu^2\Rightarrow\Ext_{\ca(1)}(\m_2^{**})$$

When $q \equiv 1 \pmod 4$, as computed in \cref{kqC}, 
$$\Ext_{\ca(1)}^{***}((\m_2^\c)^{**}) \cong \dfrac{\f_2[\tau,h_0,h_1,a,b]}{h_0h_1,\tau h_1^3, h_1 a, a^2 = h_0^2 b},$$
the $E_1$-term of the $\mu$-Bockstein spectral sequence follows. Since $\mu$ is central in $\m_2^{**}\f_q$, the spectral sequence collapses at the $E_1$-page, giving

\[\Ext_{\ca(1)}^{***}(\m_2^{**}) \cong \dfrac{\f_2[\tau,h_0,h_1,a,b, \mu]}{h_0h_1,\tau h_1^3, h_1 a, a^2 = h_0^2 b, \mu^2}
\]

\begin{figure}[H]
\includegraphics[width=\linewidth]{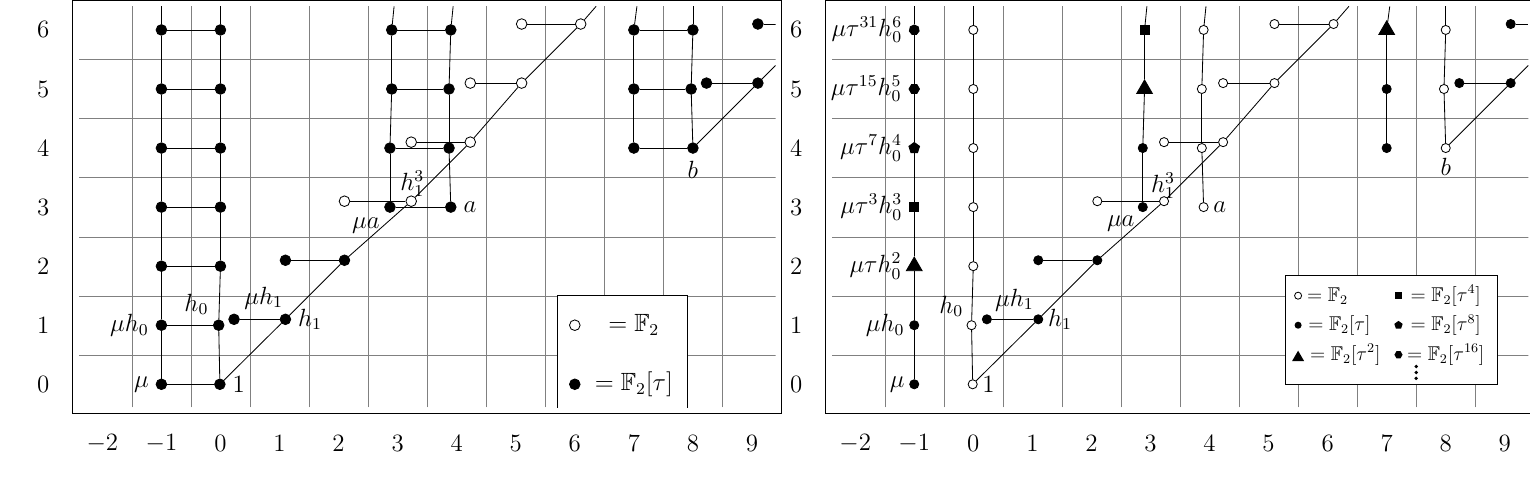}
\includegraphics[width=\linewidth]{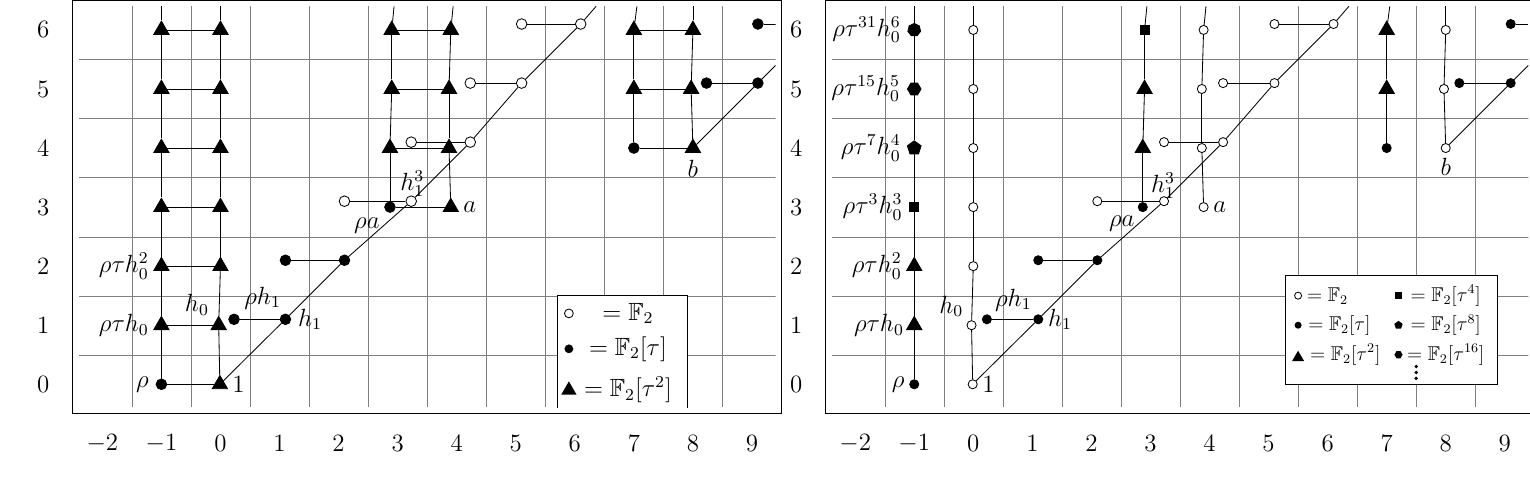}
\caption{The $E_2$- and $E_\infty$-pages of the motivic Adams spectral sequence for $kq_{**}$ over $\f_5$ (top two charts) and $\f_3$ (bottom two charts).}
\label{kqF5}
\end{figure}

When $q \equiv 3 \pmod 4$, the $\rho$-Bockstein spectral sequence
$$E_1 = \Ext_{\ca(1)}^{***}((\m_2^\c)^{**})[\rho]/\rho^2 \Rightarrow \Ext_{\ca(1)}^{***}(\m_2^{**}\f_q)$$
has nontrivial differentials generated by the Leibniz rule and
$$d_1(\tau) = \rho h_0.$$
The $\rho$-Bockstein spectral sequence collapses at the $E_2$-page, giving the $E_2$-term of the motivic Adams spectral sequence depicted in the bottom left of \cref{kqF5}. 

As in the Adams spectral sequences for $H\z_{**}$ and $kgl_{**}$ over $\f_q$, the differentials in the Adams spectral sequence for $kq_{**}$ are generated by 
$$d_{\nu(q-1)+s}(\tau^{2^s}) = \mu \tau^{2^s-1} h_0^{\nu(q-1)+s}, \quad s \geq 0, \quad q \equiv 1 \pmod 4,
$$
$$d_{\nu(q^2-1)+s-1}(\tau^{2^s}) = \rho \tau^{2^s-1} h_0^{\nu(q^2-1)+s-1}, \quad s \geq 1, \quad q \equiv 3 \pmod 4.$$
These Adams differentials imply the Toda brackets in $kq_{**}$
\begin{equation}\label{eqn:TBF1}
\langle \mu \tau^{2^s-1} 2^{\nu(q-1)+s-1}, 2, \eta \rangle = \tau^{2^s} \eta, \quad s \geq 0, \quad q \equiv 1 \pmod 4,
\end{equation}

$$\langle \rho, 2, \eta \rangle = \tau \eta, \quad q \equiv 3 \pmod 4,$$
$$\langle \rho \tau^{2^s-1} 2^{\nu(q^2-1)+s-2}, 2, \eta \rangle = \tau^{2^s} \eta, \quad s \geq 1, \quad q \equiv 3 \pmod 4.$$

We also note the presence of hidden $\eta$-extensions 
$$\eta \cdot \tau^{i+1} \eta^2 \beta^j = 2 \cdot 2^{\nu(q-1)+\nu(i+1)-1}\tau^i \alpha \beta^j$$
for all $i \geq 0$ when $q \equiv 1 \pmod 4$, and for all $i \equiv 1 \pmod 2$ when $q \equiv 3 \pmod 4$. If $q \equiv 3 \pmod 4$ and $i \equiv 0 \pmod 2$, the relation becomes
$$\eta \cdot \tau^{i+1} \eta^2 \beta^j = 2 \tau^i \alpha \beta^j.$$ These relations are visible in the slice spectral sequence; see \cite[Fig. A5]{KQ24}). 

We are now ready to analyze the Atiyah--Hirzebruch spectral sequences and Adams spectral sequences. 

\subsubsection{Computations over $\f_q$ with $q \equiv 1 \pmod 4$}

\begin{figure}[H]
\includegraphics[width=9cm]{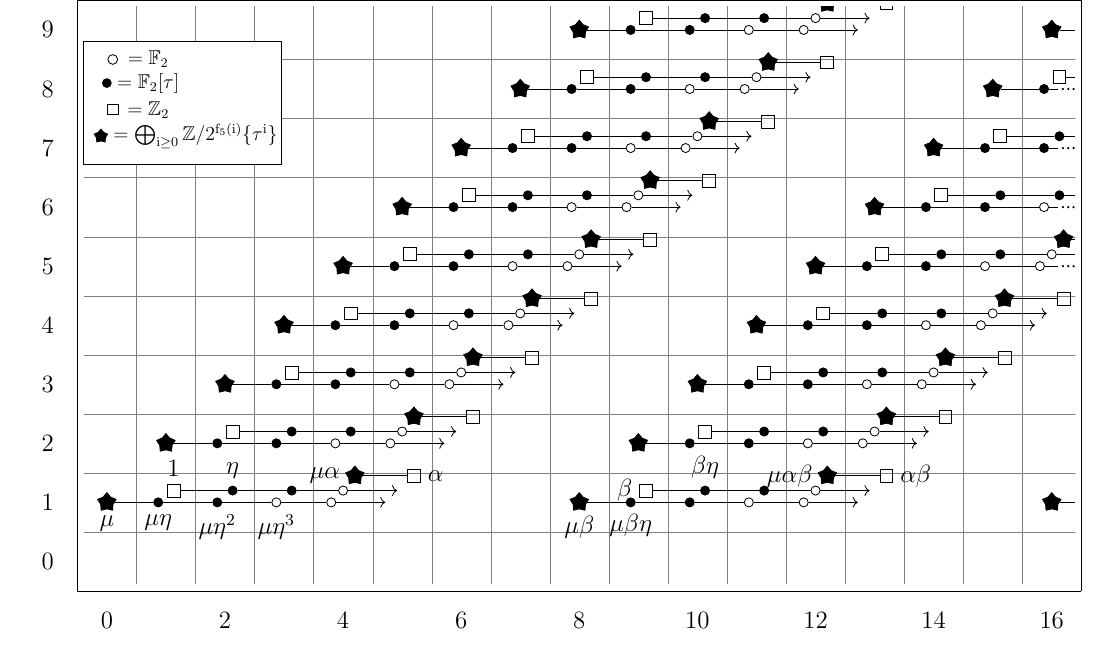}
\includegraphics[width=9cm]{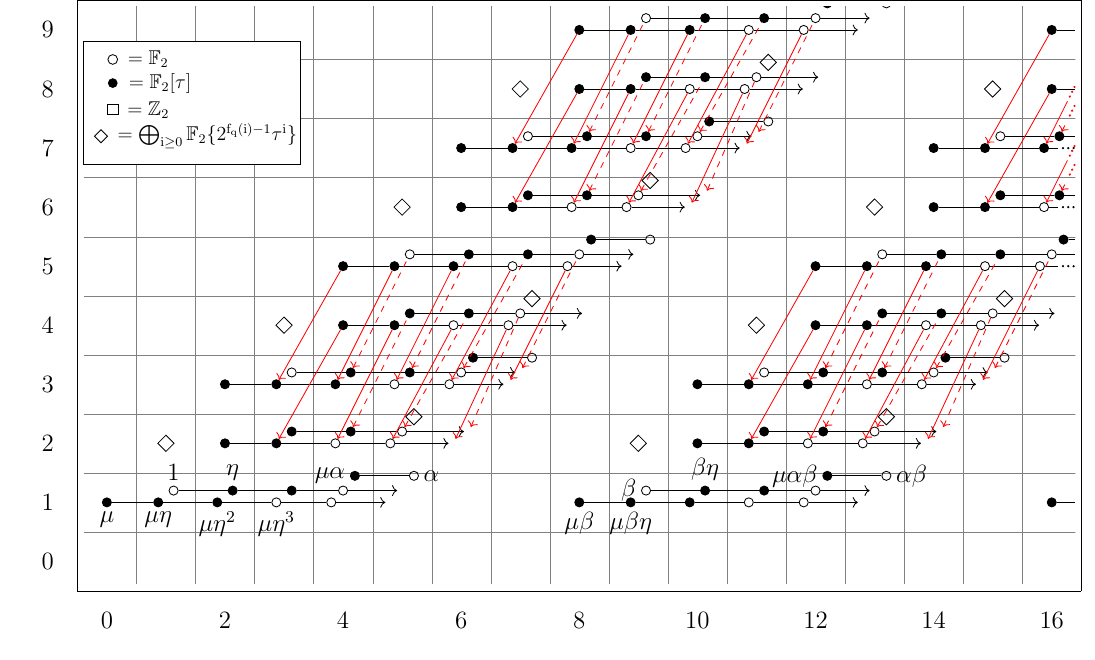}\vspace{-0.4cm}
\caption{The $E_1$-page (above) and $E_2$-page (below) of the Atiyah--Hirzebruch spectral sequence converging to $kq_{**}(BC_2)$ over $\f_5$. Red lines represents the differentials given by multiplication by $\eta$.}\label{Fig:kqAHSS_Fq_E2}
\end{figure}

We first compute the $gr_{AH}^*kq_{**}(BC_2)$ using the Atiyah--Hirzebruch spectral sequence.  The $E_1$-term, depicted in \cref{Fig:kqAHSS_Fq_E2}, contains shifted copies of $kq_{**}$, depicted in \cref{kqF5}.

The only nontrivial $d_1$-differentials are multiplication by $2$ between $\square$'s and $\bigstar$'s:
$$\cdot 2 : \z_2 \to \z_2, \quad \cdot 2: \bigoplus_{i \geq 0} \z/2^{f_q(i)} \{\tau^i\} \to \bigoplus_{i \geq 0} \z/2^{f_q(i)} \{\tau^i\},$$
where $f_q(i)$ is some positive function depending on $q$ and the 2-adic valuation of $i$, with $f_q(i) > 1$ for all $i \geq 0$. 

In the first case, the kernel is zero and the cokernel is $\z/2$. In the second case, the kernel and cokernel are
$$\bigoplus_{i \geq 0} \z/2\{2^{f_q(i)-1} \tau^i\}, \quad \bigoplus_{i \geq 0} \z/2 \{ \tau^i\},$$ respectively, since $f_q(i) > 1$ for all $i \geq 0$. In terms of the charts, we have 
\begin{equation}
f_5(i)=
\begin{cases}
    2 & i \neq 2^k-1, \\
    k+2 & i = 2^k-1.
\end{cases}
\end{equation}
The $E_2$-page is obtained by deleting the classes that support $h_0$-multiplication, as well as the classes that are $h_0$-multiples. The resulting $E_2$-term is depicted in \cref{Fig:kqAHSS_Fq_E2}. 

\begin{figure}[H]
\includegraphics[width=9cm]{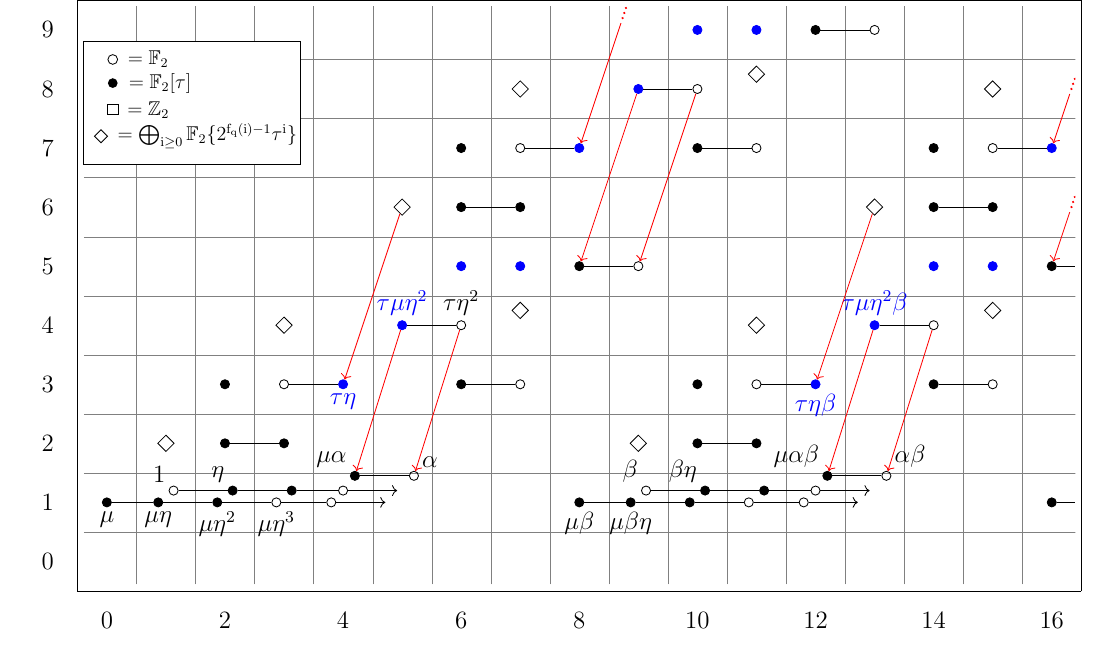}\vspace{-0.5cm}
\caption{The $E_3$-page of the Atiyah--Hirzebruch spectral sequence converging to $kq_{**}(BC_2)$ over $\f_5$. Differentials are drawn in red.}\label{Fig:kqAHSS_Fq_E3}
\end{figure}

Again, the $d_2$-differentials follow from considering the action of $\eta$, as shown in \cref{Fig:kqAHSS_Fq_E2}. The resulting $E_3$-term is depicted in \cref{Fig:kqAHSS_Fq_E3}. The nontrivial $d_3$-differentials follow from the Toda bracket
$$\alpha \in \langle \tau \eta^2, \eta, 2 \rangle$$
and the Toda brackets \eqref{eqn:TBF1}. 

\begin{figure}[H]
\includegraphics[width=12cm]{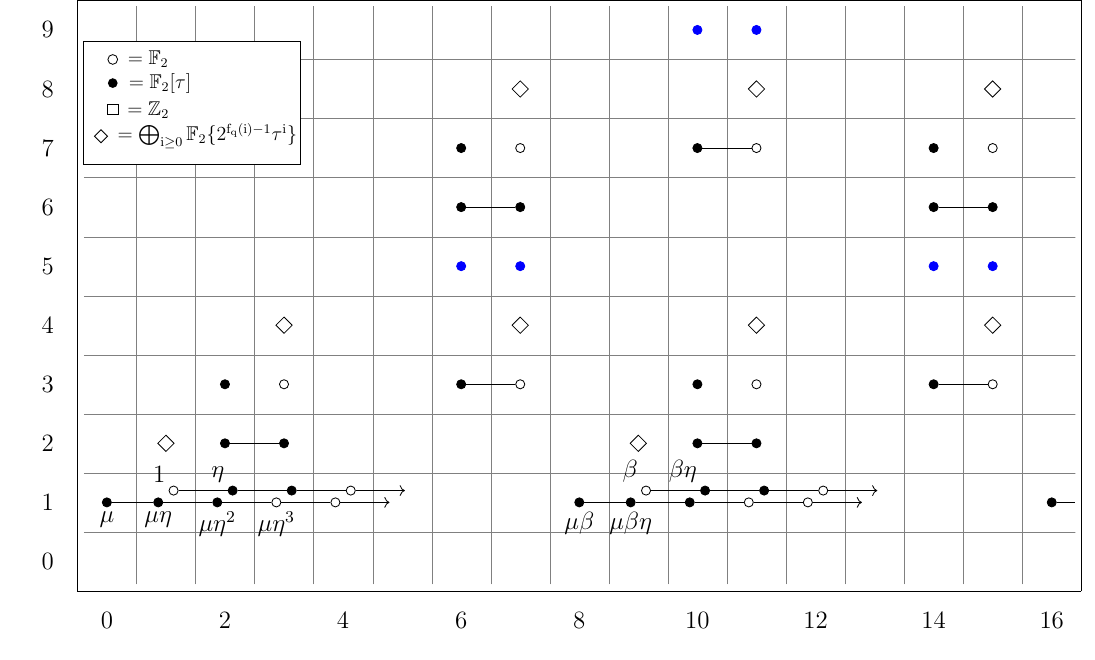}
\caption{The $E_4=E_\infty$-page of the Atiyah--Hirzebruch spectral sequence converging to $kq_{**}(BC_2)$ over $\f_5$.}\label{Fig:kqAHSS_Fq_E4}
\end{figure}

The $E_4$-page is depicted in \cref{Fig:kqAHSS_Fq_E4}. From our analysis of the motivic Adams spectral sequence, it follows that the Atiyah--Hirzebruch spectral sequence collapses at $E_4$.

\begin{figure}[H]
\includegraphics[width=9cm]{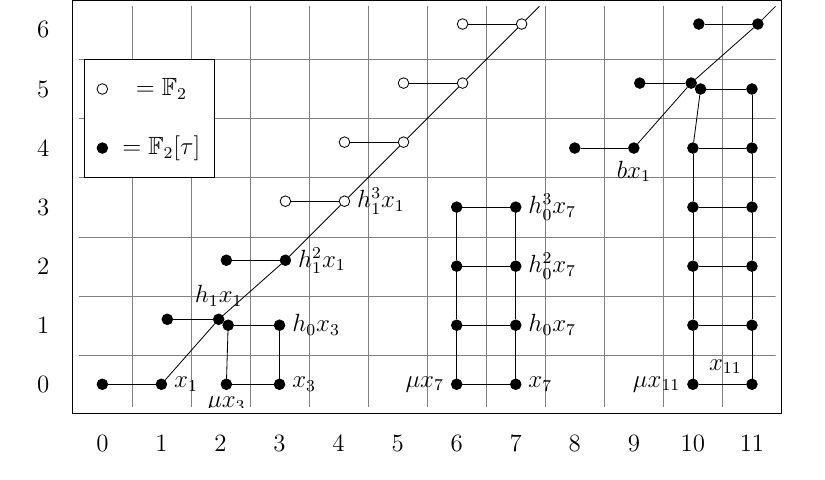}
\caption{The $E_2=E_\infty$-page of the motivic Adams spectral sequence converging to $kq_{**}(BC_2)$ over $\f_q$ when $q \equiv 1 \pmod 4$. The element $x_i$ denotes the generator corresponding to the shifted copy $\Sigma^{i, \lfloor \tfrac{i + 1}{2} \rfloor}$. }\label{Fig:kqASS_F1}
\end{figure}

Again, the combinations of hollow diamonds and hollow circles represent hidden $\tau$-extensions, as verified by comparison with the motivic Adams spectral sequence, see \cref{Fig:kqASS_F1}.

We now analyze the motivic Adams spectral sequence
$$E_2 = \Ext_{\ca(1)}^{***}(H^{**}BC_2) \Rightarrow kq_{**}BC_2.$$
One approach is to follow the strategy described in \cref{kq_{**}BC_2}, beginning with a computation of the $E_2$-page of the Adams spectral sequence arising from the short exact sequence \eqref{eqn: les for kq general strategy}. However, since $\mu$ is central in $\m_2^{**}$, the $\Ext$-term over $\f_q$ is simply the $\Ext$-term over $\c$ tensored with $\f_2[\mu]/\mu^2$, where $|\mu| = (0,-1,-1)$ \cite[7.1]{WO17}. See \cref{Fig:kqASS_F1}.

We claim that the spectral sequence collapses at the $E_2$-page. To see this, we first observe that the only potential Adams differentials originate on the  truncated $h_0$-towers in stems congruent to $2$ or $3$ modulo $4$. Since the towers in stems congruent to $2$ mod $4$ are $\mu$-multiples of the classes in stems congruent to $3$ mod $4$, it suffices to show that the classes in stems congruent to $3$ mod $4$ are permanent cycles. 

This follows from base change to the algebraic closure: if any of the classes in stem $4k+3$ over $\f_q$ supported a nontrivial differential, then the degree of $2$-torsion in the  group $kq^{\f_q}_{4k+3,2k+2}(BC_2)$ would be less than the degree of $2$-torsion $kq^{\bar{\f}_q}_{4k+3,2k+2}(BC_2)$. On the other hand, inspection of the Adams spectral sequences implies that base change sends generators to generators, leading to a contradiction.

\begin{thm}
The very effective Hermitian \(K\)-theory of \(BC_2\) over \(\f_q\) with \(q \equiv 1 \pmod 4\) , as a module over the \(E_\infty\)-page of \(kq_{**}\), is given by
\[
kq_{**}(BC_2)\cong \Sigma^{1,1}
\frac{\f_2[\tau,\mu,h_1,b]}
     {(\tau h_1^3,\mu^2)} \oplus 
\bigoplus_{i\ge 0}
\Sigma^{4i+3,2i+2}
\frac{\f_2[\tau,\mu,h_0]}
     {(h_0^{2i+2},\mu^2)} .
\]
\end{thm}

\subsubsection{Computations over $\f_q$ with $q \equiv 3 \pmod 4$}

The analysis of the Atiyah--Hirzebruch spectral sequence when $q \equiv 3 \pmod 4$ is similar to the case $q \equiv 1 \pmod 4$ described above, but there are enough subtle differences to merit a detailed description of the calculation. 

The $E_1$-term contains a copy of $kq_{**}$ in each filtration; we refer to \cite[Fig.~A5]{KQ24} for these coefficients. Alternatively, the $E_1$-term can be represented by the same chart (see \cref{Fig:kqAHSS_Fq_E2}) as in the case $q \equiv 1 \pmod 4$, except that now each $\bigstar$ represents $\bigoplus_{i \geq 0} \z/2^{g_q(i)}\{\tau^i\}$ for some positive function $g_q(i)$ depending on $q$ and the 2-adic valuation of $i$. In particular, over $\f_3$, one has
\begin{equation}
g_3(i)=
\begin{cases}
    1 & i \neq 2^k-1, \\
    k+2 & i = 2^k-1.
\end{cases}\label{f35}
\end{equation}

Observe that, unlike the functions $f_q(i)$ appearing for $q \equiv 3 \pmod 4$, we have $g_q(i)=1$ whenever $i \equiv 0 \pmod 2$. These groups arise from the differential $d_1(\tau) = \rho h_0$ in the $\rho$-Bockstein spectral sequence for $\Ext_{\ca(1)}(\m_2^{**})$ over $\f_q$ with $q \equiv 3 \pmod 4$. 

\begin{figure}[H]
\includegraphics[width=11cm]{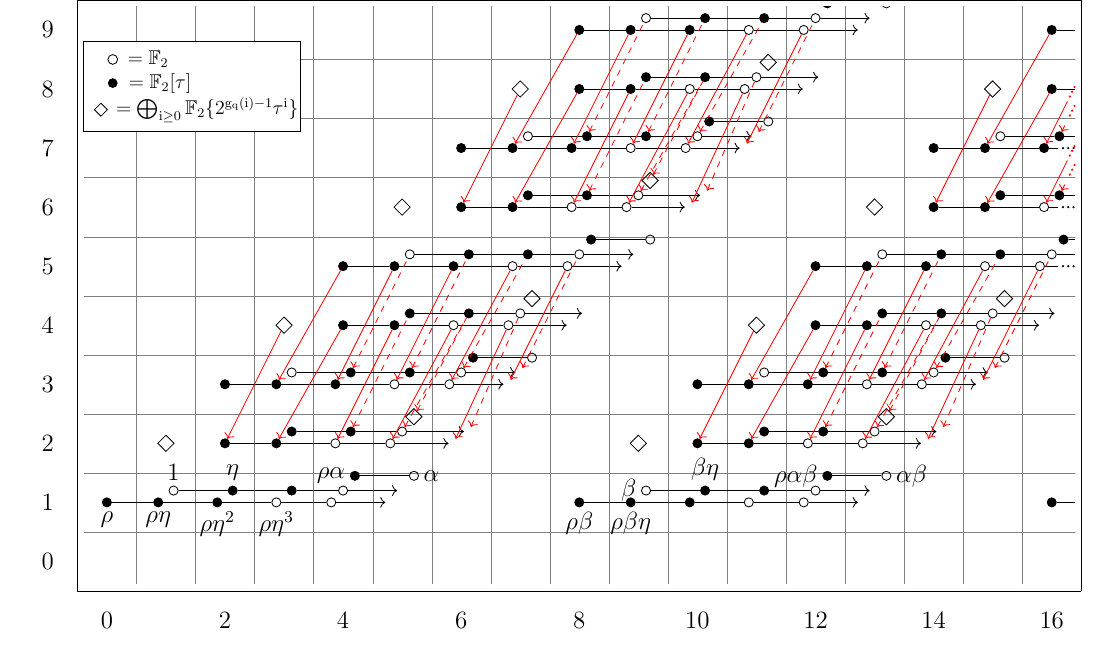}
\caption{The $E_2$-page of the Atiyah--Hirzebruch spectral sequence converging to $kq_{**}(BC_2)$ over $\f_3$. Differentials are shown in red.}\label{Fig:kqAHSS_F3_E2}
\end{figure}

After calculating $d_1$-differentials, which again are given by multiplication by 2, we obtain the $E_2$-page depicted in  \cref{Fig:kqAHSS_F3_E2}. 

\begin{figure}[H]
\includegraphics[width=7.5cm]{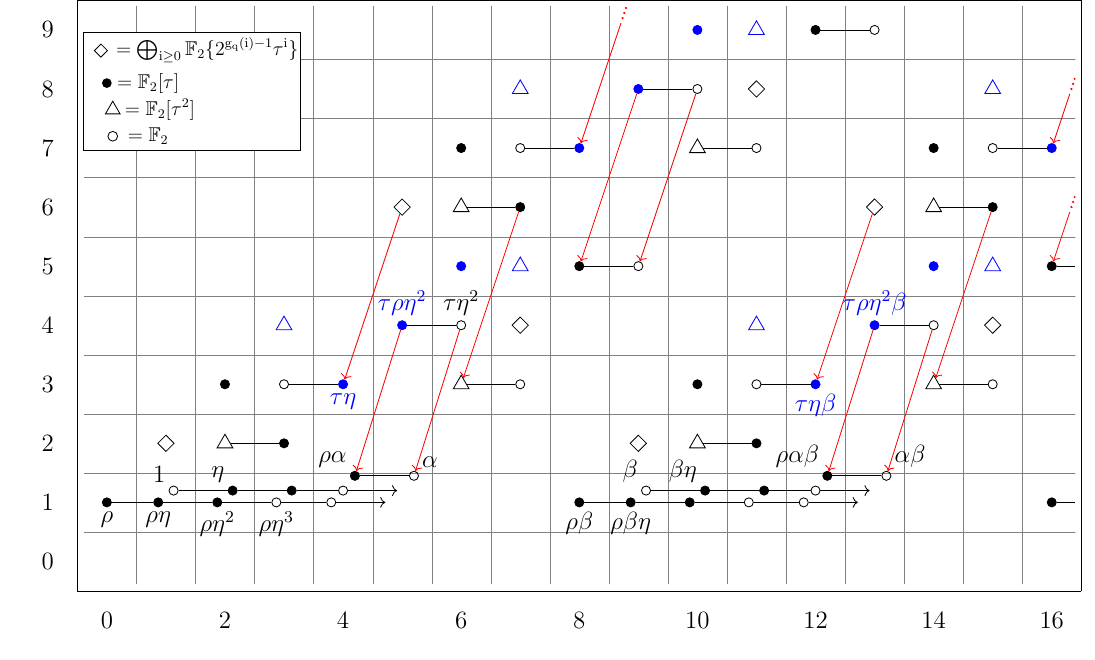}
\includegraphics[width=7.5cm]{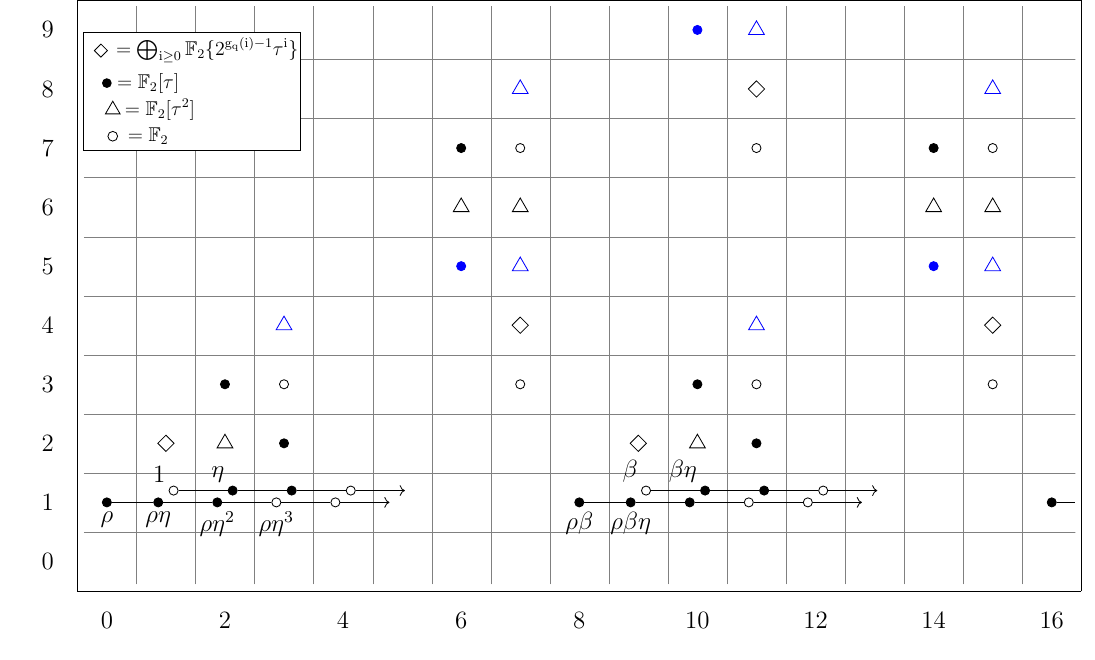}
\caption{The $E_3$-page (left) and $E_4=E_\infty$-page (right) of the Atiyah--Hirzebruch spectral sequence converging to $kq_{**}(BC_2)$ over $\f_3$. The red arrows represent nontrivial $d_3$-differentials and blue classes represents that their generators are $\tau$-multiples. }\label{Fig:kqAHSS_F3_E3}
\end{figure}

\begin{figure}[H]
\includegraphics[width=10cm]{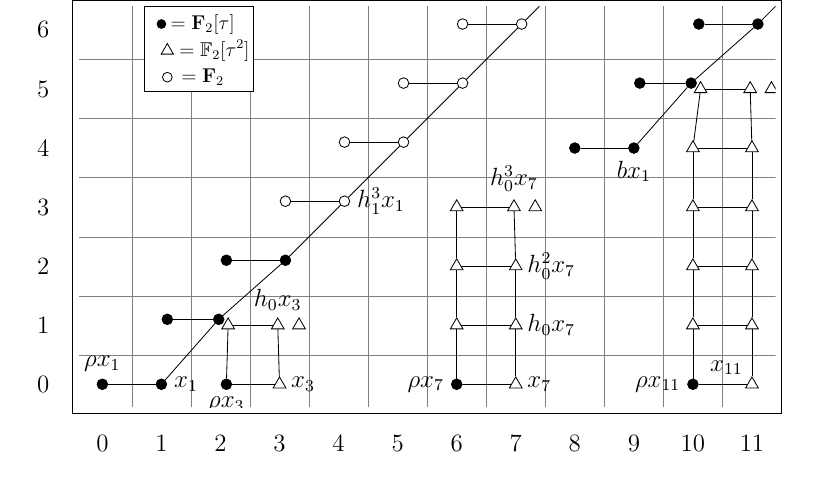}
\caption{The $E_2=E_\infty$-page of the motivic Adams spectral sequence converging to $kq_{**}(BC_2)$ over $\f_3$. }\label{Fig:kqASS_F3}
\end{figure}

From here, we compute kernels and cokernels of multiplication by $\eta$, yielding the $E_3$-page depicted in \cref{Fig:kqAHSS_F3_E3}. The $E_3$-page depends only on the filtration modulo $4$, with the exception of filtration $1$, which remains unchanged from the $E_2$-page.

We assemble the different filtrations into a single chart in \cref{Fig:kqAHSS_F3_E3}. The $d_3$-differentials follow from the Toda brackets listed at the beginning of the section, together with the cell structure of $BC_2$. The resulting $E_4$-term is shown in the same figure. From our analysis of the motivic Adams spectral sequence, it follows that the Atiyah--Hirzebruch spectral sequence collapses at the $E_4$-page.

We now analyze the motivic Adams spectral sequence. After incorporating the $\rho$-Bockstein differential $d_1(\tau) = \rho h_0$ \cite[Section 7.2]{WO17}, the computation of the $E_2$-term of the Adams spectral sequence is similar to the $q \equiv 1 \pmod 4$ case. As in that case, the spectral sequence collapses at the $E_2$-page. The result can now be read off from the Adams $E_\infty$-page.
\begin{thm}
The $E_\infty$-page of the motivic Adams spectral sequence converging to the very effective Hermitian \(K\)-theory of \(BC_2\) over \(\f_q\) with \(q \equiv 3 \pmod 4\) , as a module over the \(E_\infty\)-page for \(kq_{**}\), is given by
\[
E_\infty \cong \Sigma^{1,1}
\frac{\f_2[\tau,\rho,h_1,b]}
     {(\tau h_1^3,\rho^2)} \oplus 
\bigoplus_{i\ge 0}
\Sigma^{4i+3,2i+2}
\frac{\f_2[\tau^2,\rho,h_0]}
     {(h_0^{2i+2},\rho^2,\rho h_0)}\{1,\tau \rho \}.
\]
\end{thm}
\begin{remark}
    Again, the combinations of hollow diamonds and hollow circles represent hidden $\tau$-extensions in the Atiyah--Hirzebruch spectral sequence.
\end{remark}

\begin{remark}
Our computations over $\f_q$ can be adapted with minimal effort to $\q_q$ for $q$ odd. Up to a slight modification in the degrees of $2$-torsion, the coefficients $H\z_{**}^{\q_q}$, $kgl_{**}^{\q_q}$, and $kq_{**}^{\q_q}$ can be obtained from $H\z_{**}^{\f_q}$, $kgl_{**}^{\f_q}$, and $kq_{**}^{\f_q}$ by tensoring with an exterior algebra $\f_2[\pi]/(\pi^2)$, where $|\pi| = (-1,-1)$. See \cite[Sec. 4.3]{KQ24} for details. This ``tensoring up" can be extended through all of our computations to produce $E_{**}^{\q_q}(BC_2)$ for $E = H\z, \ kgl, \ kq$. 
\end{remark}

\bibliographystyle{alpha}
\bibliography{master}

\end{document}